\documentclass[10pt,letterpaper]{article}

\usepackage[centertags]{amsmath}
\usepackage{amsthm}
\usepackage{amssymb}
\usepackage{amscd}
\usepackage{eucal}
\usepackage{epsfig}
\usepackage[matrix,arrow,curve]{xy}
\CompileMatrices

\title{A long exact sequence for\\ symplectic Floer cohomology}
\author{Paul Seidel}
\date{Revised version, May 2002}

\begin{document}
\maketitle

\newcommand{\hatx}{\hat{x}}
\newcommand{\Sympe}{\mathrm{Symp}^e}
\newcommand{\crit}{{crit}}
\newcommand{\hatE}{\widehat{E}}
\newcommand{\hatO}{\widehat{\Omega}}
\newcommand{\hatQ}{\widehat{Q}}
\newcommand{\hatJ}{\widehat{J}}
\newcommand{\JJ}{\mathcal{J}}
\newcommand{\MM}{\mathcal{M}}
\newcommand{\WW}{\mathcal{W}}
\newcommand{\JJreg}{\JJ^{reg}}
\newcommand{\TT}{\mathcal{T}}
\newcommand{\BB}{\mathcal{B}}
\newcommand{\EE}{\mathcal{E}}
\newcommand{\FF}{\mathcal{F}}
\newcommand{\GG}{\mathcal{G}}
\newcommand{\loc}{{loc}}
\newcommand{\PP}{\mathcal{P}}
\newcommand{\gen}[1]{\leftsc #1 \rightsc}
\newcommand{\dist}{{dist}}


\newcommand{\R}{\mathbb{R}}
\newcommand{\Z}{\mathbb{Z}}
\newcommand{\Q}{\mathbb{Q}}
\newcommand{\C}{\mathbb{C}}
\newcommand{\N}{\mathbb{N}}
\newcommand{\half}{{\textstyle\frac{1}{2}}}
\newcommand{\quarter}{{\textstyle\frac{1}{4}}}

\newcommand{\iso}{\cong}           
\newcommand{\htp}{\simeq}          
\newcommand{\smooth}{C^\infty}
\newcommand{\CP}[1]{\C {\mathrm P}^{#1}}
\newcommand{\RP}[1]{\R {\mathrm P}^{#1}}
\newcommand{\leftsc}{\langle}
\newcommand{\rightsc}{\rangle}
\newcommand{\suchthat}{\; : \;}

\newcommand{\id}{\mathrm{id}}
\newcommand{\ind}{\mathrm{ind}}
\newcommand{\re}{\mathrm{re}}
\newcommand{\im}{\mathrm{im}}
\renewcommand{\ker}{\mathrm{ker}}
\newcommand{\coker}{\mathrm{coker}}
\newcommand{\mymod}{\quad\text{mod }}
\newcommand{\Hom}{\mathrm{Hom}}
\newcommand{\End}{\mathrm{End}}

\newcommand{\mo}{(M,\omega)}
\renewcommand{\o}{\omega}
\renewcommand{\O}{\Omega}
\newcommand{\Diff}{\mathrm{Diff}}


\theoremstyle{plain}
\newtheorem{itheorem}{Theorem}
\newtheorem{thm}{Theorem}[section]
\newtheorem{theorem}[thm]{Theorem}
\newtheorem{cor}[thm]{Corollary}
\newtheorem{corollary}[thm]{Corollary}
\newtheorem{lemma}[thm]{Lemma}
\newtheorem{prop}[thm]{Proposition}
\newtheorem{proposition}[thm]{Proposition}
\newtheorem{defn}[thm]{Definition}
\newtheorem{definition}[thm]{Definition}
\newtheorem{definitions}[thm]{Definitions}
\newtheorem{rem}[thm]{Remark}
\newtheorem{remarks}[thm]{Remarks}
\newtheorem{remark}[thm]{Remark}
\newtheorem{example}[thm]{Example}
\newtheorem*{acknow}{Acknowledgments}


\newenvironment{condensedlist}%
{\renewcommand{\theenumi}{(\roman{enumi})}
\renewcommand{\labelenumi}{\theenumi}
\parskip0em
\begin{enumerate} \parsep0em
\parskip0em
\itemsep0em}{\vspace{-0.5em}\end{enumerate}}

\newenvironment{condensedprimelist}%
{\renewcommand{\theenumi}{(\roman{enumi}')}
\renewcommand{\labelenumi}{\theenumi}
\parskip0em
\begin{enumerate} \parsep0em
\parskip0em
\itemsep0em}{\vspace{-0.5em}\end{enumerate}}

\newenvironment{romanlist}%
{\renewcommand{\theenumi}{(\roman{enumi})}
\renewcommand{\labelenumi}{\theenumi}
\begin{enumerate} \parsep0em
\itemsep0em \parskip0.5em}{\end{enumerate}}

\newenvironment{Romanlist}%
{\renewcommand{\theenumi}{(\Roman{enumi})}
\renewcommand{\labelenumi}{\theenumi}
\begin{enumerate} \parsep0em
\itemsep0em \parskip0.5em}{\end{enumerate}}

\newenvironment{theoremlist}%
{\begin{list}{{\rm(\roman{enumi}) }}{\usecounter{enumi}
\renewcommand{\theenumi}{(\alph{enumi})}
\renewcommand{\labelenumi}{\theenumi}
\leftmargin0cm \labelsep0cm \rightmargin0cm \parsep1em \listparindent0em
\itemsep0em \topsep0em \parskip0.5em \setlength{\labelwidth}{\fill}}}
{\parskip0em \end{list}}


\newcommand{\printwarning}[1]{%
\typeout{ #1 }%
}

\newcommand{\FIGUREGAP}[1]{%
\printwarning{Figure #1 missing}%
\begin{figure}[htb] \begin{center}
\setlength{\unitlength}{1cm} \framebox[10cm]{\begin{picture}(6,4)
\end{picture}}
\caption{#1%
\label{fig:#1}}
\end{center} \end{figure}}

\newcommand{\includefigure}[3]{%
\begin{figure}[#3]
\begin{center}
\epsfig{file=#2} \\ \caption{\label{fig:#1}}
\end{center}
\end{figure}}

\section*{Introduction}

Let $(M^{2n},\o,\alpha)$ be a compact symplectic manifold with contact type
boundary: $\alpha$ is a contact one-form on $\partial M$ which satisfies
$d\alpha = \o|\partial M$ and makes $\partial M$ convex. Assume in addition
that $[\o,\alpha] \in H^2(M,\partial M;\R)$ is zero, so that $\alpha$ can be
extended to a one-form $\theta$ on $M$ satisfying $d\theta = \o$. After fixing
such a $\theta$ once and for all, one can talk about exact Lagrangian
submanifolds in $M$. The Floer cohomology of two such submanifolds is
comparatively easy to define, since the corresponding action functional has no
periods, so that bubbling is impossible. The aim of this paper is to prove the
following result, which was announced in \cite{seidel00} (with an additional
assumption on $c_1(M)$ that has been removed in the meantime).

\begin{itheorem} \label{th:main}
Let $L$ be an exact Lagrangian sphere in $M$ together with a preferred
diffeomorphism $f: S^n \rightarrow L$. One can associate to it an exact
symplectic automorphism of $M$, the Dehn twist $\tau_L = \tau_{(L,[f])}$. For
any two exact Lagrangian submanifolds $L_0,L_1 \subset M$, there is a long
exact sequence of Floer cohomology groups
\begin{equation} \label{eq:exact} \xymatrix@C=-4em{
 {HF(\tau_L(L_0),L_1)} \ar[rr] && {HF(L_0,L_1)} \ar[dl] \\
 & {HF(L,L_1) \otimes HF(L_0,L).} \ar[ul]
}
\end{equation}
\end{itheorem}

The original inspiration for this came from the exact sequence in
Donaldson-Floer theory \cite{braam-donaldson94}, which can be translated into
symplectic geometry using various versions, proved \cite{dostoglou-salamon94}
and unproved, of the Atiyah-Floer conjecture. That line of thought should have
a Seiberg-Witten sibling, starting from \cite{carey-marcolli-wang98}, but the
corresponding Atiyah-Floer type relationships are only just beginning to be
understood \cite{salamon00}, \cite{ozsvath-szabo01}. In any case, the exact
sequences obtained from such speculations differ somewhat in generality from
that stated above. This reflects the fact that the first motivation has been
largely superseded by different ones, coming from mirror symmetry. Kontsevich's
homological mirror conjecture \cite{kontsevich94} for Calabi-Yau varieties
implies a relation between symplectic automorphisms and self-equivalences of
derived categories of coherent sheaves; see the survey \cite{aspinwall01} or
the papers \cite{seidel-thomas99}, \cite{horja01}. A particular class of
self-equivalences, ``twist functors along spherical objects'', is expected to
correspond to Dehn twists. By definition, twist functors give rise to an exact
sequence of the same form as \eqref{eq:exact}, with Floer cohomology replaced
by Ext-groups, so that the expected correspondence fits in well with our
result. With respect to the whole of Kontsevich's conjecture, this is a rather
peripheral issue. To see the exact sequence take on a more important role, one
has to pass to a related context, namely mirror symmetry for Fano varieties.
The derived categories of coherent sheaves on such varieties are often
generated by exceptional collections, which are subject to transformations
called mutations \cite{rudakov90}. The mirror dual notion is that of
distinguished basis of vanishing cycles, which is well-known in
Picard-Lefschetz theory. A rigorous connection between the two concepts is
established by \cite[Theorem 3.3]{seidel00} whose proof relies strongly on
Theorem \ref{th:main}; at present, this would seem to be its main application.
We conclude our discussion with some more concrete remarks about the statement
of the theorem:

(i) In this paper, Floer cohomology groups are treated as ungraded groups. One
can of course assume that $L,L_0,L_1$ are oriented, and then the groups become
$\Z/2$-graded. Different conventions are in use, but if one adopts that in
which the Euler characteristic of Floer cohomology is $(-1)^{n(n+1)/2}$ times
the intersection number, the degrees mod two of the maps in \eqref{eq:exact}
are
\begin{equation} \label{eq:grading}
\xymatrix{
 {\bullet} \ar[rr]^{0} && {\bullet} \ar[dl]^{n} \\
 & {\bullet} \ar[ul]^{1-n}
}
\end{equation}
Under more restrictive assumptions, one can introduce $\Z$-gradings. There are
several roughly equivalent ways of doing this; we adopt the approach of
\cite{kontsevich94} and \cite{seidel99}, in which one fixes a trivialization of
the bicanonical bundle $K^2_M$, thus establishing a notion of ``graded
Lagrangian submanifold''. Suppose that preferred gradings have been chosen for
$L,L_0,L_1$. Through the grading of $\tau_L$ itself, this induces a grading of
$\tau_L(L_0)$. All Floer cohomology groups in the exact sequence are then
canonically $\Z$-graded, and the degrees of the maps are as in
\eqref{eq:grading}. This is not difficult to show, it just requires a few
Maslov index computations.

(ii) We use Floer cohomology with $\Z/2$-coefficients. Inspection of the
discussion of coherent orientation in \cite{fukaya-oh-ohta-ono} suggests that
at least when $L_0,L_1$ are spin, the exact sequence should exist with
$\Z$-coefficients (with the $\otimes$ replaced by the cohomology of the
underlying tensor product of cochain groups, to avoid K{\"u}nneth terms).
However, I have not checked all the details.

(iii) The assumption $[\o,\alpha] = 0$ is the main limitation of the theorem as
it stands. We have adopted this ``exact'' framework with a view to the
application in \cite{seidel00}, and also because it simplifies a number of
technical issues, thereby hopefully allowing the basic ideas to stand out.
Almost the same proof goes through in a few other cases, such as when suitable
``monotonicity'' conditions hold. On the other hand, a considerable amount of
work remains to be done to extend the exact sequence to the most general
situation where one would want to have it; a version for closed manifolds with
$c_1(M) = 0$ seems particularly desirable.

(iv) The map $\nwarrow$ in \eqref{eq:exact} is obtained by composing the
canonical isomorphism $HF(L_0,L) \iso HF(\tau_L(L_0),L)$ which exists because
$\tau_L(L) = L$, with a pair-of-pants product (a.k.a. Donaldson product). In
the spirit of \cite{piunikhin-salamon-schwarz94}, this can be seen as a kind of
relative Gromov invariant. A more general version of the same formalism,
involving pseudo-holomorphic sections of fibrations with singularities, yields
the second map $\rightarrow$. In contrast, $\swarrow$ appears as connecting map
in our construction, and is therefore defined only indirectly. A symmetry
consideration using the duality $HF(L_0,L_1) \iso HF(L_1,L_0)^\vee$ suggests
that this map should actually be a pair-of-pants coproduct. That is in fact
true, but we will not prove it here.

The exposition in the body of the paper follows a slightly indirect course, in
that we try to familiarize the reader with each ingredient separately, before
they all get mixed up into the main argument. There is even a small amount of
material which is not necessary for our immediate purpose, but which is closely
related and useful for further development. Thus, the whole first chapter is
elementary symplectic geometry, concentrating on topics related to
Picard-Lefschetz theory; the second one deals with pseudo-holomorphic curves,
which means setting up the relative invariants mentioned above, and introducing
certain techniques for partially computing them based on symplectic curvature;
and the only the third chapter addresses the actual proof.

\begin{acknow}
Much of this work was originally done in 1996. It goes without saying that I am
heavily indebted to Simon Donaldson, who was my advisor at that time.
Conversations with Michael Callahan, Mikhail Khovanov, Dietmar Salamon, and
Richard Thomas have been helpful. A preliminary version of this paper was
presented in a series of talks at Ecole Polytechnique in 1999; I thank my
colleagues there for their patience and comments. The final touches were put on
during a visit to the University of Michigan, which provided a hospitable
environment.
\end{acknow}

\newpage
\numberwithin{equation}{section}

\section[Dehn twists]{Dehn twists, and all that\label{ch:one}}

This chapter takes a look at basic Picard-Lefschetz theory from the symplectic
viewpoint. It has been known since Arnold's note \cite{arnold95} that such a
viewpoint makes sense, and it has been used for various purposes, see e.g.\
\cite{seidel98b}, \cite{khovanov-seidel98}. Still, the present paper seems to
be the most systematic attempt at an exposition so far. The reader may be
surprised by the exactness assumptions built into our framework. As far as the
elementary theory is concerned, there is no need to make such assumptions.
However, they greatly simplify the pseudo-holomorphic theory to be introduced
later on, and in order to keep the setup coherent, we have chosen to impose
them from the start.

\subsection{Exact symplectic geometry and fibrations\label{sec:basic}}

By an exact symplectic manifold we mean a compact manifold $M$ with boundary,
together with a symplectic form $\o$ and a one-form $\theta$ satisfying
$d\theta = \o$, such that $\theta|\partial M$ is a contact one-form and makes
$\partial M$ convex. An isomorphism of exact symplectic manifolds is a
diffeomorphism $\phi: M \rightarrow M'$ which is symplectic, satisfies
$\phi^*\theta' = \theta$ in some neighbourhood of $\partial M$, and such that
$[\phi^*\theta' - \theta] \in H^1(M,\partial M;\R)$ is zero. This means that
there is a unique function $K_\phi \in \smooth_c(M \setminus \partial M,\R)$
such that $\phi^*\theta' = \theta + dK_\phi$. We denote by $\Sympe(M)$ the
group of those exact symplectic automorphisms of $M$ which are the identity
near $\partial M$. Its Lie algebra consists of vector fields $X$ such that
$\o(\cdot,X) = dH$ for some $H$ which vanishes near $\partial M$, and is thus
identified with $\smooth_c(M \setminus \partial M,\R)$.

An exact Lagrangian submanifold in $M$ is a pair consisting of a Lagrangian
submanifold $L \subset M$ (always assumed to be disjoint from $\partial M$) and
a function $K_L$ on it such that $dK_L = \theta|L$. The image of $L$ under an
isomorphism $\phi: M \rightarrow M'$ of exact symplectic manifolds is again an
exact Lagrangian submanifold, in a canonical way; the associated function is
\begin{equation} \label{eq:induced-k}
K_{\phi(L)} = (K_L + K_\phi|L) \circ \phi^{-1}.
\end{equation}
In particular, $\Sympe(M)$ acts on the set of exact Lagrangian submanifolds. A
special situation which will occur later on is that $\phi \in \Sympe(M)$
satisfies $\phi(L) = L$ in the ordinary sense. Then (supposing $L$ to be
connected) $K_{\phi(L)} = K_L + c$ for some constant $c = K_\phi|L$, which
needs not be zero. This means that $\phi$ may not map $L$ to itself as an exact
Lagrangian submanifold, instead ``shifting'' it by some amount.

A notion of fibre bundle suitable for exact symplectic geometry is as follows.
Let $S$ be a smooth connected manifold, possibly with boundary (one could also
allow corners), and $\pi: E \rightarrow S$ a differentiable fibre bundle whose
fibres are compact manifolds with boundary. Write $\partial_hE \subset E$ for
the union of the boundaries of all the fibres. If the boundary of $S$ is empty,
$\partial_hE = \partial E$; otherwise $\partial E$ has another face
$\partial_vE = \pi^{-1}(\partial S)$, and the two faces meet at a codimension
two corner. An exact symplectic fibration is such an $(E,\pi)$ equipped with
$\Omega \in \Omega^2(E)$ and $\Theta \in \Omega^1(E)$, satisfying $d\Theta =
\Omega$, such that each fibre $E_z$ with $\o_z = \Omega|E_z$ and $\theta_z =
\Theta|E_z$ is an exact symplectic manifold. There is an additional condition
of triviality near $\partial_hE$, by which we mean the following: choose some
$z \in S$ and consider the trivial fibration $\tilde{\pi} : \tilde{E} = S
\times E_z \rightarrow S$, with the forms $\tilde{\Omega},\tilde{\Theta}$ which
are pullbacks of $\o_z,\theta_z$. There should be a fibrewise diffeomorphism
\begin{equation} \label{eq:boundary-trivialization}
\xymatrix{
 {N} \ar[rr] \ar[dr]_{\pi} && {\tilde{N}} \ar[dl]^{\tilde{\pi}} \\
 & {S} &
}
\end{equation}
betweeen neighbourhoods $N \subset E$, $\tilde{N} \subset \tilde{E}$ of
$\partial_hE$ resp.\ $\partial_h\tilde{E}$, which maps $\partial_hE$ to
$\partial_h\tilde{E}$, equals the identity on the fibre over $z$, and sends
$\Omega,\Theta$ to $\tilde{\Omega},\tilde{\Theta}$. Note that the choice of $z$
and the diffeomorphism are not considered to be part of the data defining an
exact symplectic fibration; only their existence is assumed.

\begin{lemma} \label{th:fibre-bundles}
Take a point $z \in S$ and a chart $\psi: U \rightarrow S$, with $U \subset
\R^k$ a contractible neighbourhood of $0$, such that $\psi(0) = z$. Then there
is a trivialization $\Psi: U \times E_z \rightarrow \psi^*E$, such that
$\Psi\,|\,\{0\} \times E_z = \id$ and
\[
\Psi^*\Theta = \theta_z + \sum_{i=1}^k H_i dt_i + dR.
\]
Here $H_1,\dots,H_k,R \in \smooth(U \times E_z,\R)$ are functions which vanish
near $U \times \partial E_z$, and $t_i$ the coordinates on $\R^k$. Moreover,
the difference between any two such trivializations $\Psi_1,\Psi_2$ is a map
$\Psi_2^{-1} \circ \Psi_1: (U,0) \rightarrow (\Sympe(E_z),\id)$. \qed
\end{lemma}

The proof is by a standard argument involving Moser's Lemma. The result means
first of all that any exact symplectic fibration $(E,\pi)$ has a structure of
$\Sympe(E_z)$-fibre bundle, where $E_z$ is any fibre. In addition to that, $E$
carries a preferred connection, which gives rise to canonical parallel
transport maps $\rho_c: E_{c(a)} \rightarrow E_{c(b)}$ over smooth paths $c:
[a;b] \rightarrow S$; these are exact symplectic isomorphisms, and lie in
$\Sympe(E_{c(a)})$ if $c$ is closed. To define the preferred connection, one
can use a local trivialization as before and set $A = \sum_i H_i\,dt_i$, which
is a one-form on $U$ with values in $\smooth_c(E_z \setminus \partial E_z,\R)$.
It is easy to check that this transforms in the proper way. A more intrinsic
approach is to observe that $TE_x = TE^h_x \oplus TE^v_x$ splits into a
horizontal and a vertical piece, given by $TE^v_x = \ker(D\pi_x)$ and
\begin{equation} \label{eq:connection}
TE^h_x = \{X \in TE_x \suchthat \Omega(X,\cdot)|TE^v_x = 0\}.
\end{equation}
Each $Z \in TS_z$ has a unique lift $Z^h \in \smooth(TE^h|E_z)$, and these
vectors define the same connection as before (conversely, one can show that
given a $\Sympe(M)$-fibre bundle and a compatible connection, one can equip its
total space with the structure of an exact symplectic fibration). In a local
trivialization as in Lemma \ref{th:fibre-bundles}, the curvature of the
connection is
\begin{equation} \label{eq:curvature}
F_A = \sum_{i<j} \left( -\frac{\partial H_i}{\partial t_j} + \frac{\partial
H_j}{\partial t_i} - \o_z(X_i,X_j) \right) dt_i \wedge dt_j,
\end{equation}
where $X_i$ is the family of Hamiltonian vector fields on $E_z$ corresponding
to $H_i(t,\cdot)$. If one takes the more intrinsic view, the curvature is a
two-form on $S$ with values in functions on the fibres, and is given by
$(Z_1,Z_2) \mapsto \Omega(Z_1^h,Z_2^h)$. In the case where the base $S$ is an
oriented surface, we say that $(E,\pi)$ is nonnegatively curved if for any
oriented chart $\psi$ and trivialization $\Psi$, the function in front of $dt_1
\wedge dt_2$ in \eqref{eq:curvature} is nonnegative; or equivalently, if
$\Omega|TE^h$ is nonnegative for the induced orientation of $TE^h$. To see what
this means, consider an exact symplectic fibration $(E,\pi)$ over the closed
unit disc $\bar{D}(1) \subset \C$, and the monodromy $\rho = \rho_{\partial
\bar{D}(1)} : E_1 \rightarrow E_1$ around the boundary in positive sense. If
the curvature is nonnegative, one can write $\rho$ as time-one map of some
(time-dependent) Hamiltonian on $E_1$ which vanishes near $\partial E_1$ and is
$\leq 0$ everywhere.

\begin{example} \label{ex:mapping-torus}
Let $E^\rho = \R \times M / (t,x) \sim (t-1,\rho(x))$ be the mapping
torus of $\rho \in \Sympe(M)$. To make this into an exact symplectic
fibration over $S^1 = \R/\Z$, one chooses a function $R_\rho \in
\smooth(\R \times M,\R)$ such that $R_\rho(t-1,\rho(x)) = R_\rho(t,x)
- K_\rho(x)$, and then sets $\Omega_{E^\rho} = \o$, $\Theta_{E^\rho}
= \theta + dR_\rho$. The preferred connection is the one induced from
the trivial connection on $\R \times M$, and its monodromy is $\rho$
itself. It is easy to prove that all exact symplectic fibrations over
a circle are of this form.
\end{example}

For future use, we make an observation on the compatibility of symplectic
parallel transport with Lagrangian submanifolds. Let $c: [a;b] \rightarrow S$
be a smooth embedded path. Suppose that we have an exact Lagrangian submanifold
in each fibre $E_{c(t)}$, depending smoothly on $t$; by this we mean a
subbundle $Q \subset E|\im(c)$ such that each fibre $Q_{c(t)} \subset E_{c(t)}$
is Lagrangian, together with a $K_Q \in \smooth(Q,\R)$ whose restrictions
$K_{Q_{c(t)}} = K_Q|Q_{c(t)}$ make the $Q_{c(t)}$ exact.

\begin{lemma} \label{th:submanifold}
Assume that all the $Q_{c(t)}$ are connected. Then the following conditions are
equivalent:
\begin{condensedlist}
\item
$\Omega|Q = 0$;
\item
$\Theta|Q = dK_Q + \pi^*\kappa_Q$ for some $\kappa_Q \in \Omega^1(\im(c))$;
\item \label{item:carry-along}
the maps $\rho_{c|[s;s']}: E_{c(s)} \rightarrow E_{c(s')}$ satisfy
$\rho_{c|[s;s']}(Q_{c(s)}) = Q_{c(s')}$ for all $s,s'$. \qed
\end{condensedlist}
\end{lemma}

The proof is straightforward. To be precise, \ref{item:carry-along} concerns
$Q_{c(s)}$ as Lagrangian submanifolds only. Taking the functions into account
and using \eqref{eq:induced-k} yields
\[
K_{\rho_{c|[s;s']}(Q_{c(s)})} = K_{Q_{c(s')}} + \int_{c|[s;s']} \kappa_Q.
\]
Hence $\rho_{c|[s;s']}(Q_{c(s)}) = Q_{c(s')}$ holds in the sense of exact
Lagrangian submanifolds iff $\kappa_Q = 0$, or what is the same, $\Theta|Q =
dK_Q$.

Just like any kind of fibre bundle with connection, exact symplectic fibrations
can be manipulated by cut-and-paste methods. As an example, take two oriented
surfaces $S^k$, $k = 1,2$, and boundary circles $C^k \subset \partial S^k$; and
let $S$ be the surface obtained by identifying $C^1,C^2$
orientation-reversingly. Suppose that we have exact symplectic fibrations
$(E^k,\pi^k,\Omega^k,\Theta^k)$ over $S^k$ whose monodromies around $C^k$,
taken in opposite senses, coincide. Then, assuming additionally that
$(E^k,\pi^k)$ is flat (has zero curvature) near $C^k$, one can construct from
them an exact symplectic fibration $(E,\pi)$ over $S$. It is maybe useful to
give some details of this. To start, take oriented collars $\psi^1: (-1;0]
\times \R/\Z \rightarrow S^1$, $\psi^2: [0;1) \times \R/\Z \rightarrow S^2$
around $C^1,C^2$ respectively, so that $S$ can be defined by using $\psi^2
\circ (\psi^1)^{-1}$ to identify the two circles. Our main assumption is that
there should be an exact symplectic isomorphism between the fibres of $E^k$
over $\psi^k(0,0)$, say, which relates the monodromies around $\psi^k(\{0\}
\times \R/\Z)$. Because of flatness, an equivalent formulation is that there is
some mapping torus $E^\rho$ as in Example \ref{ex:mapping-torus} and
diffeomorphisms
\[
\xymatrix{
 {(-1;0] \times E^\rho} \ar[d] \ar[r]^-{\Psi^1} & {E^1} \ar[d]^{\pi^1} \\
 {(-1;0] \times \R/\Z} \ar[r]^-{\psi^1} & {S^1}
} \qquad \xymatrix{
 {[0;1) \times E^\rho} \ar[d] \ar[r]^-{\Psi^2} & {E^2} \ar[d]^{\pi^2} \\
 {[0;1) \times \R/\Z} \ar[r]^-{\psi^2} & {S^2}
}
\]
such that $(\Psi^k)^*\Omega^k = \Omega_{E^\rho}$ and $(\Psi^k)^*\Theta^k =
\Theta_{E^\rho} + dR^k$ for some functions $R^k$. Clearly, one can introduce
modified forms $\tilde{\Theta}^k = \Theta^k - d\tilde{R}^k$ with suitable
functions $\tilde{R}^k$, in such a way that $(\Psi^k)^*\tilde{\Theta}^k =
\Theta_{E^\rho}$ near $\{0\} \times E^\rho$. Gluing together the $(E^k,\pi^k)$
along $C^k$ via $\Psi^2 \circ (\Psi^1)^{-1}$ yields a smooth fibration
$(E,\pi)$ over $S$, and the $\Omega^k, \tilde{\Theta}^k$ match up to forms
$\Omega,\Theta$ on it, making it an exact symplectic fibration. Since we have
not changed the symplectic connection, nonnegativity of the curvature of
$(E^k,\pi^k)$ implies the same for $(E,\pi)$.

\begin{remark} \label{re:squeeze-boundary}
The flatness condition on $(E^k,\pi^k)$ near $C^k$ can be removed. Namely,
suppose that it is not satisfied for $k = 1$. What one does then is to choose a
function $g \in \smooth([-1;0],\R)$ such that $g(s) = s$ for $s$ close to $-1$,
$g(s) = 0$ for $s$ close to $0$, and $g'(s) \geq 0$ everywhere; and define a
self-map $p$ of the surface $S^1$ by $p(z) = z$ for $z \notin \im(\psi^1)$,
$p(\psi^1(s,t)) = \psi^1(g(s),t)$. This collapses a small neighbourhood of
$C^1$ onto that boundary circle, so if we replace $(E^1,\pi^1)$ by its pullback
under $p$, it becomes flat near $C^1$, and the previous construction goes
through. It is noteworthy that this still preserves nonnegative curvature,
because $Dp$ has determinant $\geq 0$ everywhere.
\end{remark}

The basic objects of Picard-Lefschetz theory are fibrations over surfaces,
where the fibres are allowed to have certain particularly simple singularities.
Let $S$ be a connected oriented surface, possibly with boundary. An {\em exact
Lefschetz fibration}\footnote{Called ``exact Morse fibration'' in
\cite{seidel00}. The present terminology is more in line with general usage,
since the notion is closely related to, even though not quite the same as, the
symplectic Lefschetz fibrations considered by Donaldson, Gompf, and others.}
over $S$ consists of data $(E,\pi,\Omega,\Theta,J_0,j_0)$ as follows. $E$ is a
$(2n+2)$-dimensional manifold whose boundary is the union of two faces
$\partial_hE$ and $\partial_vE$, meeting at a codimension two corner. $\pi: E
\rightarrow S$ is a proper map with $\partial_vE = \pi^{-1}(\partial S)$ (so
this may be empty), and such that both $\pi|\partial_hE: \partial_hE
\rightarrow S$ and $\pi|\partial_vE: \partial_vE \rightarrow \partial S$ are
smooth fibre bundles. $\pi$ can have at most finitely many critical points, and
no two may lie on the same fibre (moreover, because of the previous
assumptions, they must lie in the interior of $E$). Denote by $E^\crit \subset
E$, $S^\crit \subset S$ the set of critical points resp.\ of critical values.
$J_0$ is an almost complex structure on a neighbourhood of $E^\crit$, and $j_0$
a positively oriented complex structure on a neighbourhood of $S^\crit$. These
are such that $\pi$ is $(J_0,j_0)$-holomorphic near $E^\crit$, and the Hessian
$D^2\pi$ at any critical point is nondegenerate as a complex quadratic form.
The closed two-form $\Omega \in \Omega^2(E)$ must be nondegenerate on $TE_x^v =
\ker(D\pi_x)$ for each $x \in E$, and a K{\"a}hler form for $J_0$ in some
neighbourhood of $E^\crit$. $\Theta \in \Omega^1(E)$ must satisfy $d\Theta =
\Omega$. We also require triviality near $\partial_hE$, which means the
existence of a map \eqref{eq:boundary-trivialization} with the same properties
as for exact symplectic fibrations. For brevity, exact Lefschetz fibrations
will usually be denoted by $(E,\pi)$ alone, as we have already done for exact
symplectic fibrations.

For any $x \in E$ there is a decomposition $TE_x = TE^h_x \oplus TE^v_x$ with
$TE^h_x$ defined as in \eqref{eq:connection}; the horizontal part is zero at
critical points, and projects isomorphically to $TS_z$, $z = \pi(x)$, at any
other point. We say that an exact Lefschetz fibration has nonnegative curvature
if $\Omega|TE^h_x \geq 0$ for each $x$. Note that if $x \notin E^\crit$ is
close to a critical point, $\Omega|TE^h_x$ is strictly positive anyway, because
of the K{\"a}hlerness assumption on $\Omega$. A standard argument based on this
shows

\begin{lemma} \label{th:add-base}
If $\beta \in \Omega^2(S)$ is a sufficiently positive two-form, $\Omega +
\pi^*\beta$ is a symplectic form on $S$. \qed
\end{lemma}

Symplectic parallel transport for an exact Lefschetz fibration is well-defined
as long as one avoids the critical fibres; indeed, if one removes those fibres,
the remainder is an exact symplectic fibration over $S \setminus S^\crit$. We
now take a look at the structure of the critical points. Take $z_0 \in S^\crit$
and local $j_0$-holomorphic coordinates $\xi: U \rightarrow S$, where $U
\subset \C$ is a neighbourhood of the origin, such that $\xi(0) = z_0$. By
assumption there is a unique critical point $x_0 \in E_{z_0}$. The holomorphic
Morse Lemma says that one can find a neighbourhood of the origin $W \subset
\C^{n+1}$ and a $J_0$-holomorphic chart $\Xi: W \rightarrow E$ with $\Xi(0) =
x_0$, such that
\[
(\xi^{-1} \circ \pi \circ \Xi)(x) = x_1^2 + \dots + x_{n+1}^2
\]
is the standard nondegenerate quadratic form on $\C^{n+1}$. We call $(\xi,\Xi)$
a holomorphic Morse chart. In general, it is not possible to choose $\Xi$ in
such a way that $\Xi^*\Omega$ is the standard K{\"a}hler form on $W \subset
\C^{n+1}$; however, one can remedy this by a suitable local deformation.

\begin{lemma} \label{th:local-deformation}
Let $(E,\pi,\Omega,\Theta,J_0,j_0)$ be an exact Lefschetz fibration, and $x_0$
a critical point of $\pi$. Then there are smooth families $\Omega^\mu \in
\Omega^2(E)$, $\Theta^\mu \in \Omega^1(E)$, $0 \leq \mu \leq 1$, such that
\begin{condensedlist}
\item
$\Omega^0 = \Omega$, $\Theta^0 = \Theta$;
\item
for all $\mu$, $\Omega^\mu = \Omega^0$ and $\Theta^\mu = \Theta^0$ outside a
small neighbourhood of $x_0$;
\item
each $(E,\pi,\Omega^\mu,\Theta^\mu,J_0,j_0)$ is an exact Lefschetz fibration;
\item
there is a holomorphic Morse chart $(\xi,\Xi)$ around $x_0$ such that
$\Xi^*\Omega^1$, $\Xi^*\Theta^1$ agree near the origin with the standard forms
$\o_{\C^{n+1}} = \frac{i}{2} \sum dx_k \wedge d\bar{x}_k$, $\theta_{\C^{n+1}} =
\frac{i}{4} \sum x_k d\bar{x}_k - \bar{x}_k dx_k$.
\end{condensedlist}
\end{lemma}

The proof is based on an elementary local statement about K{\"a}hler forms.

\begin{lemma} \label{th:kaehler}
Let $\o$ be a K{\"a}hler form on the ball $B = B^{2n+2}(r)$ of radius $r>0$ in
$\C^{n+1}$. Then there is another K{\"a}hler form $\o'$ which agrees with $\o$
near $\partial B$, and which close to the origin is some small multiple of
$\o_{\C^{n+1}}$.
\end{lemma}

\proof The first step is to find a K{\"a}hler form $\o''$ which is equal to
$\o$ near $\partial B$, and which has constant coefficients near the origin.
For this write\footnote{The definition of $d^c$ in this paper is such that
$\o_{\C} = -dd^c(\quarter|z|^2)$. This differs from the majority convention by
a negative constant.} $\o = \beta + dd^c f$, with $\beta$ constant and $f$
vanishing to second order at $x = 0$. Take a cutoff function $g \in
\smooth(\R^+,\R)$ such that $g(t) = 1$ for $t \leq 1$ and $g(t) = 0$ for $t
\geq 2$, and set $f_\epsilon(x) = g(||x||/\epsilon)f(x)$. A straightforward
computation shows that as $\epsilon \rightarrow 0$, the functions $f_\epsilon$
not only have increasingly small support, but also tend to $0$ in the $C^2$
topology. The desired form is, for small $\epsilon$,
\begin{equation} \label{eq:ohalf}
\o'' = \o - dd^c f_\epsilon = \beta + dd^c (f - f_\epsilon).
\end{equation}

In a second step, choose some small $t>0$ such that the constant form $\beta' =
\beta - t\o_{\C^{n+1}}$ is still K{\"a}hler. Take a two-form $\gamma$ on
$\C^{n+1}$ which is of type (1,1) and nonnegative everywhere, which vanishes
near the origin, and which is equal to $\o_{\C^{n+1}}$ outside a compact
subset; this can be obtained as $\gamma = -dd^c h(||x||^2)$ for a suitable
convex function $h$. Pulling $\gamma$ back by a linear map transforms it into
another nonnegative (1,1)-form $\gamma'$, zero near the origin and equal to
$\beta'$ outside a compact subset. Then $\gamma'' = \gamma' + t\o_{\C^{n+1}}$
is K{\"a}hler, equals $t\o_{\C^{n+1}}$ near the origin, and $\beta$ outside a
compact subset. That compact subset can be made arbitrarily small by retracting
linearly and rescaling; the two-form obtained in that way can be plugged into
$\o''$ locally near zero, yielding $\o'$. \qed

\proof[Proof of Lemma \ref{th:local-deformation}] Take some holomorphic Morse
chart $(\xi,\Xi)$ for $x_0$. Lemma \ref{th:kaehler} says that one can find a
two-form $\Omega^1$ on $E$ which agrees with $\Omega = \Omega^0$ outside
$\im(\Xi)$, such that $\Xi^*\Omega^1$ is K{\"a}hler and, near the origin,
equals $c\, \o_{\C^{n+1}}$ for some small constant $c>0$. The obstruction to
finding a $\Theta^1 \in \Omega^1(E)$ which agrees with $\Theta = \Theta^0$
outside $\im(\Xi)$ and satisfies $d\Theta^1 = \Omega^1$ lies in $H^2(B,\partial
B;\R)$, with $B$ a $(2n+2)$-dimensional ball; which is zero. An arbitrarily
chosen $\Theta^1$ needs to be modified to make $\Xi^*\Theta^1$ equal to
$c\,\theta_{\C^{n+1}}$ near $x = 0$, but that can be done by adding the
differential of some function to it. By restricting $\xi,\Xi$ to smaller
neighbourhoods, and rescaling them by $c^{-1/2}$ and $c^{-1}$ respectively, one
achieves that $\Xi^*\Omega^1 = \o_{\C^{n+1}}$ and $\Xi^*\Theta^1 =
\theta_{\C^{n+1}}$ near the origin. Finally, $\Omega^\mu$ and $\Theta^\mu$ are
defined by interpolating linearly between $\mu = 0$ and $1$; the required
properties are obvious. \qed

\subsection{The local model\label{sec:model}}

Consider $T = T^*\!S^n$ with its standard forms $\o_T \in \Omega^2(T)$,
$\theta_T \in \Omega^1(T)$. For concrete computations we use the coordinates $T
= \{ (u,v) \in \R^{n+1} \times \R^{n+1} \suchthat \leftsc u, v \rightsc = 0, \;
||v|| = 1 \}$. For each $\lambda>0$, the subspace $T(\lambda)$ of cotangent
vectors of length $\leq \lambda$ is an exact symplectic manifold. We write
similarly $T(0) \subset T$ for the zero-section. The length function $\mu: T
\rightarrow \R$, $\mu(u,v) = ||u||$, generates a Hamiltonian circle action
$\sigma$ on $T \setminus T(0)$ which, after identifying $T \iso TS^n$ via the
standard metric, can be described as the normalized geodesic flow on $S^n$. In
coordinates
\[
\sigma_t(u,v) = (\cos(t)u - \sin(t) ||u|| v,\cos(t)v + \sin(t)\frac{u}{||u||}).
\]
$\sigma_\pi$ is the antipodal involution $A(u,v) = (-u,-v)$, hence extends
continuously over $T(0)$ (unlike any $\sigma_t$, $0 < t < \pi$). This can be
used to define certain symplectic automorphisms of $T$. The construction is by
now well-known, but we repeat it here since precise control over the parameters
will be important later on.

\begin{lemma} \label{th:on-invariant}
Let $R \in \smooth(\R,\R)$ be a function which vanishes for $t \gg 0$ and which
satisfies $R(-t) = R(t) - kt$ for small $|t|$, with some $k \in \Z$. Let
$(\phi_t^H)$ be the Hamiltonian flow of $H = R(\mu)$, defined on $T \setminus
T(0)$. Then $\phi_{2\pi}^H$ extends smoothly over $T(0)$ to a compactly
supported symplectic automorphism $\phi$ of $T$. The function $K = 2\pi
\,(R'(\mu)\mu - R(\mu))$ also extends smoothly over $T(0)$, and
\begin{equation} \label{eq:k-function}
\phi^*\theta_T - \theta_T = dK.
\end{equation}
\end{lemma}

\proof Two Hamiltonians $H_1,H_2$ which are both functions of $\mu$ always
Poisson-commute, so that $\phi_t^{H_1} \phi_t^{H_2} = \phi_t^{H_1 + H_2}$.
Decomposing $H$ into $H_1 = R(\mu) - (k/2)\mu$ and $H_2 = (k/2)\mu$, one gets
\[
\phi_{2\pi}^H = \phi_{2\pi}^{H_1} \circ \sigma_{k\pi}.
\]
We may assume that $R(-t) = R(t) - kt$ holds everywhere, since that can be
achieved by modifying $R$ for negative values only, which does not affect
$\phi$. Then $R(t) - (k/2)t$ is an even function, so it can be written as a
smooth function of $t^2$. This proves that $H_1$ and its flow extend smoothly
over $T(0)$. We know that $\sigma_{k\pi} = A^k$ extends smoothly, so the same
holds for $\phi_{2\pi}^H$. Since $H(y)$ vanishes for points with $\mu(y) \gg
0$, $\phi$ is compactly supported (one can show that the compactly supported
symplectic automorphisms which are obtained in this way are precisely those
which are equivariant for the obvious $O(n+1)$-action). The function $R'(t)t -
R(t)$ is even, so $K$ extends smoothly over $T(0)$ for the same reason as
before. A computation shows that the Hamiltonian vector field $X$ of $H$
satisfies $L_{2\pi X}\theta_T = dK$. Since $H$ and $K$ are both functions of
$\mu$, $\phi_t^H$ preserves $K$, which implies \eqref{eq:k-function}. \qed

Clearly, if $supp(R) \subset (-\infty;\lambda)$ then $\phi$ restricts to a
symplectic automorphism of $T(\lambda)$ which is the identity near the
boundary. \eqref{eq:k-function} shows that this is an exact symplectic
automorphism. In the case $k = 1$ we call these automorphisms {\em model Dehn
twists}, and generally denote them by $\tau$; any two of them are isotopic in
$\Sympe(T(\lambda))$. An explicit formula is
\begin{equation} \label{eq:explicit-twist}
\tau(y) =
 \begin{cases}
 \sigma_{2\pi R'(\mu(y))}(y) & y \in T(\lambda) \setminus T(0),\\
 A(y) & y \in T(0),
 \end{cases}
\end{equation}
where the angle of rotation goes from $2\pi R'(0) = \pi$ to $2\pi R'(\lambda) =
0$. Note that $\tau$ maps $T(0)$ to itself, and is the antipodal map on it. If
one considers $T(0)$ as an exact Lagrangian submanifold, with a function
$K_{T(0)} = const.$ associated to it, then by \eqref{eq:induced-k} and
\eqref{eq:k-function}
\begin{equation} \label{eq:self-shift}
K_{\tau(T(0))} = (K_{T(0)} + K_{\tau_L}|T(0)) \circ \tau_L^{-1} = K_{T(0)} -
2\pi R(0).
\end{equation}
On occasion, it is useful to demand that the angle $R'(t)$ does not oscillate
too much. We say that $\tau$ is {\em $\delta$-wobbly} for some $0 < \delta <
1/2$ if $R'(t) \geq 0$ for all $t \geq 0$, and $R''(t) < 0$ for all $t \geq 0$
such that $R'(t) \geq \delta$.

\begin{lemma} \label{th:local-intersections}
Suppose that $\tau$ is $\delta$-wobbly. Let $F_0 = T(\lambda)_{y_0}$, $F_1 =
T(\lambda)_{y_1}$ be the fibres of $T(\lambda) \rightarrow S^n$ at points
$y_0,y_1$, whose distance in the standard metric is $\dist(y_0,y_1) \geq 2\pi
\delta$. Then $\tau(F_0)$ intersects $F_1$ transversally and at a single point
$y$, which satisfies
\[
2\pi R'(||y||) = \dist(y_0,y_1).
\]
In the special case where $y_1 = A(y_0)$ one has $y = y_1$; and then the
tangent space of $T(\lambda)$ at $y$ can be identified symplectically with
$\C^n$ in such a way that the subspaces tangent to $\tau(F_0)$, $T(0)$, $F_1$
become respectively, $\R^n$, $e^{2\pi i/3}\R^n$, and $e^{\pi i/3}\R^n$.
\end{lemma}

\proof Consider first the case when $y_1 \neq A(y_0)$. Suppose that $y \in F_1$
is a point with $\tau^{-1}(y) \in F_0$. By identifying $T \iso TS^n$ and using
the interpretation of $\sigma$ as normalized geodesic flow, one sees that $y$
must be a positive multiple of $c'(1)$, where $c: [0;1] \rightarrow S^n$ is the
minimal geodesic from $c(0) = y_0$ to $c(1) = y_1$. Moreover, the angle of
rotation must be $2\pi R'(||y||) = ||c'(1)|| = dist(y_0,y_1)$. These two
conditions are also sufficient. $\delta$-wobblyness implies that $2\pi R'(t) =
dist(y_0,y_1)$ has exactly one solution $t > 0$, which proves that there is
exactly one $y$. Combine the two conditions above into one, $c'(1) = 2\pi
R'(||y||)(y/||y||)$. Taking the derivative, one sees that a vector $Y \in
T(F_1)_y \iso T(S^n)_{y_1}$ satisfies $(T\tau)^{-1}(Y) \in T(F_0)$ iff
\[
 \textstyle R''(||y||) \leftsc \frac{y}{||y||}, Y \rightsc \frac{y}{||y||} +
 R'(||y||) \big( Y - \leftsc \frac{y}{||y||}, Y \rightsc \frac{y}{||y||} \big) = 0.
\]
We know that $R'(||y||) \geq \delta$; by $\delta$-wobblyness this implies
$R''(||y||) < 0$, which shows that $Y = 0$. Therefore $y \in \tau(F_0) \cap
F_1$ is a transverse intersection point.

Now consider the case when $y_1 = A(y_0)$. Then $y = y_1$ clearly lies in
$\tau(F_0) \cap F_1$, and because $R'(t) < 1/2$ for all $t>0$, there is no
other intersection point. In the notation from the proof of Lemma
\ref{th:on-invariant}, $\tau(F_0) = \phi^{H_1}_{2\pi}(A(F_0)) =
\phi^{H_1}_{2\pi}(F_1)$, and $y$ is a stationary point of $(\phi^{H_1}_t)$. It
follows that $T(\tau(F_0))_y$ is the image of $T(F_1)_y$ under the time $2\pi$
map of the linear Hamiltonian flow generated by the quadratic form $\half
Hess(H_1)_y$. If one identifies the tangent space to $T$ at $y$ with $T(S^n)_y
\oplus T(S^n)_y$ in such a way that the first summand is $T(F_1)_y$, then
\[
Hess(H_1)_y = \begin{pmatrix} R''(0) \cdot I & 0 \\ 0 & 0 \end{pmatrix}.
\]
Taking some isomorphism $T(S^n)_y \iso \R^n$ and its complexification $T(S^n)_y
\oplus T(S^n)_y \iso \C^n$, one finds that the tangent spaces of $\tau(F_0)$,
$T(0)$, and $F_1$ at $y$ correspond respectively to
\begin{equation} \label{eq:i-angle}
(1 + 2\pi i R''(0)) \R^n,\; i\R^n, \; \R^n \subset \C^n.
\end{equation}
In particular, since $R''(0) < 0$ by $\delta$-wobblyness, the intersection
$\tau(F_0) \cap F_1$ is transverse. At this point we need to recall a fact from
symplectic linear algebra, see e.g.\ \cite[p.\ 40]{lion-vergne}: the
classification of triples of mutually transverse linear Lagrangian subspaces,
up to the action of $Sp(2n)$, is equivalent to the classification of
nondegenerate quadratic forms on $\R^n$, up to $GL(n,\R)$. In particular there
is a finite number of equivalence classes, and deforming a triple continuously
while keeping transversality will not change its equivalence class. One can
clearly deform the three subspaces \eqref{eq:i-angle} in this way to $\R^n$,
$e^{2\pi i/3}\R^n$, $e^{\pi i/3}\R^n$, which proves the last part of the
statement. \qed

The next result links model Dehn twists to exact Lefschetz fibrations. The
connection has been known to algebraic geometers for a very long time, as
attested by the terminology ``Picard-Lefschetz transformations'' used for model
Dehn twists. But while the traditional approach ignores symplectic forms, they
are of course crucial for our purpose.

\begin{lemma} \label{th:model-fibrations}
Fix $\lambda>0$ and $r>0$, and let $\bar{D}(r) \subset \C$ be the closed disc
of radius $r$ around the origin. There is an exact Lefschetz fibration
$(E,\pi)$ over $\bar{D}(r)$, together with a diffeomorphism $\phi : E_r
\rightarrow T(\lambda)$ which respects both the symplectic forms and the exact
one-forms, such that the following holds. Denote by $\rho \in \Sympe(E_r)$ the
symplectic monodromy around $\partial \bar{D}(r)$, in positive sense. Then
$\tau = \phi \circ \rho \circ \phi^{-1} \in \Sympe(T(\lambda))$ is a model Dehn
twist.
\end{lemma}

\proof Take $\C^{n+1}$ with its standard forms $\o_{\C^{n+1}}$,
$\theta_{\C^{n+1}}$ and the function $q: \C^{n+1} \rightarrow \C$, $q(x) =
x_1^2 + \dots + x_{n+1}^2$. Even though $(\C^{n+1},q)$ is clearly not an exact
Lefschetz fibration (its fibres are not even compact), much of what was said in
the previous section carries over to it. The horizontal subspaces
\begin{equation} \label{eq:horizontal-model}
 T(\C^{n+1})_x^h =
 \{X \in \C^{n+1} \suchthat \o_{\C^{n+1}}(X,\ker\, Dq_x) = 0\} = \C \bar{x}
\end{equation}
define a symplectic connection away from the critical point $x = 0$, so that
one has parallel transport maps $q^{-1}(c(a)) \rightarrow q^{-1}(c(b))$ along
paths $c: [a;b] \rightarrow \C^*$. Consider the family of Lagrangian spheres
\begin{equation} \label{eq:model-spheres}
\Sigma_z = \sqrt{z}S^n = \{ (\sqrt{z}y_1, \dots, \sqrt{z}y_{n+1}) \suchthat y
\in S^n \subset \R^{n+1} \} \subset q^{-1}(z),
\end{equation}
$z \neq 0$, which for $z \rightarrow 0$ degenerate to $\Sigma_0 = \{0\} \subset
q^{-1}(0)$. Write $\Sigma^*$ for the union of all $\Sigma_z$, $z \neq 0$, and
$\Sigma = \Sigma^* \cup \Sigma_0$. One computes that
\begin{equation} \label{eq:theta-sigma}
 \theta_{\C^{n+1}} \, | \, \Sigma^* = q^* d^c({-\textstyle\quarter}|z|).
\end{equation}
Actually, the precise formula does not matter much for the moment. All we need
is that $\o_{\C^{n+1}} \,|\, \Sigma^*$ is the pullback by $q$ of some two-form
on $\C^*$, since that implies that parallel transport in any direction in
$\C^*$ takes the $\Sigma_z$ into each other; compare Lemma
\ref{th:submanifold}. The next observation is that if one removes $\Sigma$ then
parallel transport can be extended even to the singular fibre, so that for any
path $c$ in $\C$ one has a canonical symplectic isomorphism
\[
q^{-1}(c(a)) \setminus \Sigma_{c(a)} \rightarrow q^{-1}(c(b)) \setminus
\Sigma_{c(b)}.
\]
Since taking out $\Sigma_0$ removes the critical point, the only possible
problem is that a point in $\C^{n+1} \setminus \Sigma$ might move in horizontal
direction and converge to some point of $\Sigma$, which means that the flow of
some horizontal vector field would not be defined for all time. However, that
cannot happen since there is a function, $h(x) = ||x||^4 - |q(x)|^2$, which
satisfies $dh_x(\C\bar{x}) = 0$, hence is constant horizontally, and with
$h^{-1}(0) = \Sigma$. Therefore one can use parallel transport in radial
direction to trivialize $q: \C^{n+1} \setminus \Sigma \longrightarrow \C$
symplectically. In particular, if $\tilde{\rho}_s: q^{-1}(s) \rightarrow
q^{-1}(s)$ is the monodromy along the circle of radius $s>0$ around the origin,
$\tilde{\rho}_s\,|\,(q^{-1}(s) \setminus \Sigma_s)$ will be isotopic to the
identity in the group of all symplectic automorphisms of $q^{-1}(s) \setminus
\Sigma_s$. We have mentioned all this mainly to motivate the subsequent proof,
which is more computational.

Consider the map
\begin{equation} \label{eq:phi}
\begin{split}
 & \Phi: \C^{n+1} \setminus \Sigma \longrightarrow \C \times (T \setminus T(0)), \\
 & \Phi(x) = (q(x),\sigma_{\alpha/2}(-\im(\hatx)\,||\re(\hatx)||,
 \re(\hatx)\,||\re(\hatx)||^{-1})),
\end{split}
\end{equation}
where $s e^{i\alpha} = q(x)$ are polar coordinates on the base, and $\hatx =
e^{-i\alpha/2}x$. We claim that this is a diffeomorphism fibered over $\C$.
First of all, because of the use of polar coordinates, it is not obvious that
$\Phi$ is well-defined and smooth at $q^{-1}(0) \setminus \Sigma_0$. To dispel
any doubts about that one writes, after some manipulations,
\begin{gather*}
\Phi(x) = \big(q(x),
 - \half \im(x) \beta(x)
 - \half \im(\overline{q(x)}x) \beta(x)^{-1}, \\
 \qquad\qquad  h(x)^{-1/2} \re(x) \beta(x)
 - h(x)^{-1/2} \re(\overline{q(x)}x) \beta(x)^{-1}
 \big),
\end{gather*}
where $\beta(x) = (||x||^2 + h(x)^{1/2})^{1/2}$. The fact that $\Phi$ is a
diffeomorphism on each fibre is easy to see, either directly or by using the
symplectic forms and \eqref{eq:theta-pullback} below. Moreover, if one
restricts $\Phi$ to a fibre over $s > 0$, it extends to a diffeomorphism
$\phi_s: q^{-1}(s) \rightarrow T$. Such an extension does not exist for other
fibres, due to the non-continuity of $\sigma$ at $T(0)$. A computation yields
\begin{equation} \label{eq:theta-pullback}
 (\Phi^{-1})^*\theta_{\C^{n+1}} = \theta_T - \tilde{R}_s(\mu)\, d\alpha, \quad
 \tilde{R}_s(t) = \half t - \half\big(t^2 + s^2/4\big)^{1/2}.
\end{equation}
This implies that the restriction of $\Phi$ to any fibre is symplectic, and
actually maps the respective one-forms into each other. Of course, the same
will then be true for the continuous extensions $\phi_s$, $s>0$. In fact
\eqref{eq:theta-pullback} shows even more: if one restricts to any ray $\alpha
= const.$ in the basis then $\Phi^*\o_T = \o_{\C^{n+1}}$, which means that
$\Phi$ trivializes the symplectic parallel transport on $(\C^{n+1} \setminus
\Sigma,q)$ in radial directions, in accordance with the strategy which we set
out before. Consider $\tilde{\tau}_s = \phi_s \circ \tilde{\rho}_s \circ
\phi_s^{-1} : T \rightarrow T$, where $\tilde{\rho}_s$ is the monodromy map
introduced above. From \eqref{eq:theta-pullback} it follows that
$\tilde{\tau}_s$ restricted to $T \setminus T(0)$ is the time $2\pi$ map of the
Hamiltonian $\tilde{H}_s = \tilde{R}_s(\mu)$. Since $\tilde{R}_s(-t) =
\tilde{R}_s(t) - t$, this is quite close to the case $k = 1$ of Lemma
\ref{th:on-invariant}. The difference is that $\tilde{R}_s(t)$ does not vanish
for $t \gg 0$. Instead, it decays as follows:
\begin{equation} \label{eq:decay}
 0 > \tilde{R}_s(t) \geq -\textstyle{\frac{1}{16}} s^2t^{-1}, \quad
 0 < \textstyle{\frac{d}{dt}}
 \tilde{R}_s(t) \leq \textstyle{\frac{1}{16}} s^2t^{-2}.
\end{equation}
This is good enough to imply that $\tilde{\tau}_s$ is asymptotic to the
identity at infinity, for each $s>0$. It remains to tweak the given data
slightly, so as to produce a honest exact Lefschetz fibration, whose monodromy
is an actual model Dehn twist.

Fix $\lambda>0$, $r>0$. Choose a cutoff function $g \in \smooth(\R^+,\R)$ such
that $g'(t) \geq 0$ everywhere, $g(t) = 0$ for small $t$, and $g(t) = 1$ if $t$
is close to $\lambda$. We claim that there is a unique $\gamma \in
\Omega^1(\C^{n+1})$ with, again in polar coordinates on the base,
\begin{equation} \label{eq:cutoff-form}
(\Phi^{-1})^*\gamma = g(\mu)\tilde{R}_s(\mu) d\alpha.
\end{equation}
Since $\tilde{R}_0(t) = 0$, the function $\tilde{R}_s(t)/s$ extends smoothly to
$s = 0$, $t \neq 0$. Therefore the right hand side of \eqref{eq:cutoff-form}
can be written as $g(\mu)(\tilde{R}_s(\mu)/s)\, s d\alpha$, which means that
$\gamma$ is smooth at least on $\C^{n+1} \setminus \Sigma$. Now $\mu(\Phi(x)) =
(1/2) h(x)^{1/2}$; since $\Sigma = h^{-1}(0)$ and $g(t) = 0$ for small $t$, one
sees that $\Phi^*(g(\mu)\tilde{R}_s(\mu)d\alpha)$ vanishes near $\Sigma$, so
that $\gamma$ extends by zero over $\Sigma$. Set
\begin{equation} \label{eq:define-model}
\begin{cases}
 E = \Phi^{-1}(\bar{D}(r) \times (T(\lambda) \setminus T(0)))
 \cup (\Sigma \cap q^{-1}(\bar{D}(r))), \\
 \pi = q|E: E \longrightarrow \bar{D}(r), \\
 \Theta = (\theta_{\C^{n+1}} + \gamma) \,|\,E, \quad
 \Omega = d\Theta, \\
 \phi = \phi_r|E_r: E_r \longrightarrow T(\lambda).
\end{cases}
\end{equation}
$E \subset \C^{n+1}$ is cut out by the inequalities $h(x) \leq 4\lambda^2$,
$|q(x)| \leq r$. This makes it easy to show that it is a compact manifold with
corners, whose boundary faces are $\partial_vE = \pi^{-1}(\partial \bar{D}(r))$
and $\partial_hE = \{h(x) = 4\lambda^2\} = \{\mu(\Phi(x)) = \lambda\} =
\Phi^{-1}(\bar{D}(r) \times \partial T(\lambda))$. Because $\gamma$ vanishes
when restricted to any fibre, $\Omega|TE^v_x = \o_{\C^{n+1}} | \ker\, Dq_x$ is
nondegenerate for all $x$. In
\begin{equation} \label{eq:tilde-pullback}
\Theta\,|\,(E \setminus \Sigma) = \Phi^*(\theta_T + (g(\mu)-1)
\tilde{R}_s(\mu)\, d\alpha)
\end{equation}
the second term on the right hand side vanishes close to $\bar{D}(s) \times
\partial T(\lambda)$, so that $\Phi$ provides a trivialization near
$\partial_hE$ in the sense introduced in the previous section. Moreover, since
$\gamma$ vanishes near the critical point $x = 0$, equipping $E$ with the
standard complex structure $J_0$ near that point, and $\bar{D}(r)$ with the
standard complex structure $j_0$, turns $(E,\pi)$ into an exact Lefschetz
fibration. As for the statement about the monodromy $\rho$, one can repeat the
argument above, using \eqref{eq:tilde-pullback} instead of
\eqref{eq:theta-pullback}. This shows that $\tau = \phi \circ \rho \circ
\phi^{-1}$, when restricted to $T(\lambda) \setminus T(0)$, is the time $2\pi$
map of $R_r(\mu)$, where
\begin{equation} \label{eq:r-function}
R_r(t) = (1-g(t))\tilde{R}_r(t);
\end{equation}
by definition, this is a model Dehn twist. \qed

Let $(M,\o,\theta)$ be an exact symplectic manifold. A {\em
framed}\footnote{This has little or nothing to do with the usual topological
notion of framed manifold.} {\em exact Lagrangian sphere} is an exact
Lagrangian submanifold $L \subset M$ together with an equivalence class $[f]$
of diffeomorphisms $f: S^n \rightarrow L$. Here $f_1,f_2$ are equivalent iff
$f_2^{-1}f_1$ can be deformed inside $\Diff(S^n)$ to an element of $O(n+1)$. To
any such $(L,[f])$ one associates a Dehn twist $\tau_{(L,[f])} \in \Sympe(M)$
as follows. Choose a representative $f$ and extend it to a symplectic embedding
$\iota: T(\lambda) \rightarrow M$ for some $\lambda>0$. Take a model Dehn twist
$\tau$ which is supported in the interior of $T(\lambda)$, and define
\[
\tau_{(L,[f])} = \begin{cases}
 \iota \circ \tau \circ \iota^{-1} & \text{on $\im(\iota)$,} \\
 \id & \text{elsewhere.}
\end{cases}
\]
The exactness of $\tau_{(L,[f])}$ follows from that of $L$; moreover, the
analogue of \eqref{eq:self-shift} holds. It is not difficult to show that the
isotopy class $[\tau_{(L,[f])}] \in \pi_0(\Sympe(M))$ is independent of the
choices made in the definition. In contrast, it is unknown whether a change of
the framing $[f]$ can affect $[\tau_{(L,[f])}]$; if the answer is negative the
notion of framing could be dropped altogether, but while the question is open
one cannot do without it. Still, for the sake of brevity we will often omit
framings from the notation and write $\tau_L$ instead of $\tau_{(L,[f])}$. We
will say that $\tau_L$ is $\delta$-wobbly if the local model it is constructed
out of has this property.

\begin{prop} \label{th:standard-fibrations}
Let $(L,[f])$ be a framed exact Lagrangian sphere in $M$. Fix some $r>0$. There
is an exact Lefschetz fibration $(E^L,\pi^L)$ over $\bar{D}(r)$ together with
an isomorphism $\phi^L: E^L_r \rightarrow M$ of exact symplectic manifolds,
such that if $\rho^L$ is the symplectic monodromy around $\partial \bar{D}(r)$,
then $\tau_L = \phi^L \circ \rho^L \circ (\phi^L)^{-1}$ is a Dehn twist along
$(L,[f])$.
\end{prop}

We will prove this under the assumption that there is an embedding $\iota:
T(\lambda) \rightarrow M$ as before, with $\iota^*\theta = \theta_T$. This is
not really a restriction, since one can always satisfy it by adding the
derivative of a function to $\theta$, which does not change $M$ up to exact
symplectic isomorphism. On the other hand, it allows us to make the statement
slightly sharper: $\phi^L$ will map the one-forms on $E^L_r$ and $M$ into each
other, and the Dehn twist obtained from the monodromy of $(E^L,\pi^L)$ will be
one constructed using the given embedding $\iota$.

\proof Take $(E,\pi)$ and $\phi$ from Lemma \ref{th:model-fibrations}, with the
given $r$ and $\lambda$. We will construct $(E^L,\pi^L)$ by attaching a trivial
piece to $(E,\pi)$. By construction, there is a neighbourhood $N \subset E$ of
$\partial_hE$, a neighbourhood $V$ of $\partial T(\lambda)$ in $T(\lambda)$,
and a diffeomorphism $\Phi: N \rightarrow \bar{D}(r) \times V$ fibered over
$\bar{D}(r)$, such that $\Theta|N = \Phi^*\theta_T$. Moreover, $\Phi$ agrees
with $\phi$ on $N \cap E_r$. Set
\[
E^L = E \cup_{\sim} \bar{D}(r) \times (M \setminus \iota(T(\lambda) \setminus
V)),
\]
where $\sim$ identifies $N$ with $\bar{D}(r) \times \iota(V)$ through $(\id
\times \iota) \circ \Phi$. One similarly defines $\pi^L$ from $\pi$ and the
projection $\bar{D}(r) \times M \rightarrow \bar{D}(r)$. The forms
$\Omega^L,\Theta^L$ on $E^L$ come from the corresponding ones on $E$ and the
pullbacks of $\o,\theta$ on the trivial part. $\phi^L$ is constructed from
$\iota \circ \phi$ and the identity map. The complex structure near the
critical point and critical value are inherited from $E$. All properties stated
above are obvious from the construction and the definition of Dehn twists. \qed

We call the $(E^L,\pi^L)$ {\em standard fibrations}. For future reference, we
will now state certain properties which these fibrations inherit from the local
model $(E,\pi)$, and which depend on the details of its construction.

\begin{lemma} \label{th:standard-properties}
Any standard fibration $(E^L,\pi^L)$ has the following properties.
\begin{romanlist}
\item \label{item:lagrangian-subbundle}
There is a closed subset $\Sigma^L \subset E^L$ such that
\[
\Sigma^L_z = \Sigma^L \cap E^L_z =
 \begin{cases}
 \text{is an embedded $n$-sphere} & \text{if $z \neq 0$,} \\
 \text{is the unique critical point $x_0 \in E^L_0$} &
 \text{if $z = 0$.}
 \end{cases}
\]
In fact $(\Sigma^L)^* = \Sigma^L \setminus \Sigma^L_0$ is a smooth $n$-sphere
bundle over $\bar{D}(r) \setminus \{0\}$, and satisfies
\begin{equation} \label{eq:theta-sigma-l}
 \Theta^L \, | \, (\Sigma^L)^* =
 (\pi^L)^* d^c(-\textstyle\quarter|z|).
\end{equation}
As in the discussion following \eqref{eq:theta-sigma}, this implies that each
$\Sigma^L_z \subset E^L_z$, $z \neq 0$, is an exact Lagrangian submanifold, and
that symplectic parallel transport within $\bar{D}(r) \setminus \{0\}$ carries
these spheres into each other. Moreover, $\phi^L(\Sigma^L_r) = L$.
\item \label{item:preferred-charts}
There are holomorphic Morse charts $(\xi,\Xi)$ around the unique critical point
$x_0 \in E^L_0$, such that $\xi$ is the inclusion $U \hookrightarrow
\bar{D}(r)$ of some neighbourhood $U \subset \C$ of the origin; $\Xi^*\Theta^L
= \theta_{\C^{n+1}}$ and $\Xi^*\Omega^L = \omega_{\C^{n+1}}$; and
$\Xi^{-1}(\Sigma^L_z) = \Sigma_z$ for all sufficiently small $z \in \C$.
\item \label{item:nonnegatively-curved}
$(E^L,\pi^L)$ has nonnegative curvature.
\item \label{item:precise-twist}
$\tau_L = \phi^L \circ \rho^L \circ (\phi^L)^{-1}$ is the Dehn twist defined
using $\iota$ and the function $R_r$ from \eqref{eq:r-function}; in particular
$R_r(0) = -r/4$. By making $r$ smaller while keeping all other choices in the
construction fixed, one can achieve that $\tau_L$ is $\delta$-wobbly for an
arbitrary $\delta$.
\end{romanlist}
\end{lemma}

\proof \ref{item:lagrangian-subbundle} From \eqref{eq:define-model} one sees
that $\Sigma \cap q^{-1}(\bar{D}(r)) \subset E$, and that $\Theta =
\theta_{\C^{n+1}}$ in a neighbourhood of that subset. $\Sigma^L$ is defined to
be the image of this in $E^L$, so that \eqref{eq:theta-sigma-l} is a
consequence of \eqref{eq:theta-sigma}. \ref{item:preferred-charts} The
definition of $\xi,\Xi$ is obvious, and the claim about $\Xi^*\Theta^L$ follows
from the fact that in \eqref{eq:define-model} $\Theta = \theta_{\C^{n+1}}$ near
the critical point. \ref{item:nonnegatively-curved} One can compute the
curvature of $(E,\pi)$ from \eqref{eq:theta-pullback} and
\eqref{eq:cutoff-form}; it turns out to be
\[
\Phi^*\big( (g(\mu)-1) \textstyle\frac{\partial}{\partial s} \tilde{R}_s(\mu)
\big) ds \wedge d\alpha,
\]
which is $\geq 0$ everywhere. This implies the same property for $(E^L,\pi^L)$.
\ref{item:precise-twist} The nontrivial statement is $\delta$-wobblyness, which
requires that we take another look at the function $R_r$. There is a $t_0>0$
such that $g(t) = 0$ for $t \in [0;t_0]$, and in that interval $R_r(t) =
\tilde{R}_r(t)$, so that $(\partial/\partial t) R_r(t) > 0$,
$(\partial^2/\partial t^2) R_r(t) < 0$. For $t \in [t_0;\lambda]$ one estimates
using \eqref{eq:decay} that
\[
0 \leq {\textstyle\frac{\partial}{\partial t}}R_r(t) = -g'(t) \tilde{R}_r(t) +
(1-g(t)) {\textstyle\frac{\partial}{\partial t}} \tilde{R}_r(t) \leq
\textstyle{\frac{r^2}{16}} (||g'|| t_0^{-1} + t_0^{-2}).
\]
By choosing $r$ small while keeping $g$ and hence $t_0$ fixed, one can make
this $< \delta$ for an arbitrary $\delta$. \qed

\subsection{Vanishing cycles\label{sec:picard-lefschetz}}

We have seen that all Dehn twists can be realized as monodromy maps (of
standard fibrations). Conversely, the geometry of any exact Lefschetz fibration
can be understood in terms in Dehn twists. This is again well-known on a
topological level, and our exposition repeats the classical arguments while
paying more attention to symplectic forms. The results will not be used again
in this paper, but they are important in applications.

Let $(E,\pi)$ be an exact Lefschetz fibration over $S$. Take $z_0 \in S^\crit$,
and the unique critical point $x_0 \in E_{z_0}$. Let $c: [a;b] \rightarrow S$
be a smooth embedded path with $c(b) = z_0$, $c^{-1}(S^\crit) = \{b\}$; and let
$\rho_{c|[s,s']}: E_{c(s)} \rightarrow E_{c(s')}$ be the parallel transport
maps along it, which is defined for all $s \leq s' < b$. Following a suggestion
of Donaldson, we define
\begin{equation} \label{eq:vanishing-ball}
B_c = \{x \in E_{c(s)}, \; a \leq s < b \suchthat \lim_{s' \rightarrow b}
\rho_{c|[s;s']}(x) = x_0 \} \cup \{ x_0 \} \subset E.
\end{equation}

\begin{lemma} \label{th:vanishing-ball}
$B_c$ is an embedded closed $(n+1)$-ball, with $\partial B_c = B_c \cap
E_{c(a)}$. The function $p = c^{-1} \circ \pi: B_c \rightarrow [a;b]$ has $x_0$
as its unique critical point, which is a nondegenerate local maximum. Moreover,
$\Omega|B_c = 0$.
\end{lemma}

\proof Put a symplectic form $\Omega + \pi^*\beta$ on $E$ as in Lemma
\ref{th:add-base}. Choose an oriented embedding $\tilde{c}:
(a-\epsilon;b+\epsilon) \times (-\epsilon;\epsilon) \rightarrow S$ such that
$\tilde{c}(0,t) = c(t)$. Let $h$ be the function defined on $\im(\tilde{c})$
with $h(\tilde{c}(s,t)) = -t$. The Hamiltonian vector field $X$ of $H = h \circ
\pi$ has the following properties:
\begin{condensedlist}
\item \label{item:h-one}
it is horizontal everywhere, $X_x \in TE^h_x$;
\item \label{item:h-two}
for each $x \neq x_0$, $D\pi(X)_x$ is a positive multiple of
$\partial\tilde{c}/\partial s$;
\item \label{item:h-three}
$x_0$ is a hyperbolic stationary point of $X$, with $n+1$ positive and negative
eigenvalues.
\end{condensedlist}
\ref{item:h-one} and \ref{item:h-two} are straightforward; \ref{item:h-three}
can be seen by looking at $X$ in a holomorphic Morse chart, where it is
$J_0\nabla H$. Let $\tilde{B} \subset E$ be the stable manifold of $x_0$. It is
an open $(n+1)$-ball, lies in $\pi^{-1}\tilde{c}((a-\epsilon;b] \times \{0\})$,
and the projection
\[
 \tilde{p} = \tilde{c}^{-1}\pi : \tilde{B} \longrightarrow
 (a-\epsilon;b] \times \{0\}
\]
is a proper map. The tangent space of $\tilde{B}$ at $x_0$ is the negative
eigenspace of $DX_{x_0}$; one shows easily that $D^2\tilde{p}|T\tilde{B}_{x_0}$
differs from $DX|T\tilde{B}_{x_0}$ only by a positive constant, which means
that $x_0$ is a nondegenerate maximum of $\tilde{p}$. There are no other
critical points, because elsewhere $d\tilde{p}(X)_x>0$ by \ref{item:h-two}.
Finally, since the flow $X$ is symplectic and contracts the tangent spaces of
$\tilde{B}$, these must be Lagrangian subspaces, so $(\Omega +
\pi^*\beta)|\tilde{B} = \Omega|\tilde{B} = 0$.

We claim that $B_c = \tilde{B} \cap \pi^{-1}(\im(c))$, which implies all
desired properties of $B_c$. Let $Y$ be the vector field on
$\pi^{-1}(\im(\tilde{c})) \setminus \{x_0\}$ which is horizontal and satisfies
$D\pi(Y) = \partial\tilde{c}/\partial s$. From \ref{item:h-one} and
\ref{item:h-two} above one sees that $Y = gX$ for some function $g$, bounded
from below by a positive constant, and which goes to $\infty$ as one approaches
$x_0$. Therefore the orbits of $X$ and $Y$ coincide, except of course for
$\{x_0\}$. Since $Y|\pi^{-1}(\im(c))$ defines parallel transport along $c$, the
claim follows. \qed

\begin{lemma} \label{th:vanishing-cycle}
$V_c = \partial B_c$ is an exact Lagrangian sphere in $E_{c(a)}$, and comes
with a canonical framing.
\end{lemma}

\proof Since $B_c$ is a ball and $\Theta|B_c$ is closed, there is a unique $K
\in \smooth(B_c,\R)$ with $K(x_0) = 0$ and $dK = \Theta|B_c$. The restriction
$K_{V_c} = K|V_c$ makes $V_c$ into an exact Lagrangian submanifold. It remains
to explain the framing. In any chart on $B_c$ around $x_0$, the level sets
$p^{-1}(s')$ of the function from Lemma \ref{th:vanishing-ball}, for $s'$ close
to $b$, will be strictly convex hypersurfaces, so that one can map them to
$S^n$ by radial projection. On the other hand, $V_c = p^{-1}(a)$ can be
identified with $p^{-1}(s')$ by the gradient flow with respect to some metric.
Combining these two maps gives a diffeomorphism $V_c \rightarrow S^n$, which is
unique up to isotopy and action of $O(n+1)$; its inverse is our framing. \qed

The framed exact Lagrangian sphere $V_c \subset E_{c(a)}$ is called the {\em
vanishing cycle} associated to $c$. It exists more generally for any path $c:
[a;b] \rightarrow S$ which is smooth, not necessarily embedded, but still
satisfies
\begin{equation} \label{eq:path}
c^{-1}(S^\crit) = \{b\}, \quad c'(b) \neq 0.
\end{equation}
To construct it in this situation, one first extends $c$ to a map
$(a-\epsilon;b+\epsilon) \times (-\epsilon;\epsilon) \rightarrow S$ which is a
local oriented diffeomorphism at $(b,0)$, pulls back $(E,\pi)$ by that map, and
then applies Lemma \ref{th:vanishing-ball} to the pullback exact Lefschetz
fibration. Note also that deforming $c$ smoothly, rel endpoints, within the
class \eqref{eq:path} yields an exact Lagrangian isotopy of the corresponding
vanishing cycles, which is compatible with their framings.

\includefigure{doubling}{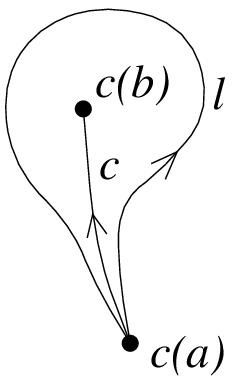}{hb}%
\begin{prop} \label{th:picard-lefschetz}
Let $c$ be a path satisfying \eqref{eq:path}, and $l$ the loop in $S \setminus
S^\crit$ obtained by ``doubling'' $c$, as in Figure \ref{fig:doubling}. Then
the monodromy around $l$ is isotopic to the Dehn twist along the vanishing
cycle $V_c$:
\begin{equation} \label{eq:picard-lefschetz}
[\rho_l] = [\tau_{V_c}] \in \pi_0(\Sympe(E_{c(a)})).
\end{equation}
\end{prop}

\proof We start with a rather special case. Let $(E,\pi)$ be an exact Lefschetz
fibration with base $\bar{D}(r)$ for some $r>0$, and which has exactly one
critical point $x_0 \in E_0$. Write $M = E_r$. In addition, $x_0$ should admit
a holomorphic Morse chart $(\xi,\Xi)$ where $\Xi$ is defined on $W = \{x \in
\C^{n+1} \suchthat ||x|| \leq 2\sqrt{r}, \; |q(x)| \leq r\}$, where $\xi =
\id_{\bar{D}(r)}$, and such that $\Xi^*\Omega$, $\Xi^*\Theta$ are standard.
This means that the symplectic geometry of $(E,\pi)$ equals that of the local
model $(\C^{n+1},q)$ discussed in Lemma \ref{th:model-fibrations}, at least in
a suitable neighbourhood of $x_0$. In particular, each fibre $E_z$, $z \neq 0$,
contains an exact Lagrangian sphere $\Sigma_z^E = \Xi(\Sigma_z)$, degenerating
to $\Sigma_0^E = \{x_0\}$. Write $L = \Sigma_r^E \subset M$; this inherits an
obvious framing from $\Sigma_r = \sqrt{r}S^n$. We claim that the monodromy
around $\partial \bar{D}(r)$, denoted by $\rho \in \Sympe(M)$, is isotopic to
$\tau_L$.

From \eqref{eq:theta-sigma} one deduces that parallel transport in $\bar{D}(r)
\setminus \{0\}$ carries the $\Sigma_z^E$ into each other. Moreover, by the
same argument as in Lemma \ref{th:model-fibrations}, if one removes $\Sigma^E =
\bigcup_z \Sigma^E_z$ from $E$, parallel transport can be extended over the
singular fibre. Using parallel transport in radial directions one constructs a
trivialization $\Phi^E: E \setminus \Sigma^E \longrightarrow \bar{D}(r) \times
(M \setminus L)$ which is the identity on the fibre over $r$; this is such
that, in radial coordinates $z = se^{i\alpha}$ on the basis,
$((\Phi^E)^{-1})^*\Theta = \theta - R^E \wedge d\alpha + dS^E$ for some
functions $R^E = R^E(s,\alpha,x)$ and $S^E = S^E(s,\alpha,x)$ which vanish for
$x$ close to $\partial M$. This implies that $\rho$ restricted to $M \setminus
L$ is the time-$2\pi$ map of the flow generated by the time-dependent
Hamiltonian $(\alpha,x) \mapsto R^E(r,\alpha,x)$. Clearly, $[\rho] \in
\pi_0(\Sympe(M))$ depends only on the behaviour of this function in an
arbitrarily small neighbourhood of $L$, or equivalently on $(E,\pi)$ and
$\Omega$ close to $\Sigma^E$; since $\Sigma^E \subset \Xi(W)$ by definition,
this can be determined from the local model $(\C^{n+1},q)$. It remains to spell
out the computation.

Define an embedding $\iota: T(\lambda) \rightarrow M$, for some $\lambda>0$, by
combining $\Xi_r = \Xi|(W \cap q^{-1}(r)) \rightarrow M$ with the inverse of
the isomorphism $\phi_r: q^{-1}(r) \rightarrow T$ from the proof of Lemma
\ref{th:model-fibrations}. This satisfies $\iota^*\theta = \theta_T$,
$\iota^*\o = \o_T$, and $\iota(T(0)) = L$. Because both $\Phi^E$ and the
trivialization $\Phi$ from \eqref{eq:phi} are defined by radial parallel
transport, there is a commutative diagram
\[
\xymatrix{
 {\C^{n+1} \setminus \Sigma} \ar[r]^-{\Phi} &
 {\C \times (T \setminus T(0))} \\
 {W \setminus \Sigma} \ar@{^{(}->}[u] \ar[d]_{\Xi\,|\,W \setminus \Sigma} &
 {\bar{D}(r) \times (T(\lambda) \setminus T(0))} \ar@{^{(}->}[u]
 \ar[d]^{\id \times (\iota \,|\, T(\lambda) \setminus T(0))} \\
 {E \setminus \Sigma^E} \ar[r]^-{\Phi^E} & {\bar{D}(r) \times (M \setminus L).}
}
\]
From this, $\Xi^*\Omega = \o_{\C^{n+1}}$ and \eqref{eq:theta-pullback} it
follows that
\begin{align*}
 & (\id \times \iota)^*(\theta_T - R^E d\alpha + dS^E)
 = (\id \times \iota)^*((\Phi^E)^{-1})^*\Theta \\
 & = (\Phi^{-1})^*\theta_{\C^{n+1}} = \theta_T - \tilde{R}_s(\mu) \wedge d\alpha,
\end{align*}
and in particular that $R^E(r,\alpha,\iota(y)) = \tilde{R}_r(\mu(y))$ for all
$\alpha$. In view of the discussion above, and the fact that $\tilde{R}_r(-t) =
\tilde{R}_r(t) - t$, this implies the desired equality $[\rho] = [\tau_L]$. Now
consider the path $c: [0;r] \rightarrow \bar{D}(r)$, $c(s) = r-s$. Because
symplectic parallel transport takes the $\Sigma^E_z$ into each other, the ball
$B_c$ from Lemma \ref{th:vanishing-ball} must be the union of $\Sigma^E_z$ for
all $z \in [0;r]$, and the vanishing cycle is $V_c = \Sigma^E_r = L \subset M$,
with the same framing as before. Therefore, what we have done up to now proves
\eqref{eq:picard-lefschetz} for this special class of exact Lefschetz
fibrations $(E,\pi)$ and the particular path $c$.

More generally, let $(E,\pi)$ be a exact Lefschetz fibration with arbitrary
base $S$, and $x_0 \in E_{z_0}$ a critical point which admits a holomorphic
Morse chart $(\xi,\Xi)$ such that $\Xi^*\Omega$, $\Xi^*\Theta$ are standard. By
restricting to a suitably small disc around $z_0$ in the base, and making the
domains of $\xi,\Xi$ smaller, one can arrive at the situation considered
before. This means that \eqref{eq:picard-lefschetz} is true at least for one
(short) path $c$ with endpoint $z_0$. But from that it follows easily for all
other paths with the same endpoint.

It remains to remove the assumption concerning $\Xi$. Let $(E,\pi)$ be an
arbitrary exact Lefschetz fibration, $c$ a path as in \eqref{eq:path}, and $l$
a corresponding loop. With respect to the critical point in $E_{c(b)}$, take
smooth families $\Omega^\mu$, $\Theta^\mu$ as in Lemma
\ref{th:local-deformation}. These can be chosen such that in a neighbourhood of
$E_{c(a)}$, $\Omega^\mu = \Omega$ and $\Theta^\mu = \Theta$ for all $\mu$. The
corresponding vanishing cycles $V^\mu_c \subset E_{c(a)}$ form a smooth isotopy
of framed exact Lagrangian spheres, and similarly there is a smooth family of
monodromies $\rho_l^{\mu}$. For $\mu = 1$ our previous assumption about
holomorphic Morse charts is satisfied, and hence $[\rho_l] = [\rho_l^0] =
[\rho_l^1] = [\tau_{V^1_c}] = [\tau_{V^0_c}] = [\tau_{V_c}] \in
\pi_0(\Sympe(E_{c(a)}))$. \qed

\newpage
\section[Sections]{Pseudo-holomorphic sections\label{ch:two}}

Before, we have considered exact Lefschetz fibrations as geometric
objects in the sense of elementary symplectic geometry; now we will
apply the theory of pseudo-holomorphic curves to them. By today's
standards, the necessary analysis is rather unsophisticated. The
basic Gromov-type invariant which will be introduced first uses only
the most familiar techniques, on the level of the book
\cite{mcduff-salamon}. Of course ours is a Lagrangian boundary value
problem, but the only part of the analysis which specifically
concerns the boundary is bubbling off of holomorphic discs, which was
addressed by Floer \cite{floer88c} and Oh \cite{oh92} (some more
recent expositions are \cite[Section A.4.3]{ivashkovich-shevchishin},
\cite{frauenfelder01}). Later, when considering relative invariants,
we will use essentially the same analysis as in \cite{schwarz95},
\cite{piunikhin-salamon-schwarz94}, \cite{desilva98}, even though the
geometric setting is rather different. The remaining material
(Sections \ref{sec:horizontal}, \ref{sec:vanishing}, and
\ref{sec:horizontal-two}) is more specifically designed for
application to the exact sequence, and has a greater claim to
originality.

\subsection{A simple invariant\label{sec:simple-invariant}}

Let $(E^{2n+2},\pi)$ be an exact Lefschetz fibration over $S$. A {\em
Lagrangian boundary condition} is an $(n+1)$-dimensional submanifold $Q \subset
E|\partial S$ which is disjoint from $\partial_hE$ and such that $\pi|Q: Q
\rightarrow \partial S$ is a submersion, together with $\kappa_Q \in
\Omega^1(\partial S)$ and $K_Q \in \smooth(Q,\R)$, satisfying $\Theta|Q =
\pi^*\kappa_Q + dK_Q$. This implies that each $Q_z = Q \cap E_z$, $z \in
\partial S$, is an exact Lagrangian submanifold of $E_z$, with function
$K_Q|Q_z$ (one could consider a more general situation in which the condition
on $\pi|Q$ is dropped, so that the $Q_z$ could have singularities; but we will
not do that). From Lemma \ref{th:submanifold} one sees that parallel transport
along $\partial S$ carries the $Q_z$ into each other, in the ordinary sense of
the word; this holds in the ``exact'' sense, that is to say as exact Lagrangian
submanifolds, iff $\kappa_Q = 0$. Note also that if one equips $E$ with a
symplectic form $\Omega + \pi^*\beta$ as in Lemma \ref{th:add-base}, $Q$ itself
becomes a Lagrangian submanifold. The aim of this section is to introduce, in
the case where $S$ is compact with $\partial S \neq \emptyset$, a Gromov-type
invariant
\[
\Phi_1(E,\pi,Q) \in H_*(Q_\zeta;\Z/2),
\]
where $\zeta$ is some point of $\partial S$. In a nutshell, this is the cycle
represented by the values at $\zeta$ of pseudo-holomorphic sections of $E$ with
boundary in $Q$.

To begin, we remind the reader of a class of almost complex structures suitable
for exact symplectic manifolds. The Liouville vector field $N$ on such a
manifold, $i_N\o = \theta$, defines a collar $\R^- \times \partial M
\hookrightarrow M$. Let $\sigma$ be the function on a neighbourhood of
$\partial M$ whose composition with $\R^- \times \partial M \hookrightarrow M$
is projection to the $\R^-$ factor. An $\o$-compatible almost complex $J$ is
called convex near the boundary if $\theta \circ J = d(e^\sigma)$ near
$\partial M$. This implies $d(d(e^\sigma) \circ J) = -\o$, which serves to
control the behaviour of $J$-holomorphic curves. It is well-known that the
space of these $J$ is contractible (in particular, nonempty).

Given an exact Lefschetz fibration $(E,\pi)$, one can consider the Liouville
vector field on each fibre, and this gives rise to a function $\sigma$ on a
neighbourhood of $\partial_hE$ in $E$. Choose a complex structure $j$ on the
base $S$ (whenever we do that, now or later, $j$ is assumed to be positively
oriented and equal to $j_0$ in some neighbourhood of $S^\crit$). An almost
complex structure $J$ on $E$ is called {\em compatible relative to $j$} if
\begin{condensedlist}
\item
$J = J_0$ in a neighbourhood of $E^\crit$;
\item \label{item:projection}
$D\pi \circ J = j \circ D\pi$;
\item \label{item:partially-compatible}
$\Omega(\cdot,J\cdot) | TE_x^v$ is symmetric and positive definite for any $x
\in E$;
\item \label{item:convexity}
in a neighbourhood of $\partial_hE$, $J(TE^h) = TE^h$ and $\Theta \circ J =
d(e^\sigma)$.
\end{condensedlist}
To see more concretely the meaning of these conditions, take $x \notin E^\crit$
and split $TE_x = TE_x^h \oplus TE_x^v \iso TS_z \oplus TE_x^v$, $z = \pi(x)$.
One can then write
\begin{equation} \label{eq:compatible-j}
 J_x = \begin{pmatrix} j_z & 0 \\ J_x^{vh} & J_x^{vv} \end{pmatrix}
\end{equation}
where $J_x^{vv} \in \End(TE_x^v)$ is a complex structure compatible with
$\Omega|TE_x^v$, and $J_x^{vh} \in \Hom(TS_z,TE_x^v)$ is $\C$-antilinear with
respect to $j$ and $J^{vv}$. This is a reformulation of \ref{item:projection},
\ref{item:partially-compatible}. An immediate consequence is the following
result, which sharpens Lemma \ref{th:add-base}:

\begin{lemma} \label{th:tame}
Let $J$ be compatible relative to $j$. Then for any sufficiently positive
$\beta \in \Omega^2(S)$, $\Omega + \pi^*\beta$ tames $J$. \qed
\end{lemma}

Because of the lack of antisymmetry in \eqref{eq:compatible-j}, $J$ will not be
compatible with $\Omega + \pi^*\beta$ in the ordinary sense of the word, unless
$J^{vh} = 0$; we will return to this more restricted class of almost complex
structures in the next section. Continuing with the analysis of the conditions
above, suppose now that $x \in E$ is sufficiently close to $\partial_hE$. More
precisely, we require that $\sigma(x)$ is defined and that \ref{item:convexity}
applies. There is then a further splitting
\begin{equation} \label{eq:fine-splitting}
TE^v_x \iso \R N \oplus \R R \oplus (\ker\, \Theta \cap \ker\, d\sigma \cap
TE^v_x),
\end{equation}
where $N$ is the Liouville vector field and $R$ is the Hamiltonian vector field
of $e^\sigma|E_z$. The two parts of \ref{item:convexity} say that $J_x^{vh} =
0$ and that with respect to \eqref{eq:fine-splitting},
\begin{equation} \label{eq:boundary-j}
 J_x^{vv} = \begin{pmatrix} 0 & -1 & 0 \\ 1 & 0 & 0 \\ 0 & 0 & * \end{pmatrix}.
\end{equation}
The pointwise analysis which we have just carried out can be recast in terms of
sections of fibre bundles, and one sees then that the space $\JJ(E,\pi,j)$ of
almost complex structures which are compatible relative to $j$ is contractible.

Suppose now that we have a Lagrangian boundary condition $Q$. The theory of
pseudo-holomorphic sections with boundary in $Q$ fits into a familiar framework
of infinite-dimensional manifolds and maps. We will now review this, on a
formal level, that is to say using $\smooth$ spaces and without assuming that
$S$ is compact. In that sense, $\JJ(E,\pi,j)$ is an infinite-dimensional
manifold; its tangent space at $J$ consists of sections $Y \in
\smooth(\End(TE))$ which are zero near $E^\crit$ and can be written, in
parallel with \eqref{eq:compatible-j}, as
\begin{equation} \label{eq:tangent-compatible}
Y_x = \begin{pmatrix} 0 & 0 \\ Y^{vh}_x & Y^{vv}_x \end{pmatrix}
\end{equation}
with $Y^{vv}_x$ an infinitesimal deformation of the compatible complex
structure $J^{vv}_x$, and $Y^{vv}_xJ^{vh}_x + J^{vv}_x Y^{vh}_x = - Y^{vh}_x
j_z$. There is a further requirement about $Y$ near $\partial_hE$, the
linearization of \ref{item:convexity}, which we leave to the reader to write
down. The space $\BB$ of sections $u: S \rightarrow E$ satisfying $u(\partial
S) \subset Q$ is also an infinite-dimensional manifold, with $T\BB_u = \{ X \in
\smooth(u^*TE^v) \suchthat X_z \in T(Q_z) \text{ for all } z \in \partial S\}$;
note that $u^*TE^v \rightarrow S$ is really a vector bundle, since $u$ as a
smooth section of $\pi$ avoids $E^\crit$. Consider the infinite-dimensional
vector bundle $\EE \rightarrow \BB \times \JJ(E,\pi,j)$ whose fibre at $(u,J)$
is $\Omega^{0,1}(u^*TE^v)$, the space of $(0,1)$-forms on $(S,j)$ with values
in $u^*(TE^v,J|TE^v)$. It has a canonical section $\bar\partial^{univ}(u,J) =
\half (Du + J \circ Du \circ j)$, and the zero set $\MM^{univ} =
(\bar\partial^{univ})^{-1}(0)$ consists of pairs $(u,J)$ such that $u$ is
$(j,J)$-holomorphic. We denote by $\bar\partial_J$, $\MM_J$ the restrictions of
$\bar\partial^{univ}$, $\MM^{univ}$ to a fixed $J \in \JJ(E,\pi,j)$. The
derivative of $\bar\partial_J$ at $u \in \MM_J$ is a map $D_{u,J}: T\BB_u
\rightarrow \EE_{u,J}$. An explicit formula is
$
 D_{u,J}(X) = \half (L_{\tilde{X}}J) \circ Du \circ j
$
where $\tilde{X}$ is any section of $TE^v$, defined on a neighbourhood of
$\im(u)$ in $E$, such that $u^*\tilde{X} = X$. After choosing a torsion-free
connection $\nabla$ on $TE$ (away from $E^\crit$) which preserves the
integrable subbundle $TE^v$, one transforms this into the more familiar
expression
\begin{equation} \label{eq:connection-formula}
 D_{u,J}(X) = \bar{\partial}_{J,u^*\nabla} X +
 \half (\nabla_X J) \circ Du \circ j,
\end{equation}
in which $\bar{\partial}_{J,u^*\nabla} = (u^*\nabla)^{0,1}$ is the
$\bar\partial$-operator associated to the pullback connection $u^*\nabla$ on
$u^*(TE^v,J|TE^v)$. The derivative of $\bar\partial^{univ}$, $D^{univ}_{u,J} :
T\BB_u \times T\JJ(E,\pi,j)_J \longrightarrow \EE_{u,J}$, is obviously
\begin{equation} \label{eq:universal-operator}
\begin{split}
 & D^{univ}_{u,J}(X,Y) = D_{u,J}(X) + \half Y \circ Du \circ j.
\end{split}
\end{equation}

From this point onwards, we pass to a more realistic situation, and assume that
$S$ is compact with $\partial S \neq \emptyset$. Our first observation is a
consequence of \ref{item:convexity}. Informally speaking, it says that for the
purposes of pseudo-holomorphic sections, the boundary $\partial_hE$ of the
fibres can be ignored.

\begin{lemma} \label{th:convexity}
For every $J \in \JJ(E,\pi,j)$ there is a compact subset $K \subset E \setminus
\partial_hE$ such that all $u \in \MM_J$ satisfy $u(S) \subset K$.
\end{lemma}

\proof Let $W \subset E$ be a closed neighbourhood of $\partial_hE$ which,
under the collar embedding provided by the Liouville vector fields on the
fibres, corresponds to $[-\epsilon;0] \times \partial_hE$ for some
$\epsilon>0$. By definition $\sigma(W) = [-\epsilon;0]$. After possibly making
$W$ and $\epsilon$ smaller, we may assume that $W \cap Q = \emptyset$, that
$\Omega|TE^h$ vanishes on $W$, and that \ref{item:convexity} holds there. From
$\Omega|TE^h_x = 0$ and $J(TE^h_x) = TE^h_x$ it follows that $\Omega(X,JX) \geq
0$ for all $X \in TE_x$, $x \in W$. Take $u \in \MM_J$ and consider the
function $h = e^\sigma \circ u$ on $U = u^{-1}(W)$. This satisfies $h|\partial
U \equiv e^{-\epsilon}$ and is subharmonic, because $d(dh \circ j) = d(\Theta
\circ J \circ Du \circ j) = -u^*\Omega \leq 0$. It follows that $h \leq
e^{-\epsilon}$ everywhere, which shows that $K = E \setminus int(W)$ has the
required property. \qed

The action of $u \in \BB$ is defined to be $A(u) = \int_S u^*\Omega$. This is
actually the same for all $u$, since
\begin{equation} \label{eq:action}
\int_S u^*\O = \int_{\partial S} u^*\Theta = \int_{\partial S} \kappa_Q.
\end{equation}
Therefore, if one equips $E$ with a symplectic form $\Omega + \pi^*\beta$
taming $J$, all the $u \in \MM_J$ become pseudo-holomorphic curves with the
same energy $\half \int_S ||Du||^2 = A(u) + \int_S \beta$.

\begin{lemma} \label{th:compactness}
$\MM_J$ is compact in any $C^r$-topology.
\end{lemma}

\proof We apply the Gromov compactness theorem for $(j,J)$-holomorphic maps $S
\rightarrow E$ with boundary in $Q$. The bubble components in the Gromov limit
appear through a reparametrization which ``magnifies'' successively smaller
parts of the domain. Since in our case all maps are sections, the bubbles are
either nonconstant $J$-holomorphic spheres in some fibre $E_z$, or nonconstant
$J$-holomorphic discs in $E_z$, $z \in \partial S$, with boundary on $Q_z$. But
both are excluded by our assumptions, since $\Omega$ is exact and $Q_z \subset
E_z$ an exact Lagrangian submanifold. \qed

A less formal version of the infinite-dimen\-sional framework introduced above
involves spaces of $W^{1,p}$-sections, $p>2$. We omit the construction itself,
and only mention its main consequence. For $(u,J) \in \MM^{univ}$, the
differential operator $D_{u,J}$ extends to a Fredholm operator $\WW^1_u
\rightarrow \WW^0_{u,J}$ from the $W^{1,p}$-completion $\WW^1_u$ of $T\BB_u$ to
the $L^p$-completion $\WW^0_{u,J}$ of $\EE_{u,J}$. We denote this extension
equally by $D_{u,J}$. If it is onto (in which case one says that $u$ is
regular), $\MM_J$ is a smooth finite-dimensional manifold near $u$. If that
holds for all $u \in \MM_J$, $J$ itself is called regular, and the space of
such $J$ is denoted by $\JJreg(E,\pi,Q,j) \subset \JJ(E,\pi,j)$.

\begin{lemma} \label{th:transversality}
$\JJ^{reg}(E,\pi,Q,j)$ is $\smooth$-dense in $\JJ(E,\pi,j)$. In fact the
following stronger statement holds: take some nonempty open subset $U \subset
S$ and a $J \in \JJ(E,\pi,j)$. Then there are $J' \in \JJ^{reg}(E,\pi,Q,j)$
arbitrarily close to $J$, such that $J = J'$ outside $\pi^{-1}(U)$.
\end{lemma}

\proof Even though this is a well-known argument, we recall part of it as a
preparation for subsequent more refined versions. Fix $U$ and $J$. By Lemmas
\ref{th:convexity} and \ref{th:compactness} there is an open subset $V \subset
E$ with $\overline{V} \cap (\partial_hE \cup E^\crit) = \emptyset$, such that
$\im(u) \subset V$ for each $u \in \MM_J$. A modified version of the
compactness argument shows that this property, with the same $V$, remains true
for all almost complex structures in $\JJ(E,\pi,j)$ which are sufficiently
close to $J$. We want to make $J$ regular by perturbing it on $V \cap
\pi^{-1}(U)$. Take $\TT \subset T\JJ(E,\pi,j)_J$ to be the subset of those $Y$
which vanish outside $V \cap \pi^{-1}(U)$, and consider the operator
\begin{equation} \label{eq:transversality-operator}
 D_{u,J}^{univ}: \WW^1_u \times \TT \longrightarrow \WW^0_{u,J}
\end{equation}
given by the same formula as in \eqref{eq:universal-operator}. By a standard
arguments surjectivity of this operator implies the desired result (strictly
speaking, what comes up in this argument is a dense subspace of those $Y$ which
have ``finite $\smooth_\epsilon$-norm'', but since the first component
$D_{u,J}$ is Fredholm, this makes no difference as far as surjectivity is
concerned). Let $F \rightarrow S$ be the bundle dual to
$\Lambda^{0,1}(u^*TE^v)$, so that $(\WW^0_{u,J})^* \iso L^q(F)$, $p^{-1} +
q^{-1} = 1$. Suppose that $\eta \in L^q(F)$ is orthogonal to the image of
\eqref{eq:transversality-operator}. It then satisfies
\begin{equation} \label{eq:adjoint}
 D_{u,J}^*\eta = 0 \; \text{ on $S \setminus \partial S$,} \quad
 \text{and} \quad \int_S \leftsc \eta, Y \circ Du \circ j \rightsc = 0
 \; \text{ for $Y \in \TT$.}
\end{equation}
The first equation implies that $\eta$ is smooth away from the boundary.
Suppose that $z \in U \setminus \partial S$ is a point where $\eta_z \neq 0$,
and set $x = u(z) \in V \cap \pi^{-1}(U)$. Take a $(j,J)$-antilinear map $Z:
TS_z \rightarrow TE^v_x$ such that $\leftsc \eta_z, Z \circ j \rightsc \neq 0$.
One can see from \eqref{eq:tangent-compatible} that there is a $Y \in
T\JJ(E,\pi,j)_J$ with $Y^{vh}_x = Z$, $Y^{vv}_x = 0$, and this will satisfy
$\leftsc \eta_z, (Y \circ Du \circ j)_z \rightsc \neq 0$. By multiplying $Y$
with a bump function supported near $x$, one can achieve that it lies in $\TT$
and that $\int_S \leftsc \eta, Y \circ Du \circ j \rightsc \neq 0$, a
contradiction. This means that $\eta\,|\, (U \setminus \partial S) = 0$. By
unique continuation $\eta\,|\,(S \setminus \partial S) = 0$, which proves that
$\eta = 0$. \qed

Take some $\zeta \in \partial S$ and consider the map $ev_\zeta: \BB
\rightarrow Q_\zeta$, $ev_\zeta(u) = u(\zeta)$. The next result belongs to a
type called ``transversality of evaluation''.

\begin{lemma} \label{th:ev-transversality}
Let $g$ be a smooth map from some arbitrary manifold $G$ to $Q_\zeta$. Then,
for any $J$ and $U$ as in the previous lemma, there are $J' \in
\JJreg(E,\pi,Q,j)$ arbitrarily close to $J$, with $J' = J$ outside
$\pi^{-1}(U)$, such that $ev_\zeta|\MM_{J'}$ is transverse to $g$.
\end{lemma}

\proof In the same setup as before, one now has to prove that for $u \in \MM_J$
and $x = u(\zeta) \in Q_\zeta$, the operator
\[
\WW^1_u \times \TT \longrightarrow \WW^0_{u,J} \times T(Q_\zeta)_x, \quad (X,Y)
\longmapsto (D^{univ}_{u,J}(X,Y),X_\zeta)
\]
is onto. Take $(\eta,\xi)$ orthogonal to the image, with $\xi \in
T(Q_\zeta)_x^\vee$. One still has \eqref{eq:adjoint} and as before it follows
that $\eta = 0$. Then $\leftsc \xi, X_{\zeta} \rightsc = 0$ for all $X \in
\WW^1_u$, so that $\xi = 0$ as well. \qed

For a given $(E,\pi,Q)$ and $\zeta$, one now proceeds as follows. After
choosing some $j$ and a $J \in \JJreg(E,\pi,Q,j)$, one obtains a smooth compact
moduli space $\MM_J$; and then one sets
\[
\Phi_1(E,\pi,Q) = (ev_\zeta)_*[\MM_J] \in H_*(Q_\zeta;\Z/2).
\]
This is independent of the choice of $j$, $J$ by a standard argument using
parame\-trized moduli spaces. The same reasoning shows that it remains
invariant under any ``smooth deformation'' of the geometric objects involved,
that is to say of $Q,\Omega,\Theta$ or of the fibration $\pi: E \rightarrow S$
itself, as long as one remains within the class of exact Lefschetz fibrations
with Lagrangian boundary conditions. It seems pointless to formalize this
notion of deformation, because the necessary conditions will be quite obviously
satisfied in all our applications.

\begin{remark} \label{re:bordism}
Since the regular spaces $\MM_J$ are actual manifolds, one can refine
$\Phi_1(E,\pi,Q)$ by regarding it as an element of the unoriented bordism group
$MO_*(Q_\zeta)$. But even with this refinement, it is far from capturing all
the information contained in $\MM_J$. To get more sophisticated invariants one
can use evaluation at several points, allowing those to move; we will now
explain the simplest version of this. Given $\zeta \in \partial S$, choose a
positively oriented path $c: [0;1] \rightarrow \partial S$, $c(0) = c(1) =
\zeta$, parametrizing the boundary component on which $\zeta$ lies. Let
$\phi_t: Q_{c(t)} \rightarrow Q_\zeta$ be the diffeomorphisms obtained from
symplectic parallel transport along $c|[t;1]$. Denote by $\Delta, \Gamma
\subset Q_\zeta^2$ the diagonal and the graph of $\phi_0$, respectively. For
regular $J$, the parametrized evaluation map
\[
\widetilde{ev}_\zeta: [0;1] \times \MM_J \rightarrow [0;1] \times Q_\zeta^2,
\quad \widetilde{ev}_\zeta(t,u) = (t,u(\zeta),\phi_t(u(c(t))))
\]
represents a class $\tilde{\Phi}_2(E,\pi,Q) \in H_*([0;1] \times Q_\zeta^2,
\{0\} \times \Gamma \cup \{1\} \times \Delta;\Z/2)$, and this (as well as its
cobordism version) is an invariant of $(E,\pi,Q)$. An example of this invariant
will be computed in Remark \ref{re:vanishing}(iii).
\end{remark}

Let $S^1,S^2$ be two compact surfaces with marked points $\zeta^k \in \partial
S^k$, and denote by $S = S^1 \#_{\zeta^1 \sim \zeta^2} S^2$ their boundary
connected sum (Figure \ref{fig:sum}). To be precise, one should choose oriented
embeddings $\psi^k: \bar{D}^+(1) \rightarrow S^k$ of the closed half-disc
$\bar{D}^+(1) = \{z \in \C \suchthat |z| \leq 1, \;\im\,z \geq 0\}$ into $S^k$,
such that $\psi^k(0) = \zeta^k$ and $(\psi^k)^{-1}(\partial S^k) = [-1;1]$.
Then, writing $D^+(\rho)$ for the open half-disc of some radius $0<\rho<1$, one
forms $S$ by taking the two $S^k \setminus \psi^k(D^+(\rho))$ and identifying
$\psi^1(z)$ with $\psi^2(-\rho^{-1}z)$. Suppose that we have exact Lefschetz
fibrations $(E^k,\pi^k)$ over $S^k$, such that their symplectic connections are
trivial on $\im(\psi^k)$; together with Lagrangian boundary conditions $Q^k
\subset E^k$, such that the one-forms $\kappa_{Q^k}$ vanish on
$\psi^k([-1;1])$; and finally, an exact syplectic manifold $M$ with an exact
Lagrangian submanifold $L$, and isomorphisms $\phi^k: M \rightarrow
(E^k)_{\zeta^k}$ satisfying $\phi^k(L) = (Q^k)_{\zeta^k}$. One can then glue
the $(E^k,\pi^k)$ and $Q^k$ to an exact Lefschetz fibration over $S$. To do
that, one first uses symplectic parallel transport to construct embeddings
\[
\xymatrix{
 {\bar{D}^+(1) \times M} \ar[r]^-{\Psi^k} \ar[d] &
 {E^k} \ar[d]^{\pi^k} \\
 {\bar{D}^+(1)} \ar[r]^-{\psi^k} & {S^k}
}
\]
such that $\Psi^k|\{0\} \times M = \phi^k$. These will satisfy
\[
(\Psi^k)^*\Omega^k = \o, \quad (\Psi^k)^*\Theta^k = \theta + dR^k, \quad
(\Psi^k)^{-1}(Q^k) = [-1;1] \times L
\]
for some functions $R^k$. As in the pasting construction described in Section
\ref{sec:basic}, one needs to introduce modified forms $\tilde{\Theta}^k$ such
that $(\Psi^k)^*\tilde{\Theta}^k = \theta$ near $\{0\} \times M$. The functions
$K_{Q^k}$ need to be modified accordingly, but that is rather straightforward,
so we will not write it down explicitly. Now take the $E^k \setminus
\Psi^k(D^+(\rho) \times M)$ and identify $\Psi^1(z,x)$ with $\Psi^2(-\rho
z^{-1},x)$. This yields a manifold $E$ with a map $\pi: E \rightarrow S$, and
the remaining data matches up, producing the structure of an exact Lefschetz
fibration together with a Lagrangian boundary condition $Q$.
\includefigure{sum}{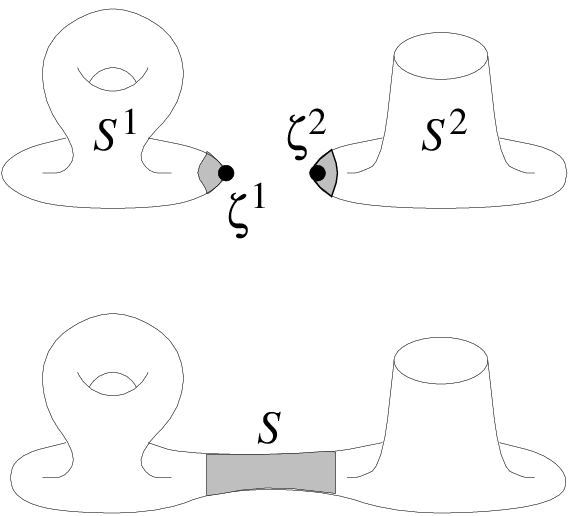}{hb}%

Assume that one has chosen complex structures $j^k$ on $S^k$ such that the
$\psi^k$ are holomorphic; these determine a complex structure $j$ on $S$. Take
$J^k \in \JJ(E^k,\pi^k,j^k)$ such that $(\Psi^k)^*(J^k)$ is the product of the
standard complex structure on $\bar{D}^+(1)$ and of some fixed $\o$-compatible
almost complex structure on $M$ (the same for both $k$). Then there is a
canonical induced $J \in \JJ(E,\pi,j)$. Note that even though we have
restricted the behaviour of $J^k$ over $\im(\Psi^k)$, it is still possible to
choose them regular, by using the more precise statement in Lemma
\ref{th:transversality}. Moreover, Lemma \ref{th:ev-transversality} says that
for suitably chosen $J^k$, the evaluation maps $ev_{\zeta^k} | \MM_{J^k}:
\MM_{J^k} \rightarrow Q^k_{\zeta_k} \iso L$ will be transverse to each other.

\begin{proposition} \label{th:gluing}
Assume that $J^k \in \JJreg(E^k,\pi^k,Q^k,j^k)$ for $k = 1,2$, and that the
$ev_{\zeta_k} | \MM_{J^k}$ are mutually transverse. Choose a sufficiently small
parameter $\rho$ for the gluing. Then $J \in \JJreg(E,\pi,Q,j)$, and there is a
diffeomorphism
\[
\MM_J \iso \MM_{J^1} \times_L \MM_{J^2},
\]
where the right hand side is the fibre product of
$(ev_{\zeta^1},ev_{\zeta^2})$. \qed
\end{proposition}

This is an average specimen of the ``gluing theorem'' type. The closest related
argument in the literature would seem to be the gluing theory for
pseudo-holomorphic discs from \cite[Section 18]{fukaya-oh-ohta-ono}, which is
far more sophisticated than what we need here; as an alternative, one can
probably adapt the proof of the more familiar gluing theorem for closed
pseudo-holomorphic curves \cite[Section 6]{ruan-tian94}, \cite[Appendix
A]{mcduff-salamon}, \cite{liu96}. The obvious next step would be to write down
the outcome as a ``gluing formula''. However, while gluing will be important
later, the situation then will be slightly different from that covered by
Proposition \ref{th:gluing}. For this reason, further discussion is postponed
to Section \ref{sec:relative-invariants}.

\subsection{Horizontality\label{sec:horizontal}}

Let $(E,\pi)$ be an exact Lefschetz fibration. Choose a complex structure $j$
on its base. $J \in \JJ(E,\pi,j)$ is called horizontal if $J_x(TE^h_x) =
TE_x^h$ for all $x \notin E^\crit$, or what is equivalent, if
$\Omega(\cdot,J\cdot)$ is symmetric. In terms of \eqref{eq:compatible-j} these
are just the $J$ with $J^{vh} = 0$, which shows that they form a contractible
subspace $\JJ^h(E,\pi,j) \subset \JJ(E,\pi,j)$. The importance of horizontal
almost complex structures is that they are sensitive to the geometry of the
symplectic connection; the following is a particularly simple instance of this.

\begin{lemma} \label{th:first-illu}
Let $(E,\pi)$ be an exact Lefschetz fibration over a compact surface $S$,
$\partial S \neq \emptyset$, with a Lagrangian boundary condition $Q$. Assume
that $(E,\pi)$ has nonnegative curvature, and that the boundary condition
satisfies $\int_{\partial S} \kappa_Q < 0$. Then $\Phi_1(E,\pi,Q) = 0$.
\end{lemma}

\proof For $J \in \JJ^h(E,\pi,j)$, take a symplectic form $\Omega + \pi^*\beta$
as in Lemma \ref{th:tame}. Then $J$ is compatible with it in the ordinary sense
of the word; we denote the associated metric by $||\cdot||$. Write $\Omega|TE^h
= f(\pi^*\beta|TE^h)$ with $f \in \smooth(E \setminus E^\crit,\R)$. For any map
$u$ from a compact Riemann surface, possibly with boundary, to $E$ there is the
familiar equality $\half \int ||Du||^2 = \int u^*(\Omega + \pi^*\beta) + \int
||\bar\partial_J u||^2$. Specialize to sections $u$ and split $Du = (Du)^h +
(Du)^v$ into horizontal and vertical parts. Since $||(Du)^h||^2 =
2(f(u)+1)\beta$, one obtains
\begin{equation} \label{eq:deficiency}
{\textstyle \half} \int_S ||(Du)^v||^2 + \int_S f(u)\beta = \int_S u^*\Omega +
\int_S ||\bar\partial_J u||^2.
\end{equation}
This implies that $\MM_J = \emptyset$. In fact, the curvature assumption is
just that $f \geq 0$, while for $u \in \MM_J$ one would have $\int_S u^*\Omega
< 0$ by \eqref{eq:action} and the second assumption. \qed

A section $u: S \rightarrow E$ is called horizontal if $Du_z(TS_z) =
(TE^h)_{u(z)}$ for all $z \in S$. To see the geometric meaning of this, it is
convenient to exclude temporarily the presence of critical points, so that
$(E,\pi)$ is an exact symplectic fibration.  If $u$ is horizontal, parallel
transport along any path $c: [a;b] \rightarrow S$ carries $u(c(a)) \in
E_{c(a)}$ to $u(c(b)) \in E_{c(b)}$. In other words, if $M$ is some fibre of
$E$ and $x \in M$ the unique point through which $u$ passes, the structure
group of the symplectic connection on $E$ is reduced from $\Sympe(M)$ to the
subgroup $\Sympe(M,x)$ of maps preserving $x$. This entails a restriction on
the curvature, namely, writing $\Omega|TE^h = f(\pi^*\beta|TE^h)$ for some
positive $\beta \in \Omega^2(S)$, one has
\begin{equation} \label{eq:critical-curvature}
d(f|E_z)_{u(z)} = 0 \quad \text{for all $z \in S$}.
\end{equation}
If $u$ is horizontal, the symplectic vector bundle $u^*TE^v \rightarrow S$ has
a preferred connection $\nabla^u$, obtained by linearizing parallel transport
around $u$. Equivalently, this is induced from the connection on $E$ by the
derivative map $\Sympe(M,x) \rightarrow Sp(TM_x)$. Explicitly $\nabla^u_ZX =
u^*([Z^h,\tilde{X}])$, where $\tilde{X}$ is any section of $TE^v$ with
$u^*\tilde{X} = X$. Using the canonical isomorphism $sp(V) \iso sym^2(V^*)$ for
any symplectic vector space $V$, one can write the curvature of $\nabla^u$ as a
two-form on $S$ with values in quadratic forms on the fibres of $u^*TE^v$. In
those terms it is given by $F_{\nabla^u} = Hess(f|E_z)_{u(z)} \beta$, which is
well-defined by \eqref{eq:critical-curvature}. We say that $\nabla^u$ is
nonnegatively curved if all these Hessians are $\geq 0$. The ``infinitesimal
deformations'' of a horizontal section $u$ are the covariantly constant
sections, $\nabla^u X = 0$. If such a section exists, it further reduces the
structure group of $(E,\pi)$ to the subgroup of maps in $\Sympe(M,x)$ which
preserve a certain tangent vector at $x$. The resulting curvature restriction
is
\begin{equation} \label{eq:hessians}
Hess(f)_{u(z)}(X(z),X(z)) = 0 \quad \text{for all $z \in S$.}
\end{equation}
Finally, if one readmits critical points, all the formulae derived above remain
valid, with essentially the same proofs, except that when talking about the
symplectic connection on $E$ it is necessary to restrict to a neighbourhood of
a fixed horizontal section.

The connection between the two notions which we have introduced so far is that
{\em a horizontal section is $(j,J)$-holomorphic for any horizontal $J$}. This
provides a useful class of pseudo-holomorphic sections with a geometric origin,
but it also raises a problem: if the space of horizontal sections has ``too
large dimension'', one cannot find an almost complex structure which is both
regular and horizontal. The rest of this section discusses this issue in more
detail. Assume from now on that $S$ is compact with $\partial S \neq
\emptyset$, and that we have a Lagrangian boundary condition $Q$. Write $\MM^h$
for the space of horizontal sections $u$ which lie in $\BB$, meaning that
$u(\partial S) \subset Q$. Since a horizontal section is determined by its
value at any point, evaluation at $\zeta \in \partial S$ identifies $\MM^h$
with a subset of $Q_\zeta$. As observed above, $\MM^h \subset \MM_J$ for all $J
\in \JJ^h(E,\pi,j)$.

\begin{lemma} \label{th:horizontal-transversality}
Let $U \subset S$ be a nonempty open subset, such that any partial section $u:
U \rightarrow E|U$ which is horizontal and satisfies $u(\partial S \cap U)
\subset Q$ is the restriction of some $u' \in \MM^h$. Then, given some $J \in
\JJ^h(E,\pi,j)$, there are $J' \in \JJ^h(E,\pi,j)$ arbitrarily close to it and
which agree with it outside $\pi^{-1}(U)$, with the property that any $u \in
\MM_{J'} \setminus \MM^h$ is regular.
\end{lemma}

\proof Take $V \subset E$ as in the proof of Lemma \ref{th:transversality}. In
parallel with the argument there, we have to show that for any $u \in \MM_J
\setminus \MM^h$,
\begin{equation} \label{eq:h-universal}
D_{u,J}^{univ}: \WW^1_u \times \TT^h \longrightarrow \WW^0_{u,J}
\end{equation}
is onto. Here $\TT^h = \TT \cap T\JJ^h(E,\pi,j)_J$ is the space of
infinitesimal deformations $Y$ of $J$ within the class of horizontal almost
complex structures, that is to say with $Y^{vh} = 0$ in
\eqref{eq:tangent-compatible}, and such that $Y = 0$ outside $V \cap
\pi^{-1}(U)$. As before, an $\eta$ which is orthogonal to the image of
\eqref{eq:h-universal} is smooth away from the boundary and satisfies
\begin{equation} \label{eq:h-orthogonal}
 \int_S \leftsc \eta, Y \circ Du \circ j \rightsc = 0
 \quad \text{for $Y \in \TT^h.$}
\end{equation}
Suppose that $u|U$ is horizontal. By assumption one could then find a $u' \in
\MM^h$ with $u'|U = u|U$; unique continuation for pseudo-holomorphic curves
would imply that $u' = u$, a contradiction. Hence there is a $z \in U$ such
that $Du_z(TS_z) \neq (TE^h)_{u(z)}$. This is an open condition, so we can
assume that $z \in U \setminus \partial S$. By choosing $Y$ suitably in
$\TT^h$, one can make $(Y \circ Du \circ j)_z$ equal to any arbitrary
$(j,J)$-antilinear homomorphism $TS_z \rightarrow (TE^v)_{u(z)}$. In
particular, if $\eta_z \neq 0$ one can achieve that $\leftsc \eta_z, (Y \circ
Du \circ j)_z \rightsc \neq 0$, which after multiplying with a cutoff function
leads to a contradiction with \eqref{eq:h-orthogonal}. The same argument
applies to all points close to $z$, proving that $\eta$ vanishes on a nonempty
open subset, from which it follows that $\eta = 0$. \qed

In particular, taking $U = S$ shows that if $\MM^h = \emptyset$ then
$\JJ^{reg,h}(E,\pi,Q,j) = \JJ^{reg}(E,\pi,Q,j) \cap \JJ^h(E,\pi,j) \subset
\JJ^h(E,\pi,j)$ is a dense subset. In the same way one proves the following
analogue of Lemma \ref{th:ev-transversality}:

\begin{lemma} \label{th:horizontal-ev-transversality}
Let $g$ be a smooth map from an arbitrary manifold $G$ to $Q_\zeta$, for some
$\zeta \in \partial S$. Let $U \subset S$ be a nonempty open subset, such that
there are no horizontal partial sections $u: U \rightarrow E|U$ with
$u(\partial S \cap U) \subset Q$. Given $J \in \JJ^h(E,\pi,j)$, there are $J'
\in \JJ^{reg,h}(E,\pi,Q,j)$ arbitrarily close to it and which agree with it
outside $\pi^{-1}(U)$, such that $ev_\zeta | \MM_{J'}$ is transverse to $g$.
\qed
\end{lemma}

Let $T_{Zar}(\MM^h)_u$ be the Zariski tangent space of $\MM^h$ at some section
$u$. It consists of those $X \in \smooth(u^*TE^v)$ that satisfy $\nabla^uX = 0$
and lie in $T\BB_u$, meaning that $X|\partial S \subset u^*(TQ \cap TE^v)$. The
next result helps to determine when $u$, considered as a pseudo-holomorphic
section for some horizontal almost complex structure, is regular.

\begin{lemma} \label{th:flat-deformations}
Take $J \in \JJ^h(E,\pi,j)$ and $u \in \MM^h$. Then $T_{Zar}(\MM^h)_u \subset
\ker\, D_{u,J}$, and if moreover $\nabla^u$ is nonnegatively curved, the two
spaces are equal.
\end{lemma}

\proof We take the second derivative of \eqref{eq:deficiency} at $u$. This is
well-defined because the action $A$ is constant on $\BB$ and the other terms
vanish to first order at $u$. The outcome is that for $X \in T\BB_u$,
\begin{equation} \label{eq:weitzenboeck}
\int_S ||\nabla^u X||^2 + \int_S Hess(f \circ u)(X,X) \,\beta = 2 \int_S
||D_{u,J}X||^2.
\end{equation}
If $\nabla^uX = 0$ then $Hess(f \circ u)(X,X) = 0$ by \eqref{eq:hessians},
which implies that $D_{u,J}X = 0$. If $\nabla^u$ has nonnegative curvature, the
converse also holds, since the second term in \eqref{eq:weitzenboeck} is $\geq
0$. \qed

Call $\MM^h$ clean if it is a smooth manifold and its tangent space is
everywhere equal to $T_{Zar}\MM^h$. The next result can be considered as a
limiting case of Lemma \ref{th:first-illu} (and is again just a sample
application, in itself without any great importance, but hopefully
instructive).

\begin{lemma}
Assume that $(E,\pi)$ has nonnegative curvature, and that the Lagrangian
boundary condition $Q$ satisfies $\int_{\partial S} \kappa_Q = 0$. In addition,
assume that $\MM^h$ is clean and that its dimension agrees at every point with
the index of $D_{u,J}$. Then $\Phi_1(E,\pi,Q) = [ev_{\zeta}(\MM^h)]$.
\end{lemma}

\proof Take $J \in \JJ^h(E,\pi,j)$, and consider \eqref{eq:deficiency}. Since
$f \geq 0$ and $\int_S u^*\Omega = 0$, it follows that $\MM_J = \MM^h$, and
moreover that every $u \in \MM^h$ satisfies $f(u) \equiv 0$, which in turn
implies $Hess(f|E_z)_{u(z)} \geq 0$. Applying Lemma \ref{th:flat-deformations}
shows that $\dim\, \ker\, D_{u,J} = \dim\,T_{Zar}\MM^h = \dim\,\MM^h =
\ind\,D_{u,J}$, from which one sees that $\coker\, D_{u,J} = 0$; hence $J$ is
regular. \qed

\subsection{A vanishing theorem\label{sec:vanishing}}

Let $L$ be a framed exact Lagrangian sphere in an exact symplectic manifold
$M$. According to Proposition \ref{th:standard-fibrations} one can associate to
it a standard fibration $(E^L,\pi^L)$ over $\bar{D}(r)$ for some $r>0$. Each
fibre $E^L_z$, $z \neq 0$, contains a distinguished Lagrangian sphere
$\Sigma_z^L$, described by Lemma
\ref{th:standard-properties}\ref{item:lagrangian-subbundle}; for $z = r$, the
isomorphism $\phi^L: E^L_r \rightarrow M$ takes $\Sigma_r^L$ to $L$. Define the
standard Lagrangian boundary condition for $(E^L,\pi^L)$ to be
\begin{equation} \label{eq:standard-boundary}
\begin{split}
& \textstyle Q^L = \textstyle\bigcup_{z \in \partial\bar{D}(r)} \Sigma_z^L, \\
& \kappa_{Q^L} =d^c(-\quarter|z|)\,|\,\partial\bar{D}(r), \qquad K_{Q^L} = 0.
\end{split}
\end{equation}
This is indeed a Lagrangian boundary condition, as one can see from
\eqref{eq:theta-sigma-l}. Our aim here is to prove the following result about
the invariant $\Phi_1(E^L,\pi^L,Q^L) \in H_*(\Sigma_r^L;\Z/2) \iso
H_*(L;\Z/2)$:

\begin{prop} \label{th:vanishing}
$\Phi_1(E^L,\pi^L,Q^L) = 0$ for all $M$ and $L$.
\end{prop}

The first part of the proof is a degeneration argument, in which one restricts
the base to successively smaller discs. For $0< s \leq r$ set $Q^{L,s} =
\textstyle \bigcup_{z \in \partial\bar{D}(s)} \Sigma^L_z$; together with
$\kappa_{Q^{L,s}}$ and $K_{Q^{L,s}}$ as before, it is a Lagrangian boundary
condition for $(E^{L,s},\pi^{L,s}) = (E^L,\pi^L)\,|\,\bar{D}(s)$. This
constitutes a deformation of $(E^L,\pi^L,Q^L)$ in a suitable sense; which means
that if one identifies $\Sigma_r^L$ with $\Sigma_s^L$ using parallel transport
along $[s;r]$, then
\begin{equation} \label{eq:deformation-equation}
\Phi_1(E^L,\pi^L,Q^L) = \Phi_1(E^{L,s},\pi^{L,s},Q^{L,s})
\end{equation}
for all $s$. Now fix some $J \in \JJ^h(E^L,\pi^L,j)$, where $j$ is the standard
complex structure on $\bar{D}(r)$. Let $\MM^{L,s}$ be the space of sections
$\bar{D}(s) \rightarrow E^{L,s}$ which are holomorphic with respect to
$j|\bar{D}(s)$ and $J|E^{L,s}$, and have boundary in $Q^{L,s}$. As $s
\rightarrow 0$, $Q^{L,s}$ shrinks to the single critical point $\Sigma^L_0 =
\{x_0\}$ of $\pi^L$, and we would like to apply a compactness argument to
elements of $\MM^{L,s}$ in the limit. This looks a bit unpleasant as it stands,
but one can modify the situation to make the degeneration of the $Q^{L,s}$ less
singular, and then standard Gromov compactness is sufficient.

\begin{lemma} \label{th:resolution}
For some $0<r'\leq r$, there is a compact almost complex manifold with corners
$(\hatE^L,\hatJ)$, together with pseudo-holomorphic maps
\[
\xymatrix{
 {\hatE^L} \ar[r]^{\eta^L} \ar[d]_{\hat\pi^L} & {E^L} \ar[d]^{\pi^L} \\
 {\bar{D}(r') \times \bar{D}(1)} \ar[r]^-{m} & {\bar{D}(r)}
}
\]
where $m(w,z) = wz$ is multiplication, such that the following properties are
satisfied.
\begin{romanlist}
\item \label{item:resolution-one}
Away from $Z^L = (\eta^L)^{-1}(x_0) \cap (\hat\pi^L)^{-1}(0,0)$, the map
$\eta^L$ identifies $(\hatE^L,\hat\pi^L)$ with the pullback of $(E^L,\pi^L)$ by
$m$. In particular, restricting to any $w = s>0$ gives a pseudo-holomorphic
diffeomorphism
\[
 \hatE^{L,s} = (\hat\pi^L)^{-1}(\{s\} \times \bar{D}(1))
 \longrightarrow E^{L,s}.
\]
\item \label{item:resolution-two}
$\hatE^L$ carries a symplectic form $\hatO^L$ which tames its almost complex
structure, and which is of the form
\[
 \hatO^L = (\eta^L)^*\Omega^L + (\hat\pi^L)^*(
 \textstyle\frac{i}{2} dw \wedge d\bar{w} +
 \frac{i}{2} dz \wedge d\bar{z} + d\gamma) + \delta,
\]
with $\gamma = \frac{1}{4} \im((r'-w)\bar{z}dz)$, and where $\delta$ is
supported in a small neighbourhood of $Z^L$.

\item \label{item:resolution-three}
$(\hatE^L,\hatO^L)$ contains a Lagrangian submanifold with boundary $\hatQ^L$
which projects to $[0;r'] \times S^1 \subset \bar{D}(r') \times \bar{D}(1)$ and
satisfies
\[
\eta^L(\hatQ^L \cap (\hat\pi^L)^{-1}(\{s\} \times S^1)) = \begin{cases} Q^{L,s}
& s>0,
\\ \{x_0\} & s = 0.
\end{cases}
\]
\end{romanlist}
\end{lemma}

It is convenient to begin with the local model for the construction. Choose
some $r'>0$ and set
\begin{equation} \label{eq:cylinder}
\begin{split}
 & E = \{x \in \C^{n+1} \suchthat ||x|| \leq 2\sqrt{r'},\; |q(x)| \leq r'\},
 \qquad \Omega = \o_{\C^{n+1}}|E, \\
 & \pi = q: E \longrightarrow \bar{D}(r'),
 \qquad Q = \textstyle \bigcup_{z \in \partial \bar{D}(r')} \Sigma_z
\end{split}
\end{equation}
where $q$ is our standard quadratic function, and $\Sigma_z$ is as in
\eqref{eq:model-spheres}. In parallel with the reasoning above, one can also
introduce $(E^s,\pi^s) = (E,\pi)|\bar{D}(s)$ and $Q^s = \bigcup_{z \in \partial
\bar{D}(s)} \Sigma_z$ for $s \in (0;r']$. Take the pullback of $(E,\pi)$ by the
multiplication map, $m^*E = \{(w,z,x) \in \bar{D}(r') \times \bar{D}(1) \times
E \suchthat q(x) = wz\} \subset \C^{n+3}$. This has one singular point
$(0,0,0)$, which can be resolved by blowing it up inside $\C^{n+3}$ and taking
the proper transform of $m^*E$. The outcome, which we denote by $\hatE$, comes
with maps
\[
\xymatrix{
 {\hatE} \ar[r]^{\eta} \ar[d]_{\hat\pi} & {E} \ar[d]^{\pi} \\
 {\bar{D}(r') \times \bar{D}(1)} \ar[r]^-{m} & {\bar{D}(r')}.
}
\]
\begin{condensedprimelist}
\item \label{item:model-resolution-one} \em
$\eta$ is a pullback map away from $Z = \eta^{-1}(0) \cap \hat\pi^{-1}(0,0)$.
\end{condensedprimelist}
That is obvious from the definition. Consider the two-form
\begin{equation} \label{eq:pre-form} \Omega + \textstyle
\frac{i}{2} dw \wedge d\bar{w} + \frac{i}{2} dz \wedge d\bar{z} + d\gamma =
\o_{\C^{n+3}} + d\gamma
\end{equation}
on $\bar{D}(r') \times \bar{D}(1) \times E$, where $\gamma$ is as in Lemma
\ref{th:resolution}\ref{item:resolution-two}. We claim that this is symplectic
and tames the obvious complex structure. All one needs to verify is that
$(d\gamma)(\cdot,i\cdot)$ is nonnegative, which one can do by decomposing
$d\gamma = (d\gamma)^{1,1} + (d\gamma)^{0,2} + (d\gamma)^{2,0}$; then
$(d\gamma)^{1,1} = \re(r'-w)\frac{i}{4} dz \wedge d\bar{z}$ is nonnegative
since $\re(r'-w) \geq 0$, while $(d\gamma)^{0,2}(\cdot,i\cdot)$ vanishes
because it is a two-form of type $(0,2)$, and similarly for $(d\gamma)^{2,0}$.
Now pull back \eqref{eq:pre-form} to the blowup, make it symplectic by adding a
two-form $\delta$ supported near the exceptional divisor, and restrict that to
$\hatE$. The outcome is:

\begin{condensedprimelist}
\setcounter{enumi}{1}
\item \label{item:model-resolution-two} \em
There is a symplectic form $\hatO$ on $\hatE$ which tames the complex
structure, and which is of the form $\eta^*\Omega + (\hat\pi)^*(\frac{i}{2} dw
\wedge d\bar{w} + \frac{i}{2} dz \wedge d\bar{z} + d\gamma) + \delta$, with
$\delta$ supported in a small neighbourhood of $Z$.
\end{condensedprimelist}

The subset
\begin{equation} \label{eq:hatq}
 \{ (w,z,x) \suchthat w \in [0;r'], \; |z| = 1, \; x \in
 \Sigma_{wz}\} \subset m^*E
\end{equation}
is a submanifold, since it is the image of the embedding $S^1 \times_{\Z/2}
\bar{B}^{n+1}(1) \rightarrow \C^{n+3}$, $(z,y) \mapsto (|y|^2,z^2,zy)$. It is
also Lagrangian with respect to \eqref{eq:pre-form}; which is in fact the
reason for our choice of symplectic form. Define $\hatQ$ to be the preimage of
\eqref{eq:hatq} in $\hatE$; provided that the support of $\delta$ has been
chosen sufficiently small, this is again a Lagrangian submanifold. By
construction one has

\begin{condensedprimelist}
\setcounter{enumi}{2}
\item \label{item:model-resolution-three} \em
$\hat\pi(\hatQ) = [0;r] \times S^1$, and
\[
\eta(\hatQ \cap \hat\pi^{-1}(\{s\} \times S^1)) = \begin{cases} Q^s & s>0,
\\ \{0\} & s = 0.
\end{cases}
\]
\end{condensedprimelist}

\proof[Proof of Lemma \ref{th:resolution}] Take a holomorphic Morse chart
$(\xi,\Xi)$ around $x_0 \in E_0^L$ as provided by Lemma
\ref{th:standard-properties}\ref{item:preferred-charts}. After restricting the
domains, we may assume that $\xi$ is the inclusion $\bar{D}(r') \hookrightarrow
\bar{D}(r)$ for some $0<r' \leq r$, that $\Xi$ is defined on the set $E \subset
\C^{n+1}$ from \eqref{eq:cylinder}, and that $Q^{L,s} = \Xi(Q^s)$ for all $0 <
s \leq r'$. Construct a diagram
\[
\xymatrix{
 {\hatE} \ar[d] \ar[r] \ar@/^1pc/[rr]^{\eta^L} &
 {m^*E} \ar[d]^{m^*\Xi} \ar[r] & {E} \ar[d]^{\Xi} \\
 {\hatE^L} \ar[dr]_{\hat\pi^L} \ar[r] &
 {m^*E^L} \ar[r] \ar[d] &
 {E^{L,r'}} \ar[d]^{\pi^{L,r'}} \\
 & {\bar{D}(r') \times \bar{D}(1)} \ar[r]^{m} &
 {\bar{D}(r')}
 }
\]
as follows: pull back $E^L$ by multiplication $m$; the map $m^*\Xi$ induced by
$\Xi$ identifies neighbourhoods of the singular points $(0,0,0) \in m^*E$ and
$(0,0,x_0) \in m^*E^L$; which means that there is a resolution \[\hatE^L
\longrightarrow m^*E^L\] modelled locally on $\hatE \rightarrow m^*E$. The
given almost complex structure $J$ induces one on $m^*E$, for which $m^*\Xi$ is
holomorphic; this and the complex structure on $\hatE$ define the almost
complex structure $\hatJ$. Part \ref{item:resolution-one} of the lemma now
follows from the corresponding statement \ref{item:model-resolution-one} in the
local model. Next, take the two-form $\Omega^L + \frac{i}{2} dw \wedge d\bar{w}
+ \frac{i}{2} dz \wedge d\bar{z} + d\gamma$ on $m^*E^L$; due to the nonnegative
curvature of standard fibrations, see Lemma
\ref{th:standard-properties}\ref{item:nonnegatively-curved}, and to the fact
that $J$ is horizontal, this will tame the almost complex structure away from
the singular point. One gets $\hatO^L$ by gluing this together with the
symplectic form $\hatO$ from the local model, and that proves
\ref{item:resolution-two}. Similarly, the definition of $\hatQ^L$ follows that
of $\hatQ$, and \ref{item:model-resolution-three} implies
\ref{item:resolution-three}. \qed

Returning to the moduli spaces $\MM^{L,s}$ of pseudo-holomorphic sections, we
can now carry out the compactness argument mentioned above:

\begin{lemma} \label{th:degeneration}
Choose some neighbourhood of the critical point $x_0 \in E^L$. For sufficiently
small $s$, the image of all $u \in \MM^{L,s}$ will lie in that neighbourhood.
\end{lemma}

\proof Let $(s_k)$ be a sequence in $(0;r']$ converging to zero, and $u_k \in
\MM^{L,s_k}$ a corresponding sequence of sections. By Lemma
\ref{th:resolution}\ref{item:resolution-one}, there is for each $k$ a unique
pseudo-holomorphic map \[\hat{u}_k: \bar{D}(1) \longrightarrow \hatE^L\] with
$\eta^L(\hat{u}_k(z)) = u_k(s_kz)$ and $\hat{\pi}^L(\hat{u}_k(z)) = (s_k,z)$.
By part \ref{item:resolution-three} of the same lemma, this maps $S^1$ to
$\hatQ^L$. Moreover, its energy is independent of $k$, since
\begin{align*}
 & \int_{\bar{D}(1)} \hat{u}_k^*\hatO^L
 = \int_{\bar{D}(s_k)} u_k^*\O^L + \int_{\bar{D}(1)}
 {\textstyle \frac{i}{2} dz \wedge d\bar{z}} + \int_{\bar{D}(1)}
 {\textstyle \frac{i}{4} (r'-s_k) \im(d\bar{z} \wedge dz)} +
 \\ & \qquad + \int_{\bar{D}(1)}
 \hat{u}^*_k\delta
 = {\textstyle \frac{\pi}{2}} s_k + \pi +
 {\textstyle \frac{\pi}{2}}(r'-s_k) + 0 = \pi(1 + \textstyle\frac{r'}{2}).
\end{align*}
Here we have used \eqref{eq:theta-sigma-l} for the first term; and the last
term is zero because, without changing its cohomology class, one can modify
$\delta$ to make its support arbitrarily close to $Z$, in which case, since
$u_k$ avoids the critical point $x_0$, $\hat{u}_k$ would not meet the support
of $\delta$. One can apply Gromov compactness to the sequence $\hat{u}_k:
(\bar{D}(1), S^1) \rightarrow (\hatE^L,\hatQ^L)$; the limit of some subsequence
will be a pseudo-holomorphic ``cusp disc'' or ``stable disc'' $\hat{u}_\infty$.
Since $\im(\hat{u}_k) \subset (\hat\pi^L)^{-1}(\{s_k\} \times \bar{D}(1))$, we
have $\im(\hat{u}_\infty) \subset (\hat\pi^L)^{-1}(\{0\} \times \bar{D}(1))$.
Composing with $\eta^L$ yields a ``cusp disc'' $u_\infty$ in $E$ whose boundary
lies in $\eta^L(\hatQ^L \cap (\hat{\pi}^L)^{-1}(\{0\} \times S^1)) = \{x_0\}$.
Because $\Omega^L$ is exact, $u_\infty$ is necessarily constant equal to $x_0$.
On the other hand, the image of $u_k$, for $k$ large and in our subsequence,
lies in an arbitrarily small neighbourhood of the image of $u_\infty$. The rest
is straightforward. \qed

The second part of the proof of Proposition \ref{th:vanishing} is an explicit
computation in the local model $(E,\pi)$ from \eqref{eq:cylinder}. Note that
while this is not an exact Lefschetz fibration, it still makes sense to
consider the spaces $\MM^s$ of holomorphic sections $w: \bar{D}(s) \rightarrow
E^s$ with boundary in $Q^s$. In fact $Q^s$ is a Lagrangian submanifold of
$\C^{n+1}$, as one can see from \eqref{eq:theta-sigma}, and the maximum
principle ensures that any holomorphic disc in $\C^{n+1}$ with boundary in
$Q^s$ is actually contained in $E^s$.

\begin{lemma} \label{th:explicit-sections}
$\MM^s$ consists of maps $w(z) = s^{-1/2} az + s^{1/2} \bar{a}$, where $a \in
\C^{n+1}$ satisfies $q(a) = 0$ and $||a||^2 = 1/2$. All these maps are regular,
in the sense that the associated Fredholm operators are surjective.
\end{lemma}

\proof The holomorphic functions $v: \bar{D}(s) \rightarrow \C$ which satisfy
$v(z) \in z^{1/2}\R$ for all $z \in \partial \bar{D}(s)$ are $v(z) = s^{-1/2} c
z + s^{1/2} \bar{c}$, for $c \in \C$. All components of $w \in \MM^s$ must be
of this form, and the conditions on $a$ come from $q(w(z)) = z$. As for
regularity, since we are dealing with the standard complex structure $J_0$ on
$\C^{n+1}$, $D_{w,J_0}$ is an actual $\bar\partial$-operator, see
\eqref{eq:connection-formula}. Its kernel consists of holomorphic maps $X:
\bar{D}(s) \rightarrow \C^{n+1}$ such that $X(z) \in \sqrt{z}\R^n$ for $z \in
\partial \bar{D}(s)$, and $Dq(w(z))X(z) = 2 \sum_k w_k(z) X_k(z) = 0$ for all
$z$. The same argument as before determines all such $X$ explicitly, the
outcome being that $\ker\,D_{w,J_0} \iso \R^{2n-1}$. Using the Riemann-Roch
formula for surfaces with boundary one computes that $\ind\,D_{w,J_0} = 2n-1$,
which shows that the cokernel is zero. \qed

The condition on $a$ can be written as
\[
||\re\,a||^2 = ||\im\,a||^2 = 1/4, \quad \leftsc \re\,a,  \im\,a \rightsc = 0
\]
so that $\MM^s$ can be identified with the sphere bundle $S(T^*S^n)$ by mapping
$a$ to $(u,v) = (-2\, \im\,a, 2\, \re\,a)$. With this and the diffeomorphism
$S^n \rightarrow \Sigma_s$, $x \mapsto \sqrt{s}x$, one can identify the
evaluation $ev_s: \MM^s \rightarrow \Sigma_s$ with the projection to the base
$S(T^*S^n) \rightarrow S^n$. This clearly represents the zero cycle in
$H_*(S^n;\Z/2)$. Hence, if one defined an invariant $\Phi_1(E^s,\pi^s,Q^s)$ in
this local model, it would vanish.

\proof[Proof of Proposition \ref{th:vanishing}] As in the proof of Lemma
\ref{th:degeneration}, we consider a holomorphic Morse chart $(\xi,\Xi)$ with
$\Xi$ defined on $E \subset \C^{n+1}$ for some $r'>0$. We know that for some
sufficiently small $s>0$, all $u \in \MM^{L,s}$ lie in $\im(\Xi)$. Recall that
$\Xi$ is holomorphic with respect to the standard complex structure $J_0$ on
$E$ and to $J$ on $E^L$; and that it takes $Q^s$ to $Q^{L,s}$. This implies
that composition with $\Xi$ yields a diffeomorphism
\[
\xymatrix{
 {\MM^s} \ar[r]^{\iso} \ar[d]_{ev_s} & {\MM^{L,s}} \ar[d]_{ev_s} \\
 {\Sigma_s} \ar[r]^{\Xi} & {\Sigma^L_s.}
}
\]
It is easy to see that, for $w \in \MM^s$ and $u = \Xi \circ w \in \MM^{L,s}$,
the kernels and cokernels of the associated operators $D_{w,J_0}$, $D_{u,J}$
coincide. This shows that $J|E^{L,s}$ is regular, and that the evaluation cycle
$\MM^{L,s} \rightarrow \Sigma^L_s$ can be identified with projection $S(T^*S^n)
\rightarrow S^n$, which as observed before is zero in homology; together with
\eqref{eq:deformation-equation} this completes the argument. \qed

\begin{remarks} \label{re:vanishing}
(i) As pointed out in Remark \ref{re:bordism}, $\Phi_1$ can be refined to an
invariant taking values in unoriented bordism. Proposition \ref{th:vanishing}
also applies to this refinement, since \eqref{eq:deformation-equation} comes
from a cobordism of moduli spaces and the projection $S(T^*S^n) \rightarrow
S^n$ is obviously cobordant to zero.

(ii) Suppose that $(E^1,\pi^1,Q^1)$ is a standard fibration with its standard
boundary condition, and that we have another exact Lefschetz fibration
$(E^2,\pi^2,Q^2)$, such that the two can be glued together to $(E,\pi,Q)$ as in
Section \ref{sec:simple-invariant}. Proposition \ref{th:vanishing} together
with Proposition \ref{th:gluing} implies that $\Phi_1(E,\pi,Q)$ will then be
zero. By pushing this reasoning a little further, one arrives at the following
result: {\em let $(E,\pi)$ be an arbitrary exact Lefschetz fibration over a
compact base $S$, with Lagrangian boundary condition $Q$. Assume that there is
a path $c: [a;b] \rightarrow S$ with $c(a) \in \partial S$, $c^{-1}(S^\crit) =
\{b\}$, and $c'(b) \neq 0$, whose vanishing cycle $V_{c(a)}$ is isotopic to
$Q_{c(a)}$ as an exact Lagrangian submanifold in $E_{c(a)}$; then
$\Phi_1(E,\pi,Q) = 0$}.

(iii) Despite the vanishing of their $\Phi_1$-invariant, the moduli spaces of
pseudo-holomorphic sections of a standard fibration are never empty. To see
this, one has to use the invariant $\tilde{\Phi}_2$ mentioned in Remark
\ref{re:bordism}. The computation of this can be reduced to the local model in
the same way as before, using Lemma \ref{th:degeneration}. The relevant
parametrized evaluation map then becomes
\begin{equation} \label{eq:para-e}
\begin{split}
 & [0;1] \times \MM^s \longrightarrow [0;1] \times
\Sigma_s \times \Sigma_s, \\
 & (t,w) \longmapsto (t,w(s),e^{\pi i(1-t)} w(s e^{2\pi i t})).
\end{split}
\end{equation}
Here $\Sigma_{s e^{2\pi i t}} \rightarrow \Sigma_s$, $x \mapsto e^{\pi
i(1-t)}x$ arises as parallel transport along $\partial\bar{D}(s)$ in positive
direction. Using Lemma \ref{th:explicit-sections} one can make
\eqref{eq:para-e} even more concrete, identifying it with
\begin{align*}
 & [0;1] \times S(T^*S^n) \longrightarrow [0;1] \times S^n \times S^n, \\
 & (t,u,v) \longmapsto (t,v,-cos(\pi t)v - \sin(\pi t)u).
\end{align*}
It is now easy to see that $\tilde{\Phi}_2(E^L,\pi^L,Q^L)$ is nonzero: it is
the image of the fundamental class $[S^n \times S^n]$ under the map
\begin{multline*}
H_*(S^n \times S^n;\Z/2) \rightarrow H_*(S^n \times S^n, \Delta \cup
\overline{\Delta};\Z/2) \iso \\ \iso H_*([0;1] \times S^n \times S^n, \{0\}
\times \overline{\Delta} \cup \{1\} \times \Delta;\Z/2)
\end{multline*}
where $\Delta$, $\overline{\Delta}$ are the diagonal and antidiagonal,
respectively.
\end{remarks}

\subsection{Relative invariants\label{sec:relative-invariants}}

A surface with strip-like ends is an oriented connected surface $S$, together
with finite sets $I^-,I^+$ and oriented proper embeddings $\{\gamma_e: \R^-
\times [0;1] \rightarrow S\}_{e \in I^-}$, $\{\gamma_e: \R^+ \times [0;1]
\rightarrow S\}_{e \in I^+}$, such that $\gamma_e^{-1}(\partial S) = \R^{\pm}
\times \{0;1\}$. The images of the $\gamma_e$ (the ends of $S$) should be
mutually disjoint, and the complement of the union of all ends should be a
relatively compact subset of $S$. We will always assume that there is at least
one end.

Define an {\em exact Lefschetz fibration trivial over the ends of $S$} to be an
exact Lefschetz fibration $(E,\pi)$, whose regular fibres are isomorphic to
some exact symplectic manifold $M$, together with smooth trivializations
$\{\Gamma_e: \R^- \times [0;1] \times M \rightarrow \gamma_e^*E\}_{e \in I^-}$,
$\{\Gamma_e: \R^+ \times [0;1] \times M \rightarrow \gamma_e^*E\}_{e \in I^+}$,
such that $\Gamma_e^*\Omega$, $\Gamma_e^*\Theta$ are equal to the pullbacks of
$\o$,$\theta$ by the projection $\R^{\pm} \times [0;1] \times M \rightarrow M$.
When considering a Lagrangian boundary condition $Q$ for such an exact
Lefschetz fibration, we will always impose the additional condition that
$\kappa_Q$ vanishes on $\im(\gamma_e) \cap \partial S$. To see the significance
of this, recall that the family $Q_z \subset E_z$ of Lagrangian submanifolds is
preserved by symplectic parallel transport along $\partial S$. Since the
symplectic connection is trivial on the ends, it follows that
$
 \Gamma_e^{-1}(Q) =
 (\R^\pm \times \{0\} \times L_{e,0}) \cup
 (\R^\pm \times \{1\} \times L_{e,1})
$
for some pair of Lagrangian submanifolds $L_{e,0},L_{e,1} \subset M$; we say
that $Q$ is modelled on $(L_{e,0},L_{e,1})$ over the end $e$. Using the
definition of a Lagrangian boundary condition, one finds that
\[
 d(K_Q \circ \Gamma_e) \,|\, \R^{\pm} \times \{k\} \times L_{e,k} =
 (\theta | L_{e,k}) - (\gamma_e^*\kappa_Q | \R^{\pm} \times \{k\}).
\]
Our assumption is that the second term on the right vanishes, which implies
that $K_Q(\Gamma_e(s,k,y))$ is independent of $s$. It follows that $L_{e,k}$ is
an exact Lagrangian submanifold in a canonical way, with associated function
$K_{L_{e,k}}(y) = K_Q(\Gamma_e(s,k,y))$.

In this situation, and under the additional assumption that the intersections
$L_{e,0} \cap L_{e,1}$ are transverse for all $e$, we will associate to
$(E,\pi)$ and $Q$ a relative invariant, which is a map between Floer cohomology
groups
\[
 \Phi_0^{rel}(E,\pi,Q) : \bigotimes_{e \in I^+} HF(L_{e,0},L_{e,1})
 \longrightarrow \bigotimes_{e \in I^-} HF(L_{e,0},L_{e,1}).
\]
This is a modified version of the frameworks described in \cite{schwarz95},
\cite{piunikhin-salamon-schwarz94} and \cite{seidel97}. The fact that
$\Phi^{rel}_0$ goes from positive to negative ends has to do with the use of
Floer {\em co}homology, which we think of as behaving contravariantly. This is
of course largely a matter of convention. In any case, there is no real
difference between positive and negative ends, as one can be turned into the
other by switching from $\gamma_e(s,t)$ to $\gamma_e(-s,1-t)$. Doing that does
not change the relative invariant, up to the ``Poincar{\'e} duality''
isomorphism $HF(L_{e,0},L_{e,1}) \iso HF(L_{e,1},L_{e,0})^\vee$. One could
therefore formulate the theory using only one kind of ends, bringing it closer
to a (1+1)-dimensional TQFT. Finally, we should mention that the transverse
intersection condition can be lifted. This goes by a standard argument, using
the same ``continuation map'' technique as the proof of isotopy invariance of
Floer cohomology. We will not discuss that further, since it is not necessary
for our immediate purpose.

To begin with, a brief review of Floer cohomology, just to remind the
reader of the special features of the ``exact'' situation. Let $M$ be
an exact symplectic manifold and $(L_0,L_1)$ a pair of transversally
intersecting exact Lagrangian submanifolds. The action functional on
the path space $\PP(L_0,L_1) = \{c \in \smooth([0;1], M) \suchthat
c(0) \in L_0, \; c(1) \in L_1 \}$ is
\[
a_{L_0,L_1}(c) = -\!\int c^*\theta \;+ K_{L_1}(c(1)) - K_{L_0}(c(0))
\]
(as pointed out to me by Oh, this form of the action functional
appears in \cite{ohnew}, and is then used crucially in his papers
with Milinkovic). Its critical points are the constant paths $c_x$ at
points $x \in L_0 \cap L_1$. We write $a_{L_0,L_1}(x) =
a_{L_0,L_1}(c_x) = K_{L_1}(x) - K_{L_0}(x)$ for their action. The
Floer cochain space $CF(L_0,L_1)$ is the vector space over $\Z/2$
with a basis given by these points; we denote the canonical basis
vectors by $\gen{x}$. Let $\JJ(M)$ be the space of smooth families $J
= (J_t)_{0 \leq t \leq 1}$ of almost complex structures on $M$, each
of which is $\o$-compatible and convex near the boundary. For $J \in
\JJ(M)$ and $x_{\pm} \in L_0 \cap L_1$, one considers the space
$\FF_J(x_-,x_+)$ of maps
\begin{equation} \label{eq:floer}
\begin{cases}
 & \!\! \sigma = \sigma(s,t): \R \times [0;1] \longrightarrow M, \\
 & \!\! \sigma(\R \times \{0\}) \subset L_0, \;
 \sigma(\R \times \{1\}) \subset L_1, \\
 & \!\! \partial \sigma/\partial t = J_t \,\partial \sigma/\partial s, \\
 & \!\! \lim_{s \rightarrow \pm \infty} \sigma(s,\cdot) = c_{x_{\pm}}
\end{cases}
\end{equation}
where the limit is understood to be in the $C^1$-topology on $\PP(L_0,L_1)$.
There is a natural $\R$-action on $\FF_J(x_-,x_+)$ by translation in
$s$-direction. We denote by $\FF_J^*(x_-,x_+)$ the subspace of maps which are
not $\R$-invariant. Solutions $\sigma$ of \eqref{eq:floer} can be interpreted
as negative gradient flow lines for $a_{L_0,L_1}$ in an $L^2$-metric on
$\PP(L_0,L_1)$; of course, this implies that $\FF_J^*(x_-,x_+) = \emptyset$
whenever $a_{L_0,L_1}(x_-) \leq a_{L_0,L_1}(x_+)$.

To each $\sigma \in \FF_J(x_-,x_+)$ one can associate an operator
$D_{\sigma,J}$ which linearizes \eqref{eq:floer}, and which is Fredholm in
suitable Sobolev spaces. There is a dense subspace $\JJ^{reg}(M,L_0,L_1)
\subset \JJ(M)$ of almost complex structures $J$ for which all $D_{\sigma,J}$
are onto, and then the spaces $\FF_J(x_-,x_+)$, as well as the quotients
$\FF_J^*(x_-,x_+)/\R$, are smooth manifolds. In that case one defines
$n_J(x_-,x_+) \in \Z/2$ to be the number mod 2 of isolated points in
$\FF_J^*(x_-,x_+)/\R$. The map
\[
d_J\gen{x_+} = \textstyle \sum_{x_-} n_J(x_-,x_+) \gen{x_-}
\]
has square zero, making $CF(L_0,L_1)$ into a differential vector space.
$HF(L_0,L_1)$ is defined to be its cohomology $\ker\, d_J/\im\, d_J$. Unlike
$d_J$ itself, Floer cohomology can be shown to be independent of $J$ up to
canonical isomorphism. It is also invariant under deformations of $L_0$ or
$L_1$ as exact Lagrangian submanifolds.

We now start constructing the relative invariant associated to $(E,\pi)$ and
$Q$. Choose a point $x_e \in L_{e,0} \cap L_{e,1}$ for each end $e$, and write
$\BB(\{x_e\})$ for the space of smooth sections $u: S \rightarrow E$ such that
$u(\partial S) \subset Q$, and which over the ends have the form $\Gamma_e^{-1}
\circ u \circ \gamma_e(s,t) = (s,t,\sigma_e(s,t))$ with maps $\sigma_e:
\R^{\pm} \times [0;1] \rightarrow M$ satisfying $\lim_{s \rightarrow \pm
\infty} \sigma_e(s,\cdot) = c_{x_e}$, in the same sense as in \eqref{eq:floer}.
The action integral $A(u) = \int_S u^*\o$ is convergent for all $u \in
\BB(\{x_e\})$, and one gets a formula analogous to \eqref{eq:action}:
\begin{equation} \label{eq:relative-action}
 A(u) = \sum_{e \in I^-} a_{L_{e,0},L_{e,1}}(x_e)
 - \sum_{e \in I^+} a_{L_{e,0},L_{e,1}}(x_e)
 + \int_{\partial S} \kappa_Q.
\end{equation}
Take a complex structure $j$ on $S$ (equal to $j_0$ near $S^\crit$ as always)
which, on the ends, is induced from the standard complex structure $j_{\pm}$ on
$\R^{\pm} \times [0;1]$. Choose also a $J_e = (J_{e,t})_{0 \leq t \leq 1} \in
\JJ(M)$ for each $e$. Then $\JJ(E,\pi,j,\{J_e\})$ denotes the contractible
space of almost complex structures $J$ on $E$ which are compatible relative to
$j$, and such that $\Gamma_e^*J$, for each $e$, is the almost complex structure
on $\R^{\pm} \times [0;1] \times M$ given by $j_{\pm} \times J_{e,t}$ at a
point $(s,t,y)$. For any such $J$, we denote by $\MM_J(\{x_e\}) \subset
\BB(\{x_e\})$ the subspace of sections $u$ which are $(j,J)$-holomorphic. The
formal picture from Section \ref{sec:simple-invariant} still applies: one has
the vector bundle $\EE_J \rightarrow \BB(\{x_e\})$, its canonical section
$\bar\partial_J$ whose zero-set is $\MM_J(\{x_e\})$, as well as their
``universal'' versions $\EE^{univ}$, $\bar\partial^{univ}$. Turning this
picture into an analytically realistic one is more complicated than in the
compact case. Still, the procedure is by now standard, and so we will not say
more about it, except to mention that the derivative $D_{u,J}$ of
$\bar\partial_J$ again extends to a Fredholm operator from a suitably defined
$W^{1,p}$-version $\WW^1_u$ of $T\BB(\{x_e\})_u$ to an $L^p$-version
$\WW^0_{u,J}$ of $\EE_{u,J}$. This allows one to define the notion of
regularity of $u \in \MM_J(\{x_e\})$, and of $J \in \JJ(E,\pi,j,\{J_e\})$, in
the same way as before. The subspace of regular $J$ is denoted by
$\JJreg(E,\pi,Q,j,\{J_e\})$.

\begin{remark}
It is maybe helpful to mention that Floer's equation \eqref{eq:floer} itself
can be made to fit into this framework. Namely, take the surface $S = \R \times
[0;1]$ and the trivial exact symplectic fibration $\pi: E = S \times M
\rightarrow S$, with the Lagrangian boundary condition $Q = (\R \times \{0\}
\times L_0) \cup (\R \times \{1\} \times L_1)$; it is a tautology that this is
modelled over the two ends of $S$ on $(L_0,L_1)$. Take the standard complex
structure $j$ on $S$, the same $J_e \in \JJ(M)$ for both ends $e$, and define
$J \in \JJ(E,\pi,j,\{J_e\})$ to be the product $j \times J_{e,t}$ at a point
$(s,t,y)$. Sections $u: S \rightarrow E$ with $u(\partial S) \subset Q$ are of
the form $u(s,t) = (s,t,\sigma(s,t))$, where $\sigma: \R \times [0;1]
\rightarrow M$ is a map satisfying the boundary condition in \eqref{eq:floer}.
Moreover, $u$ is $(j,J)$-holomorphic iff $\partial\sigma/\partial t = J_{e,t}
\partial\sigma/\partial s$, so that one can identify $\MM_J(x_-,x_+) =
\FF_{J_e}(x_-,x_+)$ for all $x_{\pm}$.
\end{remark}

In parallel with the exposition in Section \ref{sec:simple-invariant}, we will
now discuss the basic properties of the spaces $\MM_J(\{x_e\})$. The next two
results are analogues of Lemmas \ref{th:convexity} and \ref{th:transversality},
and their proofs are the same.

\begin{lemma} \label{th:relative-convexity}
For every $J \in \JJ(E,\pi,j,\{J_e\})$ there is a closed subset $K \subset E
\setminus \partial_hE$ such that $\pi|K : K \rightarrow S$ is proper, and which
over the ends has the form $\Gamma_e^{-1}(K) = \R^{\pm} \times [0;1] \times
K_e$ for some compact $K_e \subset M \setminus \partial M$, such that $u(S)
\subset K$ for all $u \in \MM_J(\{x_e\})$ and all points $\{x_e\}$. \qed
\end{lemma}

\begin{lemma} \label{th:relative-transversality}
$\JJreg(E,\pi,Q,j,\{J_e\})$ is $\smooth$-dense in $\JJ(E,\pi,j,\{J_e\})$. More
precisely, given some nonempty open subset $U \subset S$ which is disjoint from
the ends, and a $J \in \JJ(E,\pi,j,\{J_e\})$, there are $J' \in
\JJreg(E,\pi,Q,j,\{J_e\})$ arbitrarily close to $J$, such that $J = J'$ outside
$\pi^{-1}(U)$. \qed
\end{lemma}

The remaining issue is compactness. After taking a symplectic form $\Omega +
\pi^*\beta$ on $E$ as in Lemma \ref{th:tame}, and the associated metric, one
finds that $\half \int_W || Du ||^2 \leq A(u) + \int_W \beta$ for any compact
subset $W \subset S$ and $u \in \MM_J(\{x_e\})$. Repeating the argument in
Lemma \ref{th:compactness}, one derives from this that any sequence $(u_i)$ in
$\MM_J(\{x_e\})$ has a subsequence which is $C^r$-convergent on compact
subsets. It is necessary to go beyond this somewhat coarse result, and to study
sequences $(u_i)$ in the more appropriate {\em Gromov-Floer topology}. A
compactification $\overline{\MM}_J(\{x_e\})$ of $\MM_J(\{x_e\})$ in this
topology can be constructed by adding ``broken sections''. Since a very similar
notion is part of the standard analytical package underlying Floer theory, we
will not define the topology, and only describe the compactification as a set.
Each point of it consists of
\begin{condensedlist}
\item
the ``principal component'' $u \in \MM_J(\{\hat{x}_e\})$ for some
$\{\hat{x}_e\}$;
\item
for each end $e \in I^{\pm}$, a finite sequence of points $\hat{x}_{e,0},
\dots, \hat{x}_{e,l_e} \in L_{e,0} \cap L_{e,1}$, $l_e \geq 0$. If the end is
negative (positive), this should satisfy $\hat{x}_{e,0} = x_e$ and
$\hat{x}_{e,l_e} = \hat{x}_e$ ($\hat{x}_{e,0} = \hat{x}_e$ and $\hat{x}_{e,l_e}
= x_e$);
\item
Floer flow lines $\sigma_{e,m} \in
\FF^*_{J_e}(\hat{x}_{e,m-1},\hat{x}_{e,m})/\R$, for each $e$ and $1 \leq m \leq
l_e$.
\end{condensedlist}
Suppose that a sequence $(u_i)$ converges to such a limit. One sees from
\eqref{eq:relative-action} that for all $i$,
\begin{equation} \label{eq:energy-additivity}
 A(u_i) = A(u) + \sum_e \sum_{1 \leq m \leq l_e}
 (a_{L_{e,0},L_{e,1}}(\hat{x}_{e,m-1}) - a_{L_{e,0},L_{e,1}}(\hat{x}_{e,m}))
\end{equation}
Note that all terms of the $\sum$ are $>0$; therefore $A(u) \leq A(u_i)$, with
equality iff $l_e = 0$ for all ends $e$. We will need another piece of
information about the limit, which can be derived from the definition of the
Gromov-Floer topology and a ``gluing theorem'' for linear elliptic operators.
Namely, for $i \gg 0$,
\begin{equation} \label{eq:index-additivity}
 \ind\,D_{u_i,J} = \ind \, D_{u,J} + \sum_e \sum_{1 \leq m \leq l_e}
 \ind \, D_{\sigma_{e,m},J_e}.
\end{equation}

From now on assume that $J_e \in \JJreg(M,L_{e,0},L_{e,1})$ for all $e$, and
that $J \in \JJreg(E,\pi,Q,j,\{J_e\})$. Then $\ind\,D_{u,J} \geq 0$ and
$\ind\,D_{\sigma_{e,m},J_e} > 0$ on the right hand side of
\eqref{eq:index-additivity}, so that the left hand side can be zero only if
$\ind\,D_{u,J} = 0$ and $l_e = 0$ for all $e$. It follows that the
zero-dimensional part of $\MM_J(\{x_e\})$, for any $\{x_e\}$, is compact, hence
a finite set. Write $\nu_J(\{x_e\}) \in \Z/2$ for the number of points modulo
two in this set, and consider the map
\begin{align} \label{eq:relative-chain-map}
 & C\Phi_0^{rel}(E,\pi,Q,J) : \bigotimes_{e \in I^+} CF(L_{e,0},L_{e,1})
 \longrightarrow \bigotimes_{e \in I^-} CF(L_{e,0},L_{e,1}) \\
 \intertext{given by the matrix with entries $\nu_J(\{x_e\})$, that is to
 say}
\notag
 & C\Phi_0^{rel}(E,\pi,Q,J)(\underset{e \in I^+\!\!}{\otimes} \gen{x_e}) =
 \sum_{\;\;\;\;\{x_e\}_{e \in I^-}\!\!\!\!} \nu_J(\{x_e\}_{e \in I^- \cup I^+})
 (\underset{e \in I^-\!\!}{\otimes} \gen{x_e}).
\end{align}
A standard argument, involving the structure at infinity of the one-dimensional
part of $\MM_J(\{x_e\})$, shows is that there is an even number of points in
$\overline{\MM}_J(\{x_e\})$ of the following form: $l_e$ is $1$ for a single
end $e=f$, and zero for all other $e$, so that the point is a pair
$(u,\sigma_{f,1})$; and moreover $\ind\,D_{u,J} = 0$, $\ind \,
D_{\sigma_{f,1},J_f} = 1$. Algebraically, what this says is that
$C\Phi_0^{rel}(E,\pi,Q,J)$ is a chain map. The relative invariant
$\Phi_0^{rel}(E,\pi,Q)$ is defined to be the induced map on cohomology.

The next step is to show that this is independent of the choice of $j$ and $J$,
keeping the $J_e$ fixed for the moment. Any two $j^0,j^1$ can be connected by a
family $j^\mu$, $0 \leq \mu \leq 1$, which remains constant on the ends;
correspondingly, for $J^0 \in \JJreg(E,\pi,Q,j^0,\{J^e\})$, $J^1 \in
\JJreg(E,\pi,Q,j^1,\{J^e\})$ one can find a family $J^\mu$ joining them, which
is regular in the parametrized sense (this is not the same as saying that each
$J^\mu$ should itself be regular, which would be impossible to achieve in
general). The parametrized moduli spaces
\[
\MM^{para}_{(J^\mu)}(\{x_e\}) = \bigcup_{\!\!\! 0 \leq \mu \leq 1 \!\!\!} \,
\{\mu\} \times \MM_{J^\mu}(\{x_e\})
\]
are then smooth manifolds, and they have a parametrized version of the
Gromov-Floer compactification. In particular, the zero-dimensional parts are
again finite sets, so that one can use the number of points $\lambda(\{x_e\})
\in \Z/2$ in them to define a map $h(E,\pi,Q,(J^\mu))$ between the same groups
as in \eqref{eq:relative-chain-map}. Arguing along the same lines as when
proving that $C\Phi_0^{rel}$ is a chain map, one can show that $h$ is a
homotopy between $C\Phi_0^{rel}(E,\pi,Q,J^0)$ and $C\Phi_0^{rel}(E,\pi,Q,J^1)$.
The same argument can be used to show that $\Phi_0^{rel}(E,\pi,Q)$ remains
invariant under deformations of the geometric data, that is to say of $Q$ or
$(E,\pi)$ itself, as long as the structure of the ends remains unchanged.
Finally, we should prove that the relative invariant is independent of the
$J_e$, which more accurately means that it commutes with the canonical
isomorphisms between Floer cohomology groups for different $J_e$. We omit this
entirely, both because it is not important for our purpose, and because it
would require a digression concerning ``continuation maps''.

Two ways of gluing together surfaces with strip-like ends will play a role
later on. One of them is a close cousin of that considered as in Section
\ref{sec:simple-invariant}. It can be formulated in various degrees of
generality, but we will need only one special case. Suppose then that $S^1$ is
a surface with strip-like ends, and $S^2$ a compact surface, together with
points $\zeta^k \in \partial S^k$ ($\zeta^1$ should not lie on any end; that
can of course always be achieved by making the ends smaller). Let $(E^k,\pi^k)$
be exact Lefschetz fibrations over $S^k$ with Lagrangian boundary conditions
$Q^k$; $(E^1,\pi^1)$ should be trivial over the ends. We also want to have $M$,
$L$, maps $\phi^k: M \rightarrow (E^k)_{\zeta^k}$, and trivializations
$(\psi^k,\Psi^k)$, with the same properties as in Section
\ref{sec:simple-invariant}. The boundary connected sum $S = S^1 \#_{\zeta^1
\sim \zeta^2} S^2$ is a surface with the same kind of strip-like ends as $S^1$.
As before one constructs an exact Lefschetz fibration $(E,\pi)$ on it with a
Lagrangian boundary condition $Q$, modelled over the ends on the same
Lagrangian submanifolds as $Q^1$.

Choose complex structures $j^k$ on $S^k$ such that $(\psi^k)^*j^k$ is standard;
$j^1$ should also be standard on the ends. Take $\{J_e\}$ and $J^1 \in
\JJreg(E^1,\pi^1,Q^1,j^1,\{J_e\})$ as when defining the relative invariant of
$(E^1,\pi^1,Q^1)$; as an additional condition, we want $(\Psi^1)^*J^1$ to be
the product of the standard complex structure on $\bar{D}^+(1)$ and some fixed
almost complex structure on $M$. By Lemma \ref{th:relative-transversality} this
can be done while still achieving regularity. Using Lemmas
\ref{th:transversality} and \ref{th:ev-transversality} one finds a $J^2 \in
\JJreg(E^2,\pi^2,Q^2,j^2)$ with the same restriction on $(\Psi^2)^*J^2$, and
such that the evaluations
\begin{equation} \label{eq:two-evaluations}
\begin{split}
 & ev_{\zeta^1}\,|\,\MM_{J^1}(\{x_e\}) :
 \MM_{J^1}(\{x_e\}) \longrightarrow Q^1_{\zeta^1} \iso L, \\
 & ev_{\zeta^2}\,|\,\MM_{J^2}: \MM_{J^2} \longrightarrow Q^2_{\zeta^2} \iso L
\end{split}
\end{equation}
are transverse to each other for every $\{x_e\}$. Let $J \in
\JJ(E,\pi,j,\{J_e\})$ be the almost complex structure constructed from $J^1,
J^2$. As in Proposition \ref{th:gluing}, this will be regular for small values
of the parameter $\rho$, and
\begin{equation} \label{eq:zero-gluing}
\MM_J(\{x_e\})_{[0]} \iso (\MM_{J^1}(\{x_e\}) \times_L \MM_{J^2})_{[0]},
\end{equation}
where the $[0]$ denotes on both sides the zero-dimensional component of these
manifolds. We will not need the full ``gluing formula'' which one can obtain
from this, but only a special case:

\begin{lemma} \label{th:gluing-vanishing}
If $\Phi_1(E^2,\pi^2,Q^2) = 0$ then $\Phi_0^{rel}(E,\pi,Q) = 0$.
\end{lemma}

\proof For simplicity, suppose that $\Phi_1(E^2,\pi^2,Q^2)$ vanishes even when
taken in the cobordism ring $MO_*(L)$. This means that there is a compact
manifold with boundary $G$ and a smooth map $g: G \rightarrow L$, such that
$\partial G = \MM_{J^2}$ and $g|\partial G = ev_{\zeta^2}$. After perturbing it
slightly, one can assume that $g$ is transverse to all maps $ev_{\zeta^1}$ in
\eqref{eq:two-evaluations}. Then the fibre products
\begin{equation} \label{eq:fibre-product}
\GG(\{x_e\}) = \MM_{J^1}(\{x_e\}) \times_L G
\end{equation}
are smooth manifolds. The evaluation map extends continuously to the
Gromov-Floer compactification, so that one can define compactifications
$\overline{\GG}(\{x_e\})$ in the obvious way. It is not difficult to see, using
\eqref{eq:index-additivity}, that the zero-dimensional part
$\GG(\{x_e\})_{[0]}$ is a finite set; one counts the points
$\xi_{J,G,g}(\{x_e\}) \in \Z/2$ in it, and uses that to define a map $k$
between the Floer cochain groups associated to $(E^1,\pi^1,Q^1)$, as in
\eqref{eq:relative-chain-map}. We claim that this is a homotopy between
$C\Phi^0_{rel}(E,\pi,Q,J)$ and the zero map. As usual, the proof is based on
analyzing the ends of the one-dimensional moduli spaces. The closure in
$\overline{\GG}(\{x_e\})$ of the one-dimensional part $\GG(\{x_e\})_{[1]}$ is a
compact one-manifold with boundary, and its boundary points are of two kinds:
first, boundary points of $\GG(\{x_e\})_{[1]}$ itself,
$\partial\GG(\{x_e\})_{[1]} = (\MM_{J^1}(\{x_e\}) \times_L \partial G)_{[0]}
\iso \MM_J(\{x_e\})_{[0]}$; their number modulo two is $\nu_J(\{x_e\})$ by
definition. The second kind of boundary points are of the form
$(u,\sigma_{f,1}) \times q \in \overline{\MM}_J(\{x_e\}) \times_L G$. This
means that for some end $f$, and points $\{\hat{x}_e\}$ with $\hat{x}_e = x_e$
for all $e \neq f$, one has
\[
 u \times q \in \GG(\{\hat{x}_e\})_{[0]}, \quad
 \sigma_{f,1} \in \begin{cases}
 (\FF^*(x_f,\hat{x}_f)/\R)_{[0]} & \text{if $e$ is a negative end,} \\
 (\FF^*(\hat{x}_f,x_f)/\R)_{[0]} & \text{if $e$ is a positive end.}
 \end{cases}
\]
The number of such boundary points is
\[
\sum_{f \in I^-} \sum_{\hat{x}_f}
  n_{J_f}(x_f,\hat{x}_f) \,
 \xi_{J,G,g}(\{\hat{x}_e\})
 +
\sum_{f \in I^+} \sum_{\hat{x}_f}
 \xi_{J,G,g}(\{\hat{x}_e\}) \,
 n_{J_f}(\hat{x}_f,x_f)
\in \Z/2.
\]
The fact that this is equal to $\nu_J(\{x_e\})$ gives precisely the desired
equality $dk + kd = C\Phi_0^{rel}(E,\pi,Q,J)$. \qed

The other and maybe more obvious gluing process is to join together two ends.
Assume that $S^1,S^2$ are surfaces with strip-like ends, such that $S^1$ has a
single positive end $e^1$, and $S^2$ a single negative end $e^2$ (this is not
really a restriction since, as has been mentioned before, positive ends can be
turned into negative ones and vice versa). Choose some $\sigma>0$, and define
$S = S^1 \#_{e^1 \sim e^2} S^2$ by taking $S^1 \setminus
\gamma_{e^1}((\sigma;\infty) \times [0;1])$ and $S^2 \setminus
\gamma_{e^2}((-\infty;-\sigma) \times [0;1])$, and identifying
$\gamma_{e^1}(s,t)$ with $\gamma_{e^2}(s-\sigma,t)$ for $(s,t) \in [0;\sigma]
\times [0;1]$. The ends of $S$ are the negative ends of $S^1$ together with the
positive ends of $S^2$. After choosing complex structures $j^k$ on $S^k$ which
are standard over the ends, there is an obvious induced complex structure $j$
on $S$.

Let $(L_0,L_1)$ be a pair of exact Lagrangian submanifolds in $M$. Suppose that
we have exact Lefschetz fibrations $(E^k,\pi^k)$ over $S^k$, trivial over the
ends, and Lagrangian boundary conditions $Q^1,Q^2$ for them modelled on
$(L_0,L_1)$ over $e^1$ and over $e^2$, respectively. Then one can form an exact
Lefschetz fibration $(E,\pi)$ over $S$ by identifying $\Gamma_{e^1}(s,t,y) \in
E^1$ with $\Gamma_{e^2}(s-\sigma,t,y) \in E^2$, in parallel with the
construction on the base. This comes with an obvious Lagrangian boundary
condition $Q$. Choose almost complex structures $J^k$ on $E^k$ so as to define
relative invariants $C\Phi_0^{rel}(E^k,\pi^k,Q^k,J^k)$. We require that the
$J_{e^k} \in \JJreg(M,L_0,L_1)$ on which $J^k$ is modelled over the end $e^k$
should be the same for $k = 1,2$. Then $J^1$ and $J^2$ match up to an almost
complex structure $J$ on $E$, which is compatible relative to $j$.

\begin{prop} \label{th:gluing-ends}
For fixed $J^1,J^2$, if one chooses $\sigma$ to be sufficiently large, then $J$
is regular. Moreover, denoting again by $[0]$ the zero-dimensional components,
and by $I^-,I^+$ the negative and positive ends of $S$, one has
\[
\MM_J(\{x_e\})_{[0]} \iso \bigcup_{\!\!\! x \in L_0 \cap L_1 \!\!\!}
\MM_{J^1}(\{x_e\}_{e \in I^-},x)_{[0]} \times \MM_{J^2}(x,\{x_e\}_{e \in
I^+})_{[0]}
\]
for any $\{x_e\}_{e \in I^- \cup I^+}$. \qed
\end{prop}

The method used in the proof of this is the same as when setting up Floer
cohomology. A thorough exposition of a closely related result can be found in
\cite[Section 4.4]{schwarz95}. The implication for the relative invariants is
clear:
\begin{equation} \label{eq:chain-composition}
 C\Phi_0^{rel}(E,\pi,Q,J) = C\Phi_0^{rel}(E^1,\pi^1,Q^1,J^1) \circ
 C\Phi_0^{rel}(E^2,\pi^2,Q^2,J^2),
\end{equation}
which proves the main TQFT-style property of relative invariants, namely, that
they are functorial if one regards the gluing $S = S^1 \#_{e^1 \sim e^2} S^2$
as composition of the surfaces $S^1,S^2$.

\subsection{Horizontality and relative invariants\label{sec:horizontal-two}}

The aim of this section is to extend the methods of Section
\ref{sec:horizontal} to surfaces with strip-like ends. Throughout, $(E,\pi)$
will be an exact Lefschetz fibration over such a surface $S$, trivial on the
ends, together with a Lagrangian boundary condition $Q$ modelled over the ends
on pairs $L_{0,e},L_{1,e} \subset M$ of transversally intersecting exact
Lagrangian submanifolds. For $\{x_e\}_{e \in I^- \cup I^+}$ with $x_e \in
L_{0,e} \cap L_{1,e}$ write
\[
\chi(\{x_e\}) = \sum_{e \in I^-} a_{L_{0,e},L_{1,e}}(x_e) - \sum_{e \in I^+}
a_{L_{0,e},L_{1,e}}(x_e).
\]
By \eqref{eq:relative-action} $A(u) = \chi(\{x_e\}) + \int_{\partial S}
\kappa_Q$ for all $u \in \BB(\{x_e\})$. After fixing a complex structure $j$ on
$S$ which is standard on the ends, and a $J_e \in \JJ(M)$ for each $e$, we
consider the space
\[
\JJ^h(E,\pi,j,\{J_e\}) = \JJ^h(E,\pi,j) \cap \JJ(E,\pi,j,\{J_e\})
\]
of almost complex structures $J$, compatible relative to $j$, which are
horizontal and have prescribed behaviour on the ends. The two conditions do not
contradict each other, since the model on each end, the almost complex
structure $j_{\pm} \times J_e$ on the trivial fibration $\R^{\pm} \times [0;1]
\times M \rightarrow \R^{\pm} \times [0;1]$, is obviously horizontal. A little
more thought shows that $\JJ^h(E,\pi,j,\{J_e\})$ is contractible. Choosing $J$
in that space, one finds that for all $u \in \BB(\{x_e\})$,
\begin{equation} \label{eq:rel-deficiency}
{\textstyle \half} \int_S ||(Du)^v||^2 + \int_S f(u)\beta = A(u) + \int_S
||\bar\partial_J u||^2.
\end{equation}
Here $||\cdot||$ is the metric on $TE^v$ associated to $\Omega$ and $J$; $\beta
\in \Omega^2(S)$ is a positive two-form; and $f$ is the function defined by
$\Omega|TE^h = f(\pi^*\beta|TE^h)$. Because the symplectic connection is
trivial on the ends, $f(u)$ is compactly supported. \eqref{eq:rel-deficiency}
is obviously the analogue of \eqref{eq:deficiency}, and can be proved in the
same way. Recall that $(E,\pi)$ has nonnegative curvature iff $f \geq 0$. From
this we draw two conclusions:

\begin{lemma} \label{th:empty}
Assume that $(E,\pi)$ has nonnegative curvature, and choose $J \in
\JJ^h(E,\pi,j,\{J_e\})$. Then $\MM_J(\{x_e\}) = \emptyset$ for all $\{x_e\}$
with $\chi(\{x_e\}) < -\int_{\partial S} \kappa_Q$. \qed
\end{lemma} \vspace{\parskip}

\begin{lemma} \label{th:energy-compactness}
In the same situation, let $\alpha$ be the minimum of $a_{L_{0,e},L_{1,e}}(x_-)
- a_{L_{0,e},L_{1,e}}(x_+)$, taken over all $e$ and all $x_{\pm} \in L_{0,e}
\cap L_{1,e}$ for which this is $>0$. Then the spaces $\MM_J(\{x_e\})$ is
compact, in the Gromov-Floer topology, for all $\{x_e\}$ such that
$\chi(\{x_e\}) < -\int_{\partial S} \kappa_Q + \alpha$.
\end{lemma}

\proof Consider a sequence $(u_i)$ in $\MM_J(\{x_e\})$ which converges to
$(u,\{\sigma_{e,m}\}) \in \overline{\MM}_J(\{x_e\})$. The definition of
$\gamma$ implies that in \eqref{eq:energy-additivity}, each
$a_{L_{e,0},L_{e,1}}(\hat{x}_{e,m-1}) - a_{L_{e,0},L_{e,1}}(\hat{x}_{e,m}) \geq
\alpha$; therefore
\[
A(u) \leq A(u_i) - \big(\textstyle\sum_e l_e\big) \alpha.
\]
By assumption $A(u_i) = \chi(\{x_e\}) + \int_{\partial S} \kappa_Q < \alpha$,
while on the other hand, $A(u) \geq 0$ by nonnegative curvature and
\eqref{eq:rel-deficiency}. This is only possible if $l_e = 0$ for all $e$, so
that the limit actually lies in $\MM_J(\{x_e\})$. \qed

As before, denote by $\MM^h$ the space of horizontal sections with boundary
values in $Q$. The presence of strip-like ends makes this space somewhat
simpler than in the case which we have encountered before. Take an arbitrary $u
\in \BB$ and consider $u_e(s,t) = \Gamma_e^{-1} \circ u \circ \gamma_e(s,t) =
(s,t,\sigma_e(s,t))$; $u$ is horizontal over the end $e$ iff $\sigma_e \equiv x
\in M$ is constant, in which case the boundary conditions $\sigma_e(s,0) \in
L_{e,0}$, $\sigma_e(s,1) \in L_{e,1}$ imply that $x \in L_{e,0} \cap L_{e,1}$.
Since a horizontal section is determined by its value at any one point, for any
$e$ and any $x \in L_{e,0} \cap L_{e,1}$ there can be at most one $u \in \MM^h$
such that $\sigma_e \equiv x$. Thus $\MM^h$ is a finite disjoint union of the
subsets $\MM^h(\{x_e\}) = \MM^h \cap \BB(\{x_e\})$, each of which consists of
at most one element. One easily proves the following limiting case of Lemma
\ref{th:empty}:

\begin{lemma} \label{th:almost-empty}
For $(E,\pi)$ and $J$ as in Lemma \ref{th:empty}, suppose that $\{x_e\}$
satisfy $\chi(\{x_e\}) = -\int_{\partial S} \kappa_Q$; then $\MM_J(\{x_e\}) =
\MM^h(\{x_e\})$. \qed
\end{lemma}

To translate these elementary observations into results about relative
invariants, one needs to address again the question of regularity of horizontal
almost complex structures; in other words, what is required are analogues of
Lemmas \ref{th:horizontal-transversality} and \ref{th:flat-deformations}. For
the first of these, both the statement and proof as essentially the same as in
the original situation; the second needs to be adapted a little.

\begin{lemma} \label{th:relative-h-transversality}
Let $U \subset S$ be a nonempty open subset disjoint from the ends, such that
any partial section $w: U \rightarrow E|U$ which is horizontal and satisfies
$w(\partial S \cap U) \subset Q$ is the restriction of a $u \in \MM^h$. Then,
given some $J \in \JJ^h(E,\pi,j,\{J_e\})$, there are $J' \in
\JJ^h(E,\pi,j,\{J_e\})$ arbitrarily close to it and which agree with it outside
$\pi^{-1}(U)$, such that for all $\{x_e\}$, any $u \in \MM_J(\{x_e\}) \setminus
\MM^h$ is regular. \qed
\end{lemma} \vspace{\parskip}

\begin{lemma} \label{th:negative-index}
Let $(E,\pi)$ and $J$ be as in Lemma \ref{th:empty}. If $u \in \MM^h$ is a
horizontal section with $A(u) = 0$, then $ker \, D_{u,J} = 0$.
\end{lemma}

\proof Formally, taking the second derivative of \eqref{eq:rel-deficiency} at
$u$ yields the same formula \eqref{eq:weitzenboeck} as in the case of a compact
$S$, but a little care needs to be exercised about its validity. It certainly
holds for those elements of $T\BB_u = \{ X \in \smooth(u^*TE^v) \suchthat X_z
\in T(Q_z) \text{ for } z \in \partial S\}$ which are compactly supported, and
by continuity, for all $X$ in the $W^{1,2}$-completion; to be precise, this
Sobolev space is with respect to the metric $||\cdot||$ and the connection
$\nabla^u$ on $u^*TE^v$. We actually want to use the formula with $X \in
\WW^1_u$ and this is a space of $W^{1,p}$-sections with $p>2$, hence not
contained in $W^{1,2}$. However, if one assumes additionally that $D_{u,J}X =
0$ there is no problem, because any such $X$ is smooth and decays exponentially
on the ends, as do its derivatives.

From $A(u) = 0$ and \eqref{eq:rel-deficiency} it follows that $f(u) \equiv 0$,
so that the Hessians $Hess(f|E_z)_{u(z)}$ are nonnegative. One then sees from
\eqref{eq:weitzenboeck} that any $X \in \WW^1_u$ with $D_{u,J}X = 0$ satisfies
$\nabla^uX = 0$. Choose some end $e$ and write $x = x_e$. The trivialization
$\Gamma_e$ induces a trivialization of the vector bundle $\gamma_e^*(u^*TE^v)
\rightarrow \R^{\pm} \times [0;1]$, which identifies it with the trivial bundle
with fibre $TM_x$. In this way, $Y = \gamma_e^*X$ becomes a map $\R^{\pm}
\times [0;1] \rightarrow TM_x$. Since the connection $\nabla^u$ is compatible
with the trivialization, $\nabla^uX = 0$ implies that $Y$ must be constant. On
the other hand, the boundary conditions which are part of the definition of
$T\BB_u$ and of $\WW^1_u$ tell us that $Y_{s,k} \in T(L_{e,k})_x$ for $k =
0,1$. Because the $L_{e,k}$ intersect transversally, it follows that $Y = 0$,
hence that $X | \im(\gamma^e) = 0$. Since $X$ is covariantly constant, it must
be zero everywhere. \qed

The next result summarizes what progress we have made so far, as well as the
implications for the coefficients $\nu_J(\{x_e\})$ of the relative invariant.

\begin{prop} \parindent0em \label{th:nonnegative-one}
Assume that $(E,\pi)$ has nonnegative curvature, and that any $u \in \MM^h$
satisfies $A(u) = 0$ and $\ind\,D_{u,J} = 0$. Set $\kappa = \int_{\partial S}
\kappa_Q$.

(i) Let $U \subset S$ be a nonempty open subset, disjoint from the ends, such
that any partial horizontal section $w: U \rightarrow E|U$ with $w(\partial S
\cap U) \subset Q$ extends to a $u \in \MM^h$. Then for any $J \in
\JJ^h(E,\pi,j,\{J_e\})$ there are $J' \in \JJ^{reg,h}(E,\pi,Q,j,\{J_e\}) =
\JJ^h(E,\pi,j,\{J_e\}) \cap \JJ^{reg}(E,\pi,Q,j,\{J_e\})$ arbitrarily close to
$J$, and which agree with it outside $\pi^{-1}(U)$.

(ii) If $J_e \in \JJreg(M,L_{0,e},L_{1,e})$ for all $e$, and $J \in
\JJ^{reg,h}(E,\pi,Q,j,\{J_e\})$, then
\[
 \nu_J(\{x_e\}) =
 \begin{cases}
 0 & \text{if $\chi(\{x_e\}) < -\kappa$,} \\
 \#\MM^h(\{x_e\}) & \text{if $\chi(\{x_e\}) = -\kappa$.}
 \end{cases}
\]
\end{prop}

\proof (i) Take a $J'$ as given by Lemma \ref{th:relative-h-transversality}.
Then all $u \in \MM_{J'}(\{x_e\})$ are regular except possibly for the
horizontal ones. These satisfy $A(u) = 0$, so we can apply Lemma
\ref{th:negative-index} to them, showing that $\ker\,D_{u,J'} = 0$, and by
assumption on the index, that $\coker\,D_{u,J'} = 0$; which means that they are
regular as well. (ii) This follows immediately from Lemmas \ref{th:empty} and
\ref{th:almost-empty}, in view of the definition of $\nu_J(\{x_e\})$. \qed

An algebraic language suitable for encoding results of this kind is that of
$\R$-graded vector spaces, that is to say vector spaces $C$ equipped with a
splitting $C = \bigoplus_{r \in \R} C_r$. All vector spaces occurring here will
be over $\Z/2$ and finite-dimensional; in particular, their support $supp(C) =
\{r \in \R \suchthat C_r \neq 0\}$ is always a finite set. Let $I \subset \R$
be an interval. We say that $C$ has {\em gap $I$} if there are no $r,s \in
supp(C)$ with $r-s \in I$. A map $f: C \rightarrow D$ between graded
$\R$-vector spaces is said to be {\em of order $I$} if $f(C_r) \subset
\bigoplus_{s \in r + I} D_s$ for all $r$.

Floer cochain groups are obvious examples: $C = CF(L_0,L_1)$ is canonically
$\R$-graded, with $C_r$ the subspace spanned by those $\gen{x}$ with
$a_{L_0,L_1}(x) = r$. Because of its gradient flow interpretation, the Floer
differential $d_J$ is of order $(0;\infty)$. In a rather trivial way, this can
always be strengthened slightly; namely, there is an $\alpha>0$ such that $C$
has gap $(0;\alpha)$, and then $d_J$ is of order $[\alpha;\infty)$. One can
reformulate Proposition \ref{th:nonnegative-one}(ii) in this language as
follows:

\begin{lemma} \label{th:nonnegative-two}
Take $(E,\pi)$, $Q$ and $\kappa$ as in Proposition \ref{th:nonnegative-one},
with $J_e,J$ as in part (ii). Then the map $C\Phi^0_{rel}(E,\pi,Q,J)$ is of
order $[-\kappa;\infty)$. For a more precise statement, take $\alpha>0$ to be
the minimum of $\chi(\{x_e\})+\kappa$, ranging over all $\{x_e\}$ where this is
positive. Then $C\Phi_0^{rel}(E,\pi,Q,J) = \phi + (C\Phi_0^{rel}(E,\pi,Q,J) -
\phi)$, where the first summand is of order $\{-\kappa\}$ and the second of
order $[-\kappa+\alpha;\infty)$. Moreover, $\phi$ is determined by the
horizontal sections:
\[
\phi(\underset{e \in I^+\!\!}{\otimes} \gen{x_e}) =
 \sum_{\;\;\;\;\{x_e\}_{e \in I^-}\!\!\!\!} \#\MM^h(\{x_e\}_{e \in I^- \cup I^+})
 (\underset{e \in I^-\!\!}{\otimes} \gen{x_e}). \qed
\]
\end{lemma}

There is a version of this for the homotopies between the $C\Phi^0_{rel}$ for
different choices of $(j,J)$. The proof applies the same ideas as before to
parametrized moduli spaces, and is left to the reader.

\begin{lemma} \label{th:nonnegative-homotopy}
Let $(E,\pi)$, $Q$ and $\kappa$ be as in Proposition \ref{th:nonnegative-one}.
Consider two complex structures $j^k$ on $S$, $k = 0,1$, and correspondingly
two almost complex structures $J^k \in \JJ^{reg,h}(E,\pi,Q,j,\{J_e\})$; note
that the $J_e$ are supposed to be the same for both $k$. Then the maps
$C\Phi^0_{rel}(E,\pi,Q,J^k)$ are homotopic by a chain homotopy which is of
order $(-\kappa;\infty)$. \qed
\end{lemma}

It is a familiar idea that when dealing with maps between $\R$-graded vector
spaces, the ``lowest order term'' is usually the most important, and knowing it
is often sufficient to resolve a question. A particular instance of this is
relevant for our purpose.

\begin{lemma} \label{th:spectral-sequence}
Let $D$ be an $\R$-graded vector space with a differential $d_D$ of order
$[0;\infty)$. Suppose that $D$ has gap $[\epsilon;2\epsilon)$ for some
$\epsilon>0$. One can then write $d_D = \delta + (d_D-\delta)$ with $\delta$ of
order $[0;\epsilon)$, satisfying $\delta^2 = 0$, and $(d_D-\delta)$ of order
$[2\epsilon;\infty)$. Suppose that in addition, $H(D,\delta) = 0$; then
$H(D,d_D) = 0$.
\end{lemma}

\proof Thanks to the gap assumption, $supp(D)$ can be decomposed into disjoint
subsets $R_1,\dots,R_m$ such that for $r \in R_i$, $s \in R_j$,
\[
r-s \; \begin{cases}
 \leq -2\epsilon & i < j, \\
 \in (-\epsilon;\epsilon) & i = j, \\
 \geq 2\epsilon & i > j.
\end{cases}
\]
Define a descending filtration of $(D,d_D)$,
\[
F^k = \textstyle\bigoplus_{r \in R_k \cup R_{k+1} \cup \dots \cup R_m} D_r.
\]
There is a ``spectral sequence'' which takes the form of a sequence
$(E^k,\partial^k)$ of differential vector spaces, such that $E^{k+1} =
H(E^k,\partial^k)$. It starts with $E^0 = \bigoplus_i F^i/F^{i+1}$, which has a
differential $\partial^0$ induced by $d_D$. In our case this can be identified
with $(D,\delta)$, so the assumption says that $E^1 = 0$. On the other hand
$E^k \iso H(D,d_D)$ for $k \gg 0$. \qed

\begin{lemma} \label{th:low-energy}
Take three $\R$-graded vector spaces $C',C,C''$, each of them with a
differential of order $(0;\infty)$. Suppose that we have differential maps $b:
C' \rightarrow C$, $c: C \rightarrow C''$ and a homotopy $h: C' \rightarrow
C''$ between $c \circ b$ and the zero map, such that that the following
conditions are satisfied for some $\epsilon>0$:
\begin{condensedlist}
\item \label{item:gap-one}
$C',C''$ have gap $(0;3\epsilon)$, and $C$ has gap $(0;2\epsilon)$.
\item \label{item:gap-two}
For all $r \in supp(C')$ and $s \in supp(C'')$, $|r-s| \geq 4\epsilon$.
\item \label{item:gap-three}
One can write $b = \beta + (b - \beta)$ with $\beta$ of order $[0;\epsilon)$
and $(b - \beta)$ of order $[2\epsilon;\infty)$; and $c = \gamma + (c-\gamma)$
with the same properties. The low order parts (which do not need be
differential maps) fit into a short exact sequence of vector spaces
\[
0 \rightarrow C' \xrightarrow{\beta} C \xrightarrow{\gamma} C'' \rightarrow 0.
\]
\item \label{item:gap-four}
$h$ is of order $[0;\infty)$.
\end{condensedlist}
Then the maps on cohomology induced by $b,c$ fit into a long exact sequence
\begin{equation} \label{eq:induced-sequence}
\xymatrix{
 {H(C',d_{C'})} \ar[r]^{b_*} & {H(C,d_C)} \ar[r]^{c_*} & {H(C'',d_{C''}).}
 \ar@/^1.5pc/[ll]^{ }
}
\end{equation}
\end{lemma}

\proof Consider intervals $I_r = [r;r+\epsilon)$ for $r \in supp(C')$, and $I_r
= (r-\epsilon;r]$ for $r \in supp(C'')$. By \ref{item:gap-one},
\ref{item:gap-two} these are pairwise disjoint, and the distance between any
two of them is $\geq 2\epsilon$. From \ref{item:gap-three} one sees that
$supp(C)$ is contained in the union of these intervals, which shows that $D =
C' \oplus C \oplus C''$ has gap $[\epsilon;2\epsilon)$. Consider the
differential $d_D = \delta + (d_D-\delta)$,
\[
d_D =
\begin{pmatrix}
d_{C'} & 0 & 0 \\ b & d_C & 0 \\ h & c & d_{C''}
\end{pmatrix}, \qquad
\delta =
\begin{pmatrix}
0 & 0 & 0 \\ \beta & 0 & 0 \\ 0 & \gamma & 0
\end{pmatrix}.
\]
We know that $d_{C'}$, $d_C$, $d_{C''}$, $(b-\beta)$, $(c-\gamma)$ are of order
$[2\epsilon;\infty)$. Combining \ref{item:gap-two} with \ref{item:gap-four}
shows that $h$ is of order $[4\epsilon;\infty)$, so that $(d_D-\delta)$ is of
order $[2\epsilon;\infty)$. On the other hand $\delta$ is of order
$[0;\epsilon)$, and \ref{item:gap-three} says that $H(D,\delta) = 0$. Lemma
\ref{th:spectral-sequence} shows that $H(D,d_D) = 0$, which by a general fact
implies the existence of a long exact sequence \eqref{eq:induced-sequence}.
\qed

\newpage

\section{Wrapping it up\label{ch:three}}

Technically, the proof of the exact sequence is an application of Lemma
\ref{th:low-energy}. The spaces $C'$, $C$, $C''$ which occur in the lemma will
be Floer cochain spaces, and the maps $b,c$ relative invariants. While the
definition of $b$ is rather straightforward, that of $c$ uses some of the
geometry from Chapter \ref{ch:one}, specifically the standard fibrations
constructed in Section \ref{sec:model}. In both cases, the desired properties
follow from nonnegative curvature, that is to say Proposition
\ref{th:nonnegative-one} and related results. The homotopy $h$ is obtained by
comparing two different constructions of a particular exact Lefschetz
fibration. Nonnegative curvature again plays a role in analyzing it, but the
vanishing theorem of Section \ref{sec:vanishing} is also important.

\subsection{Preliminaries\label{sec:setup}}

This section sets up the framework for the whole chapter. The data are: $M$ is
an exact symplectic manifold. $L_0,L_1 \subset M$ are exact Lagrangian
submanifolds. $L \subset M$ is an exact Lagrangian sphere, which comes with a
diffeomorphism $f: S^n \rightarrow L$ and a symplectic embedding $\iota:
T(\lambda) \rightarrow M$ for some $\lambda>0$, such that $\iota|T(0) = f$. We
use the Dehn twist $\tau_L$ defined using $\iota$ and a function $R$. Also
given are small constants $\epsilon,\delta>0$. The following conditions are
required to hold:
\begin{Romanlist}
\item \label{cond:transverse}
$L \cap L_0$, $L \cap L_1$, $L_0 \cap L_1$ are transverse intersections, and $L
\cap L_0 \cap L_1 = \emptyset$.
\item \label{cond:action}
The actions $a_{L_0,L_1}(x)$ of distinct points $x \in L_0 \cap L_1$ differ by
at least $3\epsilon$. Secondly, as $({x}_0,x_1)$ runs over $(L_0 \cap L) \times
(L \cap L_1)$, the numbers $a_{L_0,L}({x}_0) + a_{L,L_1}(x_1)$ differ pairwise
by at least $3\epsilon$. Thirdly, for any $x \in L_0 \cap L_1$, $({x}_0,x_1)
\in (L_0 \cap L) \times (L \cap L_1)$ one has
\[
|a_{L_0,L_1}(x) - a_{L_0,L}({x}_0) - a_{L,L_1}(x_1)| \geq 5 \epsilon.
\]
\item \label{cond:distance}
For all $y_k \in f^{-1}(L \cap L_k)$, $k = 0,1$, the distance $dist({y}_0,y_1)$
in the standard metric on $S^n$ is $\geq 2\pi\delta$.
\item \label{cond:local}
$\iota^*\theta = \theta_T|T(\lambda)$ is the standard one-form, and the
function $K_L$ associated to $L$ is zero. Moreover, each $\iota^{-1}(L_k)
\subset T(\lambda)$ is a union of fibres; one can write this as
\begin{equation} \label{eq:union-of-fibres}
\iota^{-1}(L_k) = \bigcup_{y \in \iota^{-1}(L \cap L_k)} T(\lambda)_y.
\end{equation}
\item \label{cond:wobbly}
$R$ satisfies $0 \geq 2\pi R(0) > -\epsilon$, and is such that $\tau_L$ is
$\delta$-wobbly.
\end{Romanlist}

\begin{remark} \label{re:scope}
Since we will establish the exact sequence under these conditions, it is
necessary to convince ourselves that they do not restrict its validity in any
way. Suppose then that we are given arbitrary exact Lagrangian submanifolds
$L_0,L_1$ and a framed exact Lagrangian sphere $(L,[f])$ in $M$. After
perturbing the submanifolds slightly, one can assume that \ref{cond:transverse}
holds. Another such perturbation achieves \ref{cond:action} for some
$\epsilon>0$. This is an instance of a general fact: by moving one of two
transverse exact Lagrangian submanifolds slightly, the action of the
intersection points can be changed independently of each other, by arbitrary
sufficiently small amounts. Choose some representative $f$ of the framing.
Since $L_0 \cap L_1 \cap L = \emptyset$, \ref{cond:distance} is automatically
true for some $\delta>0$. Because the intersections $L \cap L_k$ are
transverse, one can find a symplectic embedding $\iota: T(\lambda) \rightarrow
M$, for some $\lambda>0$, which extends $f$ and such that $\iota^{-1}(L \cap
L_k)$ is a union of fibres. By replacing $\theta$ with $\theta + dH$ for a
suitable $H$, and making $\lambda$ smaller, one can ensure that $\iota^*\theta
= \theta_T$ is satisfied. Note that when one modifies $\theta$ in this way, the
functions associated to exact Lagrangian submanifolds change accordingly. One
can use this and the freedom in the choice of $H$ to arrange that $K_L$ becomes
equal to zero. In any case, the values of the action functional at intersection
points remain the same, so that this does not interfere with \ref{cond:action}.
We have now satisfied \ref{cond:local}. It is no problem to choose $R$ such
that \ref{cond:wobbly} holds for the previously obtained $\epsilon, \delta$.
None of the changes which we have made affects Floer cohomology. Therefore,
once the exact sequence is established for the modified data, it also holds for
the original ones.
\end{remark}

Next, we need to draw some elementary inferences. \ref{cond:local} implies that
$dK_{L_k} = \theta|L_k$ vanishes on $L_k \cap \im(\iota)$; in other words
$K_{L_k} \circ \iota$ is constant on each fibre in \eqref{eq:union-of-fibres}.
Let $\tau$ be the model Dehn twist from which $\tau_L$ is constructed. The
function $K_\tau$ associated to it was determined in Lemma
\ref{th:on-invariant}. Now $K_{\tau_L}$ vanishes outside $\im(\iota)$ and,
again by \ref{cond:local}, satisfies $K_{\tau_L} \circ \iota = K_\tau$.
Concretely
\begin{equation} \label{eq:kk-function}
\begin{aligned}
 K_{\tau_L}(\iota(y))
 &= 2\pi( R'(||y||)||y|| - R(||y||)) \\
 &= -2 \pi R(0) + 2\pi \int_0^{||y||} (R'(||y||) - R'(t))\, dt.
\end{aligned}
\end{equation}
In particular $K_{\tau_L}|L = -2\pi R(0)$, which by \ref{cond:wobbly} lies in
$[0;\epsilon)$. The same condition says that $R'$ decreases monotonically from
$R'(0) = 1/2$ until it reaches the value $\delta$, and thereafter takes values
in $[0;\delta)$. By combining this with \eqref{eq:kk-function} one obtains the
estimate, valid for all $y \in T(\lambda)$ with $R'(||y||) \geq \delta$,
\begin{equation} \label{eq:action-estimate}
-2\pi R(0) \geq K_{\tau_L}(\iota(y)) \geq -2\pi R(0) - 2\pi \int_0^{\infty}
R'(t) dt = 0.
\end{equation}
Now consider the $\R$-graded vector spaces
\begin{align*}
 & C' = CF(L,L_1) \otimes CF(\tau_L(L_0),L), \\
 & C'' = CF(L_0,L_1).
\end{align*}
The first part of \ref{cond:action} implies that $C''$ has gap $(0;3\epsilon)$.
Clearly, a point $\tilde{x}_0$ lies in $\tau_L(L_0) \cap L$ iff $x_0 =
\tau_L^{-1}(\tilde{x}_0)$ lies in $L_0 \cap L_1$. By definition of
$K_{\tau_L(L_0)}$ and the computation above,
\begin{equation} \label{eq:tau-acts}
 a_{\tau_L(L_0),L}(\tilde{x}_0) = a_{L_0,L}(x_0) - K_{\tau_L}(x_0) =
 a_{L_0,L}(x_0) + 2\pi R(0).
\end{equation}
Hence $C'$ can be identified with $CF(L,L_1) \otimes CF(L_0,L)$ up to a shift
in the grading, which is by a constant of size $<\epsilon$. It therefore
follows from \ref{cond:action} that $C'$ has gap $(0;3\epsilon)$, and that the
distance between the supports of $C',C''$ is at least $4\epsilon$. To
summarize, what we have shown is that $C',C''$ satisfy the assumptions
\ref{item:gap-one}, \ref{item:gap-two} of Lemma \ref{th:low-energy}.

\begin{lemma} \label{th:new-intersection-points}
$\tau_L(L_0)$, $L_1$ intersect transversally, and there are injective maps
\begin{align*}
 & p: (\tau_L(L_0) \cap L) \times (L \cap L_1)
 \longrightarrow \tau_L(L_0) \cap L_1, \\
 & q: L_0 \cap L_1 \longrightarrow \tau_L(L_0) \cap L_1
\end{align*}
such that $\tau_L(L_0) \cap L_1$ is the disjoint union of their images. These
maps have the following properties:
\begin{romanlist}
\item \label{item:q-map}
$q$ is the inclusion $q(x) = x$. It preserves the values of the action
functional, $a_{\tau_L(L_0),L_1}(x) = a_{L_0,L_1}(x)$. Moreover, for any $w \in
\tau_L(L_0) \cap L_1$ and $x \in L_0 \cap L_1$ with $w \neq q(x)$ one has
$a_{L_0,L_1}(x) - a_{\tau_L(L_0),L_1}(w) \notin [0;3\epsilon)$.

\item \label{item:p-map}
Set $\tilde{x} = p(\tilde{x}_0,x_1)$. Then
\begin{equation} \label{eq:epsilon-change}
 0 \leq a_{\tau_L(L_0),L_1}(\tilde{x}) -
 a_{\tau_L(L_0),L}(\tilde{x}_0) - a_{L,L_1}(x_1)
 < \epsilon.
\end{equation}
Moreover, for any $w \in \tau_L(L_0) \cap L_1$ and $(\tilde{x}_0,x_1) \in
(\tau_L(L_0) \cap L) \times (L \cap L_1)$ with $w \neq p(\tilde{x}_0,x_1)$ one
has $a_{\tau_L(L_0),L_1}(w) - a_{\tau_L(L_0),L}(\tilde{x}_0) - a_{L,L_1}(x_1)
\notin [0;3\epsilon)$.

\item \label{item:p-map-two}
Suppose that there are $x_k \in L \cap L_k$, $k = 0,1$, whose preimages $y_k =
\iota^{-1}(x_k)$ are antipodes on $S^n$. Since $\tau|S^n$ is the antipodal map,
$\tilde{x}_0 = \tau_L(x_0)$ is equal to $x_1$ (hence $x_1 \in \tau_L(L_0) \cap
L \cap L_1$, and these are all such triple intersection points). In that case
$p(\tilde{x}_0,x_1) = \tilde{x}_0 = x_1$, and
$a_{\tau_L(L_0),L_1}(p(\tilde{x}_0,x_1)) = a_{\tau_L(L_0),L}(\tilde{x}_0) +
a_{L,L_1}(x_1)$.
\end{romanlist}
\end{lemma}

\proof \ref{cond:transverse} and \ref{cond:local} imply that $L_0 \cap L_1 \cap
\im(\iota) = \emptyset$. Since $\tau_L$ is the identity outside $\im(\iota)$,
one has $L_0 \cap L_1 = (\tau_L(L_0) \cap L_1) \setminus \im(\iota)$, so that
$q$ can indeed be defined to be the inclusion. The equality
$a_{\tau_L(L_0),L_1}(x) = a_{L_0,L_1}(x)$ follows from the fact that
$K_{\tau_L}$ vanishes outside $\im(\iota)$.

There is a bijective correspondence between pairs $(\tilde{x}_0,x_1) \in
(\tau_L(L_0) \cap L) \times (L \cap L_1)$ and $(y_0,y_1) \in \iota^{-1}(L_0
\cap L) \times \iota^{-1}(L \cap L_1)$, given by setting $y_0 =
\iota^{-1}(\tau_L^{-1}(\tilde{x}_0))$, $y_1 = \iota^{-1}(x_1)$. As a
consequence of \ref{cond:local},
\begin{equation} \label{eq:decompose-intersection}
\iota^{-1}(\tau_L(L_0) \cap L_1) = \bigcup_{y_0,y_1} \tau(T(\lambda)_{y_0})
\cap T(\lambda)_{y_1}.
\end{equation}
Since $\tau$ is $\delta$-wobbly \ref{cond:wobbly} and $dist(y_0,y_1) \geq
2\pi\delta$ \ref{cond:distance}, one can apply Lemma
\ref{th:local-intersections} which tells us that each subset on the right hand
side of \eqref{eq:decompose-intersection} consists of exactly one point. Fix
temporarily some $(y_0,y_1)$ and write $\tau(T(\lambda)_{y_0}) \cap
T(\lambda)_{y_1} = \{\tilde{y}\}$, $\tilde{x} = \iota(\tilde{y})$. One defines
$p(\tilde{x}_0,x_1) = \tilde{x}$. Then
\begin{align*}
 a_{\tau_L(L_0),L_1}(\tilde{x})
 & = K_{L_1}(\tilde{x}) - K_{\tau_L(L_0)}(\tilde{x}) \\
 & = K_{L_1}(\tilde{x}) - K_{L_0}(\tau_L^{-1}(\tilde{x}))
 - K_{\tau_L}(\tau_L^{-1}(\tilde{x})).
\end{align*}
By construction $\tilde{y}$ lies in the fibre $T(\lambda)_{y_1}$, and since
$K_{L_1} \circ \iota$ is constant on fibres, $K_{L_1}(\tilde{x}) =
K_{L_1}(x_1)$. The same reasoning shows that $K_{L_0}(\tau_L^{-1}(\tilde{x})) =
K_{L_0}(x_0)$, where $x_0 = \tau_L^{-1}(\tilde{x}_0)$. Moreover, one sees from
\eqref{eq:kk-function} that $K_{\tau_L}$ is invariant under $\tau_L$, so that
$K_{\tau_L}(\tau_L^{-1}(\tilde{x})) = K_{\tau_L}(\tilde{x})$. With this and
\eqref{eq:tau-acts} in mind, one continues the computation
\begin{equation} \label{eq:finish-action}
\begin{aligned}
 a_{\tau_L(L_0),L_1}(\tilde{x})
 & = K_{L_1}(x_1) - K_{L_0}(x_0) - K_{\tau_L}(\tilde{x}) \\
 & = a_{L,L_1}(x_1) + a_{L_0,L}(x_0) - K_{\tau_L}(\tilde{x}) \\
 & = a_{L,L_1}(x_1) + a_{\tau_L(L_0),L}(\tilde{x}_0) - K_{\tau_L}(\iota(\tilde{y}))
 - 2\pi R(0).
\end{aligned}
\end{equation}
Lemma \ref{th:local-intersections} also says that $R'(||\tilde{y}||) \geq
\delta$. Combining this with \eqref{eq:action-estimate} and \ref{cond:wobbly}
shows that $- K_{\tau_L}(\iota(\tilde{y})) - 2\pi R(0)$ lies in $[0;\epsilon)$,
which completes our proof of \eqref{eq:epsilon-change}.

It is clear from their definitions that $p,q$ are injective. A point of
$\tau_L(L_0) \cap L_1$ falls into $\im(q)$ or $\im(p)$ depending on whether it
lies inside or outside $\im(\iota)$, hence the two images are disjoint and
cover $\tau_L(L_0) \cap L_1$. The transversality follows from Lemma
\ref{th:local-intersections} for $\im(p)$ and from that of $L_0 \cap L_1$ for
$\im(q)$. We now turn to the claim made in the last sentence of
\ref{item:q-map}. Supposing that $w$ is a point of $L_0 \cap L_1$ different
from $x$, one has $|a_{\tau_L(L_0),L_1}(w) - a_{L_0,L_1}(x)| = |a_{L_0,L_1}(w)
- a_{L_0,L_1}(x)| \geq 3\epsilon$ by \ref{cond:action}. In the remaining case,
which is when $w = p(\tilde{x}_0,x_1)$, \eqref{eq:epsilon-change} shows that
$|a_{\tau_L(L_0),L_1}(w) - a_{L_0,L_1}(x)| > |a_{\tau_L(L_0),L}(\tilde{x}_0) +
a_{L,L_1}(x_1) - a_{L_0,L_1}(x)| - \epsilon$. We already know that the supports
of $C',C''$ are at least $4\epsilon$ apart, and one concludes that
$|a_{\tau_L(L_0),L_1}(w) - a_{L_0,L_1}(x)| > 3\epsilon$. A similar argument,
paying a little more attention to signs, proves the parallel statement in
\ref{item:p-map}. Finally, the only non-obvious things in \ref{item:p-map-two}
are the fact that $p(\tilde{x}_0,x_1) = x_1$ and the statement about the
action. But these follow from Lemma \ref{th:local-intersections} and the
definition of $p$, respectively from \eqref{eq:finish-action} and $K_{\tau_L}|L
= -2\pi R(0)$. \qed

Lemma \ref{th:new-intersection-points} and \ref{cond:action} imply that the
actions $a_{\tau_L(L_0),L_1}(x)$ of different points $x \in \tau_L(L_0) \cap
L_1$ differ by at least $2\epsilon$. In fact, for two such points which lie in
$\im(q)$ the statement follows from \ref{item:q-map} in the lemma; for two
points which lie in $\im(p)$, from \ref{item:p-map}; and combining the two
parts shows it when one point lies in $\im(p)$ and the other in $\im(q)$. Set
\[
C = CF(\tau_L(L_0),L_1).
\]
We have just seen that this has gap $(0;2\epsilon)$. Define maps $\beta: C'
\rightarrow C$, $\gamma: C \rightarrow C''$ by $\beta(\gen{x_1} \otimes
\gen{\tilde{x}_0}) = \gen{p(\tilde{x}_0,x_1)}$ and
$\gamma(\gen{p(\tilde{x}_0,x_1)}) = 0$, $\gamma(\gen{q(x)}) = \gen{x}$. The
result above shows that these are of order $[0;\epsilon)$ and fit into a short
exact sequence as in Lemma \ref{th:low-energy}\ref{item:gap-three}. What
remains to be done is to realize them as ``low order parts'' of chain maps
$b,c$, and to construct the homotopy $h$.

\subsection{The first map\label{sec:b}}

Let $S$ be the surface in Figure \ref{fig:y-piece}, which has three boundary
components $\partial_kS$ and three strip-like ends (two positive and a negative
one). Take the trivial exact symplectic fibration $\pi : E = S \times M
\rightarrow S$, with $\Omega \in \Omega^2(E)$, $\Theta \in \Omega^1(E)$ pulled
back from $\o$, $\theta$. Equip this with the Lagrangian boundary condition $Q
= (\partial_1 S \times \tau_L(L_0)) \cup (\partial_2 S \times L) \cup
(\partial_3 S \times L_1)$; $\kappa_{Q}$ is zero, and $K_{Q}(z,x)$ is equal to
$K_{\tau_L(L_0)}(x)$, $K_L(x)$ or $K_{L_1}(x)$, for $z$ in the respective
component $\partial_kS$. This gives rise to a relative invariant
\[
 \Phi^{rel}_0(E,\pi,Q) : HF(L,L_1) \otimes HF(\tau_L(L_0),L) \rightarrow
 HF(\tau_L(L_0),L_1).
\]
The purpose of this section is to analyse this more closely, on the cochain
level. Fix a complex structure $j$ on $S$, trivial over the ends, and let $U
\subset S$ be the open set shaded in Figure \ref{fig:y-piece}.

\includefigure{y-piece}{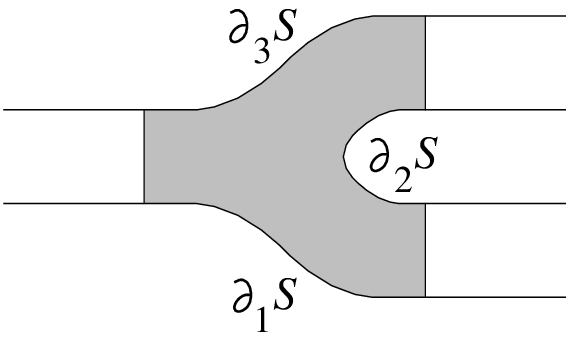}{hb}

\begin{lemma} \label{th:transversality-b}
$(E,\pi,Q)$ and $U$ satisfy the conditions of Proposition
\ref{th:nonnegative-one}(i).
\end{lemma}

\proof The curvature of $(E,\pi)$ is zero. A section $u(z) = (z,\sigma(z))$ is
horizontal iff $\sigma(z) \equiv x \in M$ is constant, and it further satisfies
$u(\partial S) \subset Q$ iff $x \in \tau_L(L_0) \cap L \cap L_1$. If $W
\subset S$ is a connected open subset which intersects all three boundary
components, the same description applies to partial horizontal sections $w: W
\rightarrow E|W$ with $w(W \cap \partial S) \subset Q$. As a consequence, any
such section can be extended to $u \in \MM^h$. This is in particular true for
$W = U$.

Fix some $x \in \tau_L(L_0) \cap L \cap L_1$ and the corresponding constant
section $u \in \MM^h$. By definition $A(u) = 0$, and it remains to prove that
$\ind\,D_{u,J} = 0$. From the description of the points $x$ in Lemma
\ref{th:new-intersection-points}\ref{item:p-map-two}, together with the
corresponding local statement in Lemma \ref{th:local-intersections}, it follows
that there is a symplectic isomorphism $TM_x \iso \C^n$ which takes the tangent
spaces to $\tau_L(L_0)$, $L$, and $L_1$ to $\R^n$, $e^{2 \pi i/3}\R^n$, and
$e^{\pi i/3}\R^n$ respectively. Since the index is independent of the almost
complex structure, we may choose $J = j \times J^M$ to be the product of $j$
and some $\o$-compatible $J^M$ on $M$. We may also assume that the isomorphism
$TM_x \iso \C^n$ takes $J^M_x$ to the standard complex structure. By
definition, the domain of $D_{u,J}$ are sections of the vector bundle
$u^*(TE^v,J|TE^v) \iso S \times TM_x$, with boundary conditions given by the
tangent spaces to $\tau_L(L_0)$, $L$, $L_1$. Its range are $(0,1)$-forms with
values in the same vector bundle. The preceding discussion allows us to
identify
\begin{equation} \label{eq:triangle-conditions}
\begin{split}
 & \WW^1_u = \{ X \in W^{1,p}(S,\C^n) \suchthat
 X_z \in e^{i(1-k)\pi/3}\R^n \text{ for $z \in \partial_kS$} \}, \\
 & \WW^0_{u,J} = L^p(\Lambda^{0,1}S \otimes \C^n).
\end{split}
\end{equation}
In \eqref{eq:connection-formula} take $\nabla = \nabla^{S} \times \nabla^M$ to
be the product of torsion-free connections on $S$ and on $M$. Then the second
term $(\nabla_X J) \circ Du \circ j$ in the formula vanishes, because $Du \circ
j$ takes values in $TE^h$ whereas $\nabla_X J$ is nontrivial only on $TE^v$;
and moreover, the pullback connection $u^*\nabla$ on $u^*TE^v$ is trivial. This
shows that $D_{u,J}$ is the standard $\bar\partial$-operator for functions $S
\rightarrow \C^n$, with boundary conditions \eqref{eq:triangle-conditions}.
\includefigure{sausage}{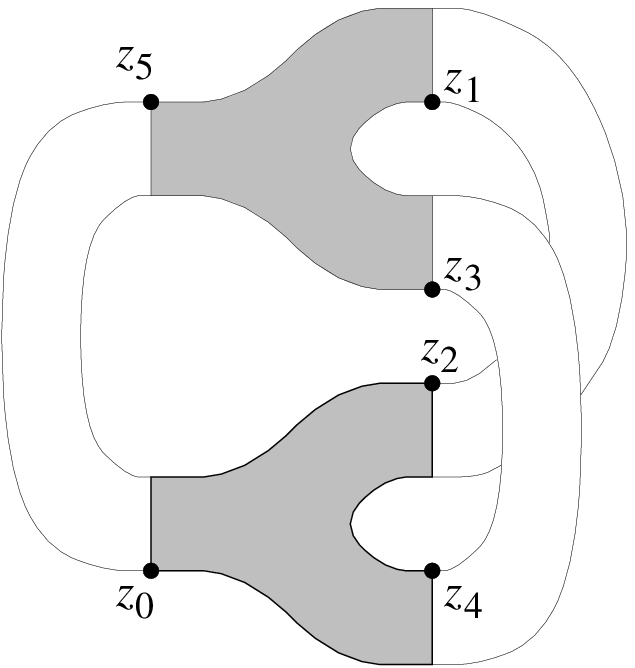}{hb}%

There is a general index formula for such operators, but we prefer to use an ad
hoc gluing trick instead. Consider the compact surface $\bar{S}$ in Figure
\ref{fig:sausage}, which is of genus one with one boundary component.
Parametrize the boundary by a closed path $l : [0;6] \rightarrow
\partial\bar{S}$, such that $l(t) = z_t$ for $t \in \{0,1,2,3,4,5\}$ are the
marked points in Figure \ref{fig:sausage}. Take a smooth nondecreasing function
$\lambda: [0;6] \rightarrow \R$ such that
\[
 \lambda(t) = \begin{cases}
 0 & 0 \leq t \leq 1, \\
 1/3 & 2 \leq t \leq 3, \\
 2/3 & 4 \leq t \leq 5. \\
 1 & t = 6.
 \end{cases}
\]
Let $\bar{D}$ be the $\bar\partial$-operator on the trivial bundle $\bar{S}
\times \C^n \rightarrow \bar{S}$, with boundary condition given by the family
of Lagrangian subspaces $\Lambda_{l(t)} = e^{\pi i \lambda(t)} \R^n \subset
\C^n$. As a loop in the Lagrangian Grassmannian, this represents $n$ times the
standard generator of the fundamental group. Riemann-Roch for compact surfaces
with boundary therefore tells us that $\ind\,\bar{D} = n\chi(\bar{S}) + n = 0$.
On the other hand, one can divide $\bar{S}$ into five pieces (two shaded and
three unshaded ones) as indicated in Figure \ref{fig:sausage}, and add
strip-like ends to each piece; the standard gluing theory for elliptic
operators says that $\ind\,\bar{D}$ is the sum of the indices of the obvious
corresponding operators on those pieces. For each shaded piece, this yields a
copy of $D_{u,J}$ (in one of the two cases, the vector space $\C^n$ has been
rotated by $e^{\pi i/3}$). The unshaded pieces give rise to operators of index
zero. This can be derived from the index theorem of \cite{robbin-salamon93}, or
else by directly deforming the operator to an invertible one. Consider for
instance the two intervals of $\partial\bar{S}$ which are contained in the
leftmost unshaded piece. The Lagrangian subspaces $\Lambda_z$ parametrized by
the points $z$ in one of these intervals are all equal to $e^{\pi i/3}\R^n$;
for the other interval, they are of the form $e^{\pi i s}\R^n$ for $2/3 \leq s
\leq 1$. Since $e^{\pi i/3} \R^n \cap e^{\pi i s}\R^n = 0$ for all such $s$,
the Maslov index for paths \cite{robbin-salamon93} is zero. The same holds for
the other unshaded pieces, and one concludes that $0 = \ind\, \bar{D} = 2\,
\ind \, D_{u,J}$. \qed

At this point, fix
\begin{align*}
 & J^{(1)} \in \JJ^{reg}(M,\tau_L(L_0),L), \\
 & J^{(2)} \in \JJ^{reg}(M,L,L_1), \\
 & J^{(3)} \in \JJ^{reg}(M,\tau_L(L_0),L_1).
\end{align*}
Proposition \ref{th:nonnegative-one}(i) tells us that there is a $J \in
\JJ(E,\pi,j,J^{(1)},J^{(2)},J^{(3)})$ which is both horizontal and regular. By
part (ii) of the same result, the coefficients $\nu_J(\tilde{x}_0,x_1,x)$ for
$\tilde{x}_0 \in \tau_L(L_0) \cap L$, $x_1 \in L \cap L_1$, $x \in \tau_L(L_0)
\cap L_1$, are zero whenever
\begin{equation} \label{eq:action-difference}
\chi(\tilde{x}_0,x_1,x) = a_{\tau_L(L_0),L_1}(x) -
a_{\tau_L(L_0),L}(\tilde{x}_0) - a_{L,L_1}(x_1)
\end{equation}
is $<0$. Lemma \ref{th:new-intersection-points}\ref{item:p-map} shows that
$\chi(\tilde{x}_0,x_1,x)$, if $\geq 0$, must be either in $[0;\epsilon)$ or
$[3\epsilon;\infty)$, and that the first case only happens for $x =
p(\tilde{x}_0,x_1)$. Hence, if one writes the relative invariant as
\begin{equation} \label{eq:split-one}
\begin{split}
 & C\Phi_0^{rel}(E,\pi,Q,J) = \phi + (C\Phi_0^{rel}(E,\pi,Q,J)-\phi), \\
 & \phi(\gen{\tilde{x}_0} \otimes \gen{x_1}) =
 \nu_J(\tilde{x}_0,x_1,p(\tilde{x}_0,x_1)) \gen{p(\tilde{x}_0,x_1)}
\end{split}
\end{equation}
then $\phi$ is of order $[0;\epsilon)$, while the second term is of order
$[3\epsilon;\infty)$. The rest of this section contains the proof of the
following result, which determines the low order part:

\begin{prop} \label{th:low-nu}
$\nu_J(\tilde{x}_0,x_1,p(\tilde{x}_0,x_1)) = 1$ for all $(\tilde{x}_0,x_1)$.
\end{prop}

There is one particular case of this which follows from the previous
considerations. Namely, suppose that there is a pair $(\tilde{x}_0,x_1)$ with
$\tilde{x}_0 = x_1 = x \in \tau_L(L_0) \cap L \cap L_1$. In that case $p(x,x) =
x$ and $\chi(x,x,x) = 0$ by Lemma
\ref{th:new-intersection-points}\ref{item:p-map-two}, and therefore
$\nu_J(x,x,x) = \# \MM^h(x,x,x)$ by Proposition \ref{th:nonnegative-one}(ii).
Now $\# \MM^h(x,x,x) = 1$ because, as we saw already when proving Lemma
\ref{th:transversality-b}, the unique element of that space is the constant
section $u(z) = (z,x)$. Our strategy will be to reduce the computation of all
the $\nu_J(\tilde{x}_0,x_1,p(\tilde{x}_0,x_1))$ to this case, by using suitable
deformations.

What ``deformation'' means here is keeping the submanifolds $L,L_0,L_1$ and the
constant $\epsilon$ fixed, while changing the remaining data. More precisely,
for $0 \leq \mu \leq 1$ we consider the following: a smooth family $f^\mu$ of
diffeomorphisms $S^n \rightarrow L$; smoothly varying positive numbers
$\lambda^\mu$, and symplectic embeddings $\iota^\mu: T(\lambda^\mu) \rightarrow
M$ with $\iota^\mu|T(0) = f^\mu$; the Dehn twists $\tau_L^\mu$ defined using
$\iota^\mu$ and functions $R^\mu$ supported in $(-\infty;\lambda^\mu)$; and
constants $\delta^\mu$. These should agree with the given data
$f,\lambda,\iota,R,\tau_L,\delta$ for $\mu = 0$. We also require that the
analogues of \ref{cond:distance}--\ref{cond:wobbly} continue to hold for any
$\mu$; where in \ref{cond:wobbly} we take the original $\epsilon$ throughout.

Choose $\tilde{x}_0 \in \tau_L(L_0) \cap L$, $x_1 \in L \cap L_1$. Given a
``deformation'' in the sense which we have just explained, one can set
$\tilde{x}_0^\mu = \tau_L^\mu(\tau_L)^{-1}(\tilde{x}_0) \in \tau_L^\mu(L_0)
\cap L$ and $x_1^\mu = x_1 \in L \cap L_1$, which ``continues'' the points
smoothly into the deformed situation. We claim that $x = p(\tilde{x}_0,x_1)$
fits similarly into a smooth family $x^\mu \in \tau_L^\mu(L_0) \cap L_1$. The
point is that since \ref{cond:transverse}--\ref{cond:wobbly} continue to hold,
Lemma \ref{th:new-intersection-points} can be applied to the situation for any
$\mu$. This ensures that the intersections $\tau_L^\mu(L_0) \cap L_1$ remain
transverse, which implies that a unique family $x^\mu$ exists. In fact, it even
provides a smooth family of injective maps $p^\mu: (\tau_L^\mu(L_0) \cap L)
\times (L \cap L_1) \rightarrow \tau_L^\mu(L_0) \cap L_1$, such that $x^\mu =
p^\mu(\tilde{x}_0^\mu,x_1^\mu)$.

\begin{lemma} \label{th:antipodes}
For any $(\tilde{x}_0,x_1)$ there is a ``deformation'' such that $\tilde{x}_0^1
= x_1^1$.
\end{lemma}

\proof From \ref{cond:distance} we know that $y_0 =
f^{-1}(\tau_L^{-1}(\tilde{x}_0))$, $y_1 = f^{-1}(x_1)$ are points on $S^n$
whose distance is $\geq 2\pi \delta$. Let $g^\mu \in \Diff(S^n)$ be an isotopy,
$g^0 = \id$, such that $g^1(y_0)$, $g^1(y_1)$ are antipodes. There are $C^\mu
\geq 1$, smoothly depending on $\mu$ and with $C^0 = 1$, with the property that
\[
C^\mu \geq ||D(g^\mu)_y||, \; ||D(g^\mu)^{-1}_y|| \quad \text{for all $y \in
S^n$.}
\]
Consider the ``deformation'' $f^\mu = f \circ (g^\mu)^{-1}$, $\lambda^\mu =
\lambda/C^\mu$, $\iota^\mu = \iota \circ G^\mu | T(\lambda^\mu)$, $\delta^\mu =
\delta/C^\mu$; here $G^\mu \in \Sympe(T)$ is induced by $g^\mu$, in the sense
that $G^\mu|T(0) = (g^\mu)^{-1}$. The bound on $D(g^\mu)$ implies that $G^\mu$
maps $T(\lambda^\mu)$ to $T(\lambda)$, so that $\iota^\mu$ is well-defined. On
the other hand, because of the bound on $D(g^\mu)^{-1}$, the distance between
any point of $(f^\mu)^{-1}(L \cap L_0) = g^\mu f^{-1}(L \cap L_0)$ and any
point of $(f^\mu)^{-1}(L \cap L_1) = g^\mu f^{-1}(L \cap L_1)$ is $\geq
2\pi\delta^\mu$, which shows that \ref{cond:distance} holds during the
deformation. For \ref{cond:local} it is sufficient to note that $G^\mu$ takes
the canonical one-form $\theta_T$ to itself and maps fibres to fibres. And it
is no problem to find functions $R^\mu$ which satisfy \ref{cond:wobbly} with
the $\delta^\mu$ defined above and the given $\epsilon$. By definition
\[
 \tau_L|L = f \circ A \circ f^{-1}, \quad
 \tau_L^\mu|L = f^\mu \circ A \circ (f^\mu)^{-1}
 = f \circ (g^\mu)^{-1} \circ A \circ g^\mu \circ f^{-1}.
\]
Since $g^1(y_1) = A(g^1(y_0))$ by construction, one sees that
\[
\tilde{x}_0^1 = \tau_L^1 (\tau_L^{})^{-1}(\tilde{x}_0)
 = f \circ (g^1)^{-1} \circ A \circ g^1 (y_0) = f(y^1)
 = x_1 = x_1^1. \qed
\]

Given a ``deformation'', one can repeat the construction at the beginning of
this section in a parametrized sense, which means equipping the trivial exact
symplectic symplectic fibration $(E,\pi)$ with a smooth family $Q^\mu$ of
Lagrangian boundary conditions, modelled over the ends on the pairs of exact
Lagrangian submanifolds $(\tau_L^\mu(L_0),L)$, $(L,L_1)$, and
$(\tau_L^\mu(L_0),L_1)$. Since their intersections are transverse for all
$\mu$, it makes sense to consider parametrized moduli spaces of
pseudo-holomorphic sections. To do that, take smooth families of almost complex
structures $J^{(1),\mu}, J^{(2),\mu}, J^{(3),\mu} \in \JJ(M)$ which reduce to
the previously chosen ones for $\mu = 0$, and similarly a family $J^\mu \in
\JJ^h(E,\pi,j,J^{(1),\mu},J^{(2),\mu},J^{(3),\mu})$. The parametrized moduli
space we are interested in is
\[
\MM^{para} = \bigcup_{\mu}\, \{\mu\} \times
\MM_{J^\mu}(\tilde{x}_0^\mu,x_1^\mu,x^\mu)
\]
for points $\tilde{x}_0^\mu,x_1^\mu,x^\mu$ as introduced above.

\begin{lemma} \label{th:no-bubbles}
$\MM^{para}$ is compact.
\end{lemma}

\proof Consider $CF(L,L_1)$, $CF(\tau_L^\mu(L_0),L)$,
$CF(\tau_L^\mu(L_0),L_1)$. We already know that the first of these $\R$-graded
vector spaces has gap $(0;3\epsilon)$; for the second one, the same is true
because it agrees with $CF(L_0,L)$ up to a constant shift in the grading; and
by repeating the considerations after Lemma \ref{th:new-intersection-points},
one can show that $CF(\tau_L^\mu(L_0),L_1)$ has gap $(0;2\epsilon)$. As another
consequence of Lemma \ref{th:new-intersection-points}\ref{item:p-map} in the
deformed situation, any $(\mu,u) \in \MM^{para}$ satisfies
\[
 A(u) = a_{\tau_L^\mu(L_0),L_1}(x^\mu) -
 a_{\tau_L^r(L_0),L}(\tilde{x}_0^\mu) - a_{L,L_1}(x_1^\mu) \in [0;\epsilon).
\]
One can now argue exactly as in Lemma \ref{th:energy-compactness}; there are no
points at infinity in the parametrized Gromov-Floer compactification, since the
principal component of any such point would have negative action, which would
contradict the fact that $(E,\pi)$ has nonnegative curvature. \qed

Suppose now that the $J^{(k),\mu}$ and $J^\mu$ for $\mu = 1$ have been chosen
regular; for $J^1$ this can be achieved without leaving the class of horizontal
almost complex structure, because Lemma \ref{th:transversality-b} equally
applies to the deformed situation for $\mu = 1$. By using a parametrized
version of the same lemma and of Proposition \ref{th:nonnegative-one}(i), one
sees that in addition, the family $(J^\mu)$ can be chosen to be regular in the
parametrized sense. Then the one-dimensional part of $\MM^{para}$ is a compact
one-manifold, and its boundary points are precisely those with $\mu = 0$ or
$1$. One concludes that
\[
\nu_J(\tilde{x}_0,x_1,x) = \nu_{J^1}(\tilde{x}_0^0,x_1^1,x^1).
\]
If the ``deformation'' is as in Lemma \ref{th:antipodes}, the situation for $r
= 1$ is exactly the special case which we have already discussed, so that
$\nu_{J^1}(\tilde{x}_0^0,x_1^1,x^1) = 1$. Since such deformations exist for all
$(\tilde{x}_0,x_1)$, Proposition \ref{th:low-nu} is proved.

From now on, we fix some $J = J^{(4)} \in
\JJ^{reg,h}(E,\pi,Q,j,J^{(1)},J^{(2)},J^{(3)})$ and write $b =
C\Phi_0^{rel}(E,\pi,Q,J^{(4)}): C' \rightarrow C$ for the relative invariant on
the cochain level. Let $\beta$ be the map defined at the end of Section
\ref{sec:setup}. What \eqref{eq:split-one} and Proposition \ref{th:low-nu} say
is that $b = \beta + (b - \beta)$, with $b-\beta$ of order
$[3\epsilon;\infty)$, which is even slightly more than required by Lemma
\ref{th:low-energy}.

\includefigure{triangle}{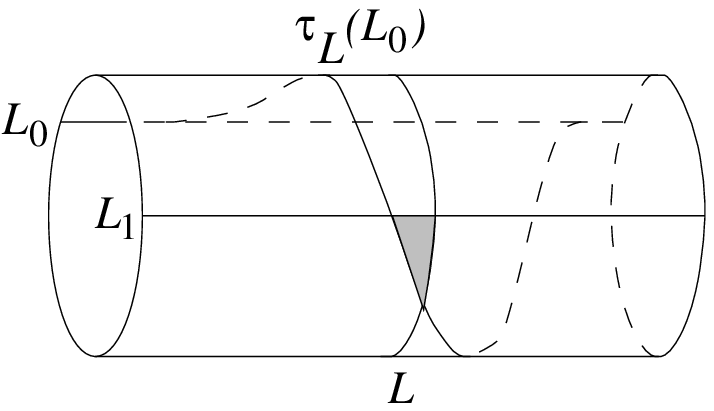}{ht}%
\begin{remark}
Our relative invariant is the well-known pair-of-pants product, or Donaldson
product; in fact, if $u(z) = (z,\sigma(z))$ lies in $\MM_J(\tilde{x}_0,x_1,x)$,
then $\sigma: S \rightarrow M$ is a ``pseudo-holomorphic triangle'' whose sides
lie on $\tau_L(L_0), L, L_1$ and whose vertices are $\tilde{x}_0,x_1,x$.
Proposition \ref{th:low-nu} asserts that there is an odd number of low-area
triangles with certain specified vertices. In the lowest dimension, $n = 1$,
one can see directly that there is precisely one such triangle (Figure
\ref{fig:triangle}). In higher dimensions it is still easy to construct
explicitly the analogue of this particular triangle, but proving that there are
no others seems more difficult. For that reason, we have preferred to take the
indirect approach via ``deformations'', which effectively meant moving
$\tau_L(L_0)$ in such a way that the area of the triangle shrinks to zero.
\end{remark}

\subsection{The second map\label{sec:c}}

From this point onwards, we add one more assumption to those in Section
\ref{sec:setup}:
\begin{Romanlist}
\setcounter{enumi}{5}
\item \label{cond:r}
$R$ is the function $R_r$ which appears in the construction of exact Lefschetz
fibrations in Lemma \ref{th:model-fibrations}, with the given $\lambda$ and
some $0<r<1/2$ which we are free to choose; see more specifically
\eqref{eq:r-function} for this function.
\end{Romanlist}

Since that severely restricts the choice of $R$, one needs to worry about
possible conflicts with the other conditions, or that it might restrict the
ultimate scope of the exact sequence. An inspection of Remark \ref{re:scope}
shows that the only issue is whether, using $R = R_r$, one can satisfy
\ref{cond:wobbly} with arbitrarily small $\delta$ and $\epsilon$. That is taken
care of by Lemma \ref{th:standard-properties}\ref{item:precise-twist}, which
shows that it suffices to choose $r$ small.

\includefigure{c-pasting-ii}{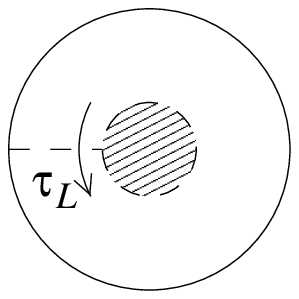}{hb} %
The obvious reason for introducing \ref{cond:r} is that there is a standard
fibration $(E^L,\pi^L)$ over a disc $\bar{D}(r)$, with $r$ small, whose
monodromy around $\partial \bar{D}(r)$ is $\tau_L$. Following Remark
\ref{re:squeeze-boundary}, we want to modify this by a pullback. Take $S^p =
\bar{D}(1/2)$ and a map $p: S^p \rightarrow \bar{D}(r)$ of the form $p(z) =
g(|z|)(z/|z|)$, where $g$ is a function with $g(t) = t$ for small $t$, $g(t) =
r$ for $t \geq r$, and $g'(t) \geq 0$ everywhere. Then
\[
(E^p,\pi^p) = p^*(E^L,\pi^L)
\]
is again an exact Lefschetz fibration. It has nonnegative curvature; this
follows from Lemma \ref{th:standard-properties}\ref{item:nonnegatively-curved}
and the fact that $det(Dp) \geq 0$. Moreover, it is flat on the annulus $S^p
\setminus D(r)$; and using the isomorphism $\phi^p: (E^p)_{1/2} \rightarrow M$
inherited from $\phi^L: (E^L)_r \rightarrow M$, one can identify its monodromy
around $\partial S^p$ with $\tau_L$. To represent this property, we draw
$(E^p,\pi^p)$ as in Figure \ref{fig:c-pasting-ii}.

\includefigure{c-pasting-i}{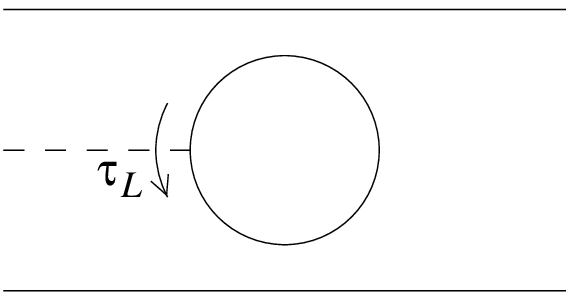}{hb} %
Now take the surface $S^f = (\R \times [-1;1]) \setminus D(1/2) \subset \R^2$,
with coordinates $(s,t)$, and divide it into two parts $S^{f,\pm} = S^f \cap \{
t \in \R^{\pm} \}$, so that $S^{f,+} \cap S^{f,-} = ((-\infty;-1/2] \cup
[1/2;\infty)) \times \{0\}$. Consider trivial fibrations $\pi^{f,\pm}:
E^{f,\pm} = S^{f,\pm} \times M \rightarrow S^{f,\pm}$ over the two parts, and
equip them with differential forms $\Omega^{f,\pm}$, $\Theta^{f,\pm}$ as
follows. $\Omega^{f,\pm}$ is the pullback of $\o \in \Omega^2(M)$, and
similarly $\Theta^{f,+}$ is the pullback of $\theta$; finally $\Theta^{f,-} =
\theta - d(\beta(s)K_{\tau_L})$, where $\beta$ is a function with $\beta(s) =
0$ for $s \geq 1/4$, $\beta(s) = 1$ for $s \leq -1/4$. Define a fibration
$(E^f,\pi^f)$ over $S^f$ by identifying the fibres
\[
E^{f,+}_{(s,0)} \longrightarrow E^{f,-}_{(s,0)}
\]
via $\id_M$ for $s \geq 1/2$, respectively via $\tau_L$ for $s \leq -1/2$.
$\Omega^{f,\pm}$ and $\Theta^{f,\pm}$ match up to forms $\Omega^f$, $\Theta^f$;
the first because $\tau_L$ is symplectic, and the second as a consequence of
our choice of $\Theta^{f,-}$. This makes $(E^f,\pi^f)$ into a (flat) exact
symplectic fibration; it is represented in Figure \ref{fig:c-pasting-i}.

\includefigure{c-pasting-iii}{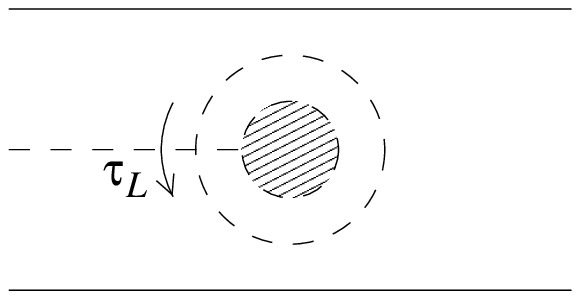}{hb} %
We now carry out a pasting construction of the kind discussed in Section
\ref{sec:basic}. If one identifies the fibres of $E^p$ and $E^f$ at the point
$1/2$ by using $\phi^p: (E^p)_{1/2} \rightarrow M = (E^f)_{1/2}$, then the
monodromies around the circle $|z| = 1/2$ coincide, being both equal to
$\tau_L$. Since the two fibrations are flat close to this circle, one can paste
them together to an exact Lefschetz fibration, denoted by $(E,\pi)$, over $S =
S^p \cup S^f = \R \times [-1;1]$; this is drawn in Figure
\ref{fig:c-pasting-iii}. Equip $(E,\pi)$ with the Lagrangian boundary condition
$Q$ which is the union of $\R \times \{1\} \times L_1 \subset E^{f,+}$ and $\R
\times \{-1\} \times \tau_L(L_0) \subset E^{f,-}$, with $\kappa_Q = 0$, and
with a function $K_Q$ which is $K_{L_1}$ on $\R \times \{1\} \times L_1$ and
$K_{\tau_L(L_0)} - \beta(s)K_{\tau_L}|\tau_L(L_0)$ on $\R \times \{-1\} \times
\tau_L(L_0)$. This is clearly modelled on $(\tau_L(L_0),L_1)$ over the positive
end of $S$; over the negative end it is modelled on $(L_0,L_1)$, as one sees
using the trivialization
\begin{align*}
 & \qquad \xymatrix{
 {(-\infty;-1] \times [-1;1] \times M} \ar[r]^-{\Gamma} \ar[d] & {E} \ar[d]^\pi \\
 {(-\infty;-1] \times [-1;1]} \ar[r]^-{\gamma} & {S}
 } \\
 &
 \gamma = \text{inclusion}, \qquad
 \Gamma(s,t,x) = \begin{cases}
 (s,t,x) \in E^{f,+} & t \geq 0, \\
 (s,t,\tau_L(x)) \in E^{f,-} & t \leq 0.
 \end{cases}
\end{align*}
We now get a relative invariant $\Phi^{rel}_0(E,\pi,Q): HF(\tau_L(L_0),L_1)
\rightarrow HF(L_0,L_1)$. Following the same pattern as in the previous section
(but with considerably less technical difficulties), we need to determine
partially the underlying cochain map.

\begin{lemma} \label{th:c}
$(E,\pi,Q)$ together with $U = (-1;1) \times [-1;1] \subset S$ satisfies the
conditions of Proposition \ref{th:nonnegative-one}(i). Moreover, given $w \in
\tau_L(L_0) \cap L_1$ and $x \in L_0 \cap L_1$, the space $\MM^h(x,w)$ contains
precisely one horizontal section if $w = q(x)$, and is empty otherwise.
\end{lemma}

\proof Because the two parts from which it is assembled have nonnegative
curvature, so does $(E,\pi)$. To any point $x \in L_0 \cap L_1$ one can
associate a horizontal section $u^f: S^f \rightarrow E^f$ satisfying
$u^f(\partial S) \subset Q$, which is defined by
\[
 u^f(s,t) = \begin{cases}
 (s,t,x) \in E^{f,+} & t \geq 0, \\
 (s,t,\tau_L(x)) \in E^{f,-} & t \leq 0.
 \end{cases}
\]
Condition \ref{cond:transverse} and \ref{cond:local} imply that $x \notin
\im(\iota)$. By construction $E^L$ contains a trivial part $\bar{D}(r) \times
(M \setminus \im(\iota))$, see Proposition \ref{th:standard-fibrations}. The
pullback $E^p$ has a corresponding property, which means that there is a
horizontal section $u^p$ of $E^p$, matching up with $u^f$ to form a $u \in
\MM^h(x,x) \subset \MM^h$.

$\int_{S^f} (u^f)^*\Omega^f = 0$ by definition of $\Omega^f$, and similarly
$\int_{S^p} (u^p)^*\Omega^p = 0$ because the image of $u^p$ lies in the trivial
part of $E^p$. It follows that the section we have constructed satisfies $A(u)
= 0$. The connection $\nabla^u$ on $u^*TE^v$ is trivial; in fact, by inspecting
the details of the construction, one can see that $(E,\pi)$ is symplectically
trivial in a neighbourhood of $\im(u)$. By linearization of basic fact that
symplectic parallel transport preserves any Lagrangian boundary condition, one
finds that the subbundle $u^*(TQ \cap TE^v) \subset u^*TE^v|\partial S$ is
preserved under $\nabla^u$; hence it is trivial. To summarize, we have found
that one can identify $u^*(TE^v) \iso S \times \C^n$ symplectically, in such a
way that $\nabla^u$ becomes trivial, and that $u^*(TQ \cap TE^v)$ is mapped to
the subbundle $(\R \times \{-1\} \times \Lambda_{-1}) \cup (\R \times \{1\}
\times \Lambda_1)$ for some Lagrangian subspaces $\Lambda_{\pm 1} \subset
\C^n$; by looking at the positive end, one sees that $\Lambda_{-1}$ and
$\Lambda_1$ are transverse. Then the index formula in terms of the Maslov index
for paths \cite{robbin-salamon93} shows that $\ind\,D_{u,J} = 0$.

Next, suppose that $u \in \MM^h$ is an arbitrary horizontal section satisfying
$u(\partial S) \subset Q$. When restricted to $S^{f,\pm}$, this is of the form
$u^{f,\pm}(z) = (z,x^{\pm})$ for points $x^+ \in L_1$, $x^- \in \tau_L(L_0)$.
The condition for the two parts to match along $S^{f,+} \cap S^{f,-}$ is that
$x^- = x^+ = \tau_L(x^+)$. In particular $x^{\pm} \in L_0 \cap L_1$. Because of
the strong unique continuation property of horizontal sections, it follows that
the construction above yields all of $\MM^h$. By the same argument, any partial
horizontal section $U \rightarrow E|U$ with boundary in $Q$ can be extended to
some element of $\MM^h$. \qed

Let $j$ be some complex structure on $S$, standard over the ends. Take the same
$J^{(3)} \in \JJreg(M,\tau_L(L_0),L_1)$ as in the previous section, and choose
an additional $J^{(5)} \in \JJreg(L_0,L_1)$. By Lemma \ref{th:c} and
Proposition \ref{th:nonnegative-one}, one can find a $J^{(6)} \in
\JJ(E,\pi,Q,j,J^{(3)},J^{(5)})$ which is both horizontal and regular. Write
\[
c = C\Phi^{rel}_0(E,\pi,Q,J^{(6)}): C \longrightarrow C''
\]
for the chain map defined by this. In view of Lemmas \ref{th:nonnegative-two}
and \ref{th:new-intersection-points}\ref{item:q-map} one can write $c = \phi +
(c-\phi)$, where $\phi$ depends only on the sections in $\MM^h$ and is of order
$\{0\}$, and the remaining term is of order $[3\epsilon;\infty)$. Again
applying Lemma \ref{th:c}, one finds that $\phi$ is precisely the map $\gamma$
defined at the end of Section \ref{sec:setup}. Again, this is marginally better
than what is needed to apply Lemma \ref{th:low-energy}.

\subsection{The homotopy\label{sec:h}}

At this point it becomes necessary to change the notation slightly, in order to
avoid conflicts. Thus, the fibration used in Section \ref{sec:b} to define the
map $b$ will be denoted by $(E^b,\pi^b)$, its base by $S^b$, and its Lagrangian
boundary condition by $Q^b$; and correspondingly we write $(E^c,\pi^c)$, $S^c$,
$Q^c$ for the objects constructed in Section \ref{sec:c} to define $c$. Over
the unique negative end of $S^b$, $Q^b$ is modelled on $(\tau_L(L_0),L_1)$, and
the same holds for $Q^c$ over the positive end of $S^c$. As described in
Section \ref{sec:relative-invariants}, one can glue these ends together to
obtain a new exact Lefschetz fibration $(E^{bc},\pi^{bc})$ over a surface
$S^{bc}$, with a Lagrangian boundary condition $Q^{bc}$. The outcome is
represented schematically in Figure \ref{fig:h-map-one}.

\includefigure{h-map-one}{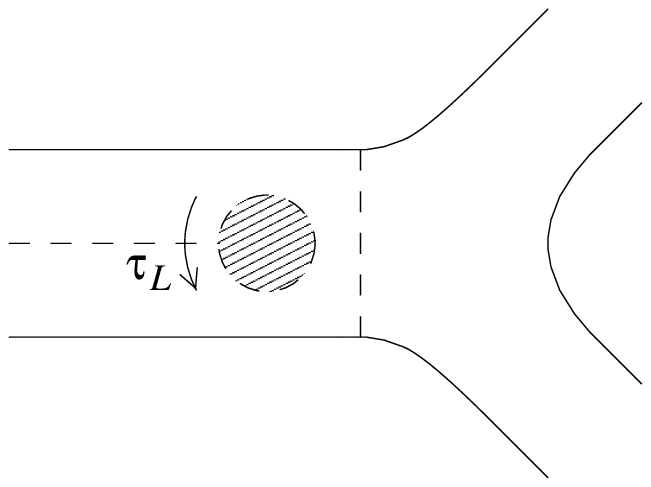}{ht}%
$(E^{bc},\pi^{bc})$ has nonnegative curvature because $(E^b,\pi^b)$ and
$(E^c,\pi^c)$ have that property. Moreover $\MM^h = \emptyset$, which means
that there are no horizontal sections $u: S^{bc} \rightarrow E^{bc}$ with
$u(\partial S^{bc}) \subset Q^{bc}$. In fact, assuming that such a $u$ exists,
one could reverse the gluing construction and obtain a horizontal section $u^b$
of $E^b$ with boundary in $Q^b$, as well as a similar section $u^c$ of $E^c$.
We have seen previously that such $u^b$ correspond to points in $\tau_L(L_0)
\cap L \cap L_1$, and $u^c$ to points in $L_0 \cap L_1$. In our case, the two
sections would have to match over the ends used in the gluing process, which
means that they would correspond to a point of $L_0 \cap L \cap L_1$; but that
is impossible by \ref{cond:transverse}. The same argument shows that for a
sufficiently large relatively compact subset $U \subset S^{bc}$, there are no
horizontal $w: U \rightarrow E^{bc}$ with $w(\partial S^{bc} \cap U) \subset
Q^{bc}$. Over the ends, $Q^{bc}$ is modelled on $(\tau_L(L_0),L)$, $(L,L_1)$
and $(L_0,L_1)$, respectively. For these pairs of submanifolds we have already
chosen almost complex structures $J^{(1)}$, $J^{(2)}$ and $J^{(5)}$,
respectively. Let $j$ be some complex structure on $S^{bc}$ which is standard
over the ends. Lemma \ref{th:nonnegative-one} ensures that one can choose a $J
\in \JJ(E^{bc},\pi^{bc},Q^{bc},j,J^{(1)},J^{(2)},J^{(5)})$ which is both
horizontal and regular.

\begin{lemma} \label{th:bc}
For any such $J$,
\[
C\Phi_0^{rel}(E^{bc},\pi^{bc},Q^{bc},J): C' \longrightarrow C''
\]
is homotopic to $c \circ b$ by a chain homotopy which is of order $(0;\infty)$.
\end{lemma}

\proof Suppose first that $j$ is induced from the complex structures on
$S^b,S^c$, and that $J$ is similarly constructed from $J^{(4)}$ and $J^{(6)}$;
this is automatically horizontal. Proposition \ref{th:gluing-ends} says that
for large values of the gluing parameter $\sigma$, $J$ is regular; and then we
have moreover $C\Phi_0(E^{bc},\pi^{bc},Q^{bc},J) = c \circ b$ by
\eqref{eq:chain-composition}. On the other hand, the observations made above
allow us to apply Lemma \ref{th:nonnegative-homotopy}, which shows that the
maps $C\Phi_0$ for any two choices of $j$ and $J$ are homotopic by a chain
homotopy of order $(0;\infty)$. \qed

\includefigure{h-map-two}{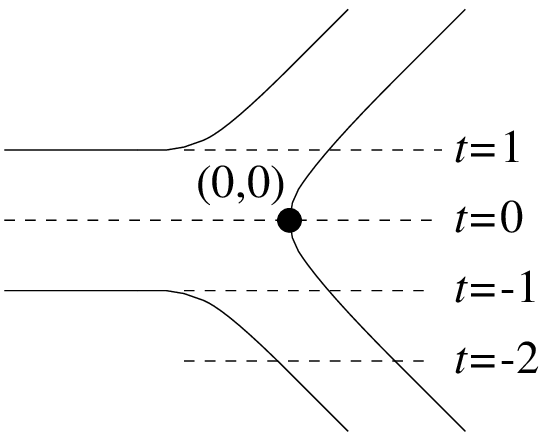}{hb} %
The proof that $c \circ b$ is chain homotopic to zero relies on an alternative
construction of the same exact Lefschetz fibration. Consider the surface $S^o$
from Figure \ref{fig:h-map-two}, embedded into $\R^2$ with coordinates $(s,t)$.
As in the previous section, we divide it into $S^{o,\pm} = S^o \cap \{t \in
\R^{\pm}\}$ and take the trivial fibrations $\pi^{o,\pm}: E^{o,\pm} = S^{o,\pm}
\times M \rightarrow S^{o,\pm}$. We equip $E^{o,\pm}$ with the forms
$\Omega^{o,\pm}$ pulled back from $\o$, and $E^{o,+}$ with the one-form
$\Theta^{o,+}$ pulled back from $\theta$; while on $E^{o,-}$ we take
$\Theta^{o,-} = \theta - d(\eta(t)K_{\tau_L})$, where $\eta$ is some function
with $\eta(t) = 1$ for $t \geq -1$, $\eta(t) = 0$ for $t \leq -2$. One now
identifies the fibres over $z \in \R^- \times \{0\} = S^{o,+} \cap S^{o,-}$ by
using $\tau_L: E^{o,+}_z \rightarrow E^{o,-}_z$, which yields an exact
symplectic fibration $(E^o,\pi^o)$ over $S^o$ (it is in fact trivial, but for
us it is convenient to think of it as being built up in this particular way).
Take the pieces of $\partial S^o$ labeled in Figure \ref{fig:h-map-three}, and
construct a Lagrangian boundary condition $Q^o$ for $(E^o,\pi^o)$ as the union
of $\partial_1S^{o,-} \times \tau_L(L_0), \,\, \partial_2S^{o,-} \times L
\subset E^{o,-}$ and $\partial_2S^{o,+} \times L, \,\, \partial_3S^{o,+} \times
L_1 \subset E^{o,+}$. The associated function $K_{Q^o}$ is equal to
$K_{\tau_L(L_0)} - \eta(t)K_{\tau_L}|\tau_L(L_0)$ over $\partial_1S^{o,-}$, to
$K_{L_1}$ over $\partial_3S^{o,+}$, and zero on the rest. $\kappa_{Q^o}$ is
equal to $-(K_{\tau_L}|L) \eta'(t) dt$ on $\partial_2S^{o,-}$ and vanishes
elsewhere; this makes sense because $K_{\tau_L}|L$ is constant, equal to $-2\pi
R(0)$ by \eqref{eq:kk-function}. Combining this with \ref{cond:wobbly} shows
that
\begin{equation} \label{eq:kappao}
\int_{\partial S^o} \kappa_{Q^o} = 2\pi R(0) \in (-\epsilon;0].
\end{equation}
$Q^o$ is modelled on $(\tau_L(L_0),L)$ and $(L,L_1)$ over the positive ends,
and on $(L_0,L_1)$ over the negative end; a suitable trivialization of
$(E^o,\pi^o)$ over that end can be defined in the same way as in Section
\ref{sec:c}. The next statement is a straightforward consequence of
\ref{cond:transverse}:

\begin{lemma} \label{th:o-sections}
Take the subset $U^o \subset S^o$ shaded in Figure \ref{fig:h-map-three}. Then
there are no horizontal sections $u: U^o \rightarrow E^o$ with $u(\partial S^o
\cap U^o) \subset Q^o$. \qed
\end{lemma}

\includefigure{h-map-three}{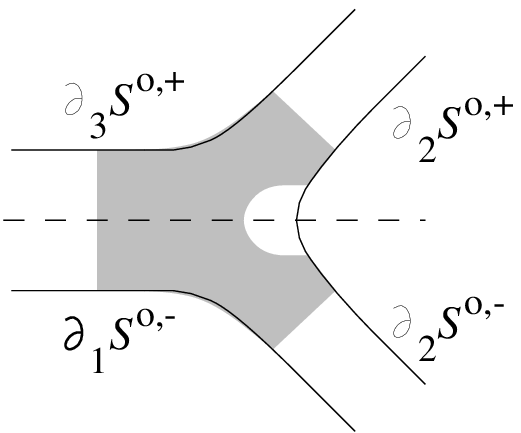}{ht} %
For the following step, we need the fibration $(E^p,\pi^p)$ over $S^p =
\bar{D}(1/2)$ from Section \ref{sec:c}, which was defined as pullback of a
standard fibration. The standard boundary condition from Section
\ref{sec:vanishing} pulls back to a Lagrangian boundary condition $Q^p$ for it.
As it stands $\kappa_{Q^p} = d^c(-\quarter|z|)$ is nonzero everywhere, but for
us it is better to modify it by some exact one-form, in such a way that it
becomes zero near $-1/2 \in \partial S^p$. This can be compensated by a change
of $K_{Q^p}$, so that the whole remains a Lagrangian boundary condition. It is
a consequence of Proposition \ref{th:vanishing} that the $\Phi_1$-invariant of
$(E^p,\pi^p,Q^p)$ vanishes, since this can be joined to $(E^L,\pi^L,Q^L)$ by a
smooth deformation of exact Lefschetz fibrations with Lagrangian boundary
conditions.

\begin{lemma} \label{th:p-sections}
Let $U^p \subset S^p$ be the complement of a sufficiently small neighbourhood
of $-1/2 \in \partial S^p$. Then there are no partial horizontal sections $w:
U^p \rightarrow E^p$ satisfying $w(\partial S^p \cap U^p) \subset Q^p$.
\end{lemma}

\proof From the standard fibration, $E^p$ inherits a smooth family of
Lagrangian spheres $\Sigma^p_z \subset E^p_z$, $z \neq 0$. These are carried
into each other by parallel transport along any path, and they degenerate to
the critical point $x_0 \in E^p_0$ as $z \rightarrow 0$. Now our boundary
condition is made up of the $\Sigma^p_z$ for $z \in \partial S^p$; therefore a
section $w$ with the properties stated above would satisfy $w(z) \in
\Sigma^p_z$ for all $z \in \partial S^p \cap U^p$, and by parallel transport
for all $z \in U^p \neq \{0\}$. In the limit this yields $w(0) = x_0$, but that
is impossible since $x_0$ is a critical point. \qed

One can identify the fibre of $E^o$ over $\zeta^o = (0,0) \in \partial S^o$
with the fibre of $E^p$ over $\zeta^p = -1/2 \in \partial S^p$ via
\begin{equation} \label{eq:identify}
(E^o)_{\zeta^o} \iso (E^{o,+})_{(0,0)} = M \xleftarrow{\phi^p} E^p_{1/2} \iso
(E^p)_{\zeta^p}
\end{equation}
where the last isomorphism is parallel transport along the upper semi-circle
$(1/2)e^{it}$, $0 \leq t \leq \pi$. By putting together the various
definitions, one sees that \eqref{eq:identify} takes $(Q^o)_{(0,0)}$ to
$(Q^p)_{-1/2}$. The fibrations $(E^o,\pi^o)$ and $(E^p,\pi^p)$ are flat near
the points $\zeta^o,\zeta^p$, and the one-forms $\kappa_{Q^o},\kappa_{Q^p}$
vanish near those points. This allows one to use the gluing construction
discussed in Section \ref{sec:simple-invariant} and again in Section
\ref{sec:relative-invariants} to produce an exact Lefschetz fibration $(E,\pi)$
over $S = S^o \#_{\zeta^o \sim \zeta^p} S^p$, together with a Lagrangian
boundary condition $Q$; see Figure \ref{fig:h-map-four}. Of course, this is
still modelled over the ends on the same Lagrangian submanifolds as $S^o$. In
what follows, we assume that the parameter $\rho$ in the gluing process has
been chosen sufficiently small.

\includefigure{h-map-four}{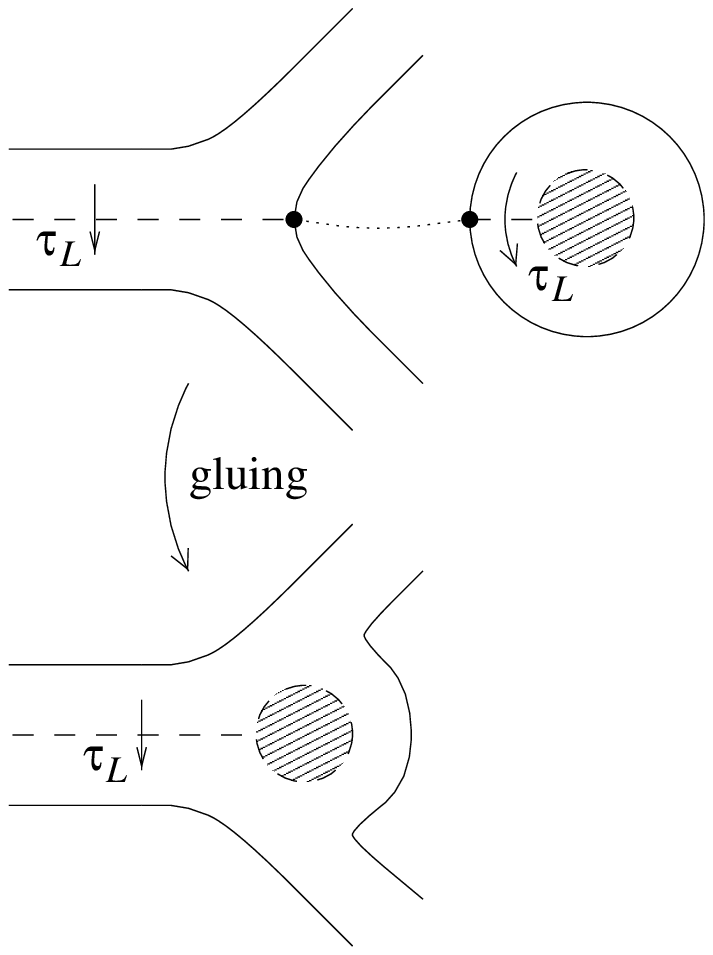}{hb} %
\begin{lemma} \label{th:order-homotopy}
There is a complex structure $j$ on $S$, standard over the ends, and a $J \in
\JJ^{reg,h}(E,\pi,Q,j,J^{(1)},J^{(2)},J^{(5)})$, such that
\[
C\Phi_0^{rel}(E,\pi,Q,J): C' \longrightarrow C''
\]
is homotopic to zero by a chain homotopy of order $[-2\pi R(0);\infty)$.
\end{lemma}

\proof Let $j^o$ be some complex structure on $S^o$, standard over the ends.
Using Lemmas \ref{th:o-sections} and \ref{th:nonnegative-one} one can find a
$J^o \in \JJ(E^o,\pi^o,j^o,J^{(1)},J^{(2)},J^{(5)})$ which is horizontal and
regular. In fact, regularity can be achieved even while prescribing what $J^o$
is outside $(\pi^o)^{-1}(U^o) \subset E^o$. This is useful because, for our
intended gluing argument, $J^o$ needs to be of a particular form close to the
fibre over $\zeta^o \notin U^o$: namely, in a local trivialization near that
point, it needs to be the product of $j^o$ and some previously fixed almost
complex structure on $M$. We now have evaluation maps, for $\tilde{x}_0 \in
\tau_L(L_0) \cap L$, $x_1 \in L \cap L_1$, $x \in \tau_L(L_0) \cap L_1$,
\begin{equation} \label{eq:ev-zetap}
ev_{\zeta^o}: \MM_{J^o}(\tilde{x}_0,x_1,x) \longrightarrow Q^o_{\zeta^o}.
\end{equation}
Take some complex structure $j^p$ on $S^p$. One can use Lemmas
\ref{th:horizontal-transversality} and \ref{th:p-sections} to find a $J^p \in
\JJ^{reg,h}(E^p,\pi^p,Q^p,j^p)$ with fixed behaviour outside
$(\pi^p)^{-1}(U^p)$; as before, since $\zeta^p \notin U^p$, one can use this to
make $J^p$ suitable for gluing. At the same time, Lemma
\ref{th:horizontal-ev-transversality} allows us to make the evaluation map
$ev_{\zeta^p}: \MM_{J^p} \rightarrow Q^p_{\zeta^p}$ transverse to any given
cycle. We take this cycle to be the disjoint union of \eqref{eq:ev-zetap} for
all $\tilde{x}_0,x_1,x$, identifying $Q^o_{\zeta^o}$ and $Q^p_{\zeta^p}$ via
\eqref{eq:identify}.

Let $j$ be the complex structure on $S$ glued together from $j^o,j^p$, and
similarly $J \in \JJ^h(E,\pi,Q,j,J^{(1)},J^{(2)},J^{(5)})$ the almost complex
structure obtained from $J^o$ and $J^p$. As discussed in Section
\ref{sec:relative-invariants}, $J$ will be regular if the gluing parameter has
been chosen sufficiently small. Moreover, the zero-dimensional spaces of
$(j,J)$-holomorphic sections can be described as fibre products of those on
both parts of the gluing, as in \eqref{eq:zero-gluing}. We know that
$\Phi_1(E^p,\pi^p,Q^p)$ is zero (as pointed out in Remark
\ref{th:vanishing}(i), this remains true even if we consider it as a cobordism
class), and using Lemma \ref{th:gluing-vanishing} one concludes that
$C\Phi_0^{rel}(E,\pi,Q,J)$ is homotopic to zero. Inspection of the proof of
that lemma will show that the homotopy constructed there is of order $[-2\pi
R(0);\infty)$. In fact, its coefficients are given by the number of points in
fibre products
\begin{equation} \label{eq:zero-product}
\MM_{J^o}(\tilde{x}_0,x_1,x) \times_{Q^o_{\zeta^o}} G
\end{equation}
where $(G,g)$ is some cycle bounding $(\MM_{J^p},ev_{\zeta^p})$. Using
\eqref{eq:kappao} and Lemma \ref{th:empty} one sees that
\eqref{eq:zero-product} is empty unless $a_{L_0,L_1}(x) \geq
a_{\tau_L(L_0),L}(\tilde{x}_0) + a_{L,L_1}(x_1) - 2\pi R(0)$, which provides
the desired estimate. \qed

To put together Lemmas \ref{th:bc} and \ref{th:order-homotopy}, one observes
that there are oriented diffeomorphisms
\[
\xymatrix{
 {E^{bc}} \ar[r]^{\Psi} \ar[d]_{\pi^{bc}} & {E} \ar[d]^{\pi} \\
 {S^{bc}} \ar[r]^{\psi} & {S}
}
\]
with the following properties: on the ends, $\psi$ and $\Psi$ relate the local
trivializations of the two fibrations. Next, $\Psi^*\Omega = \Omega^{bc}$, and
$\Psi(Q^{bc}) = Q$. Finally, $\psi$ is holomorphic near the unique critical
value, with respect to the complex structures defined there which are part of
the structure of exact Lefschetz fibrations of $(E^{bc},\pi^{bc})$ and
$(E,\pi)$; and the same holds for $\Psi$ near the unique critical point. This
is not difficult, since both fibrations contain a copy of $(E^p,\pi^p)$ and are
otherwise flat; comparing Figures \ref{fig:h-map-one} and \ref{fig:h-map-four}
shows how the bases should be identified in order for the monodromies to match.
It follows that one can take $j$ and $J$ in Lemma \ref{th:bc} to be the
pullback of almost complex structures from Lemma \ref{th:order-homotopy}, and
then $C\Phi_0^{rel}(E^{bc},\pi^{bc},Q^{bc},J)$ will be homotopic to zero by a
chain homotopy of order $[-2\pi R(0);\infty)$, with $-2\pi R(0) > 0$. On the
other hand, it is homotopic to $c \circ b$ by a homotopy of order $(0;\infty)$.
Taking the two together, one has a homotopy $h: c \circ b \htp 0$ of order
$(0;\infty)$, thus fulfilling the final requirement of Lemma
\ref{th:low-energy}.

\small

\end{document}